\documentclass[]{MyClass}

\usepackage[english]{babel}
\usepackage{url}
\usepackage[utf8x]{inputenc}
\usepackage{enumitem}
\usepackage{tikz}
\usetikzlibrary{shapes}
\usetikzlibrary{snakes}
\usetikzlibrary{cd}
\usepackage[lined,boxed,commentsnumbered]{algorithm2e}

\usepackage{rotating} 
\usepackage{shuffle}
\usepackage{cleveref}
\usepackage{multicol}
\usepackage{multirow}
\usepackage{ulem}
\newcommand{\pc}[1]{\rotatebox[origin=Bl]{90}{#1}}

\usepackage{todonotes}

\tikzset{style={anchor=base,baseline={([yshift=-1ex]current bounding box.center)}}}

\definecolor{mygreen}{RGB}{23,103,1}
\definecolor{mypurple}{RGB}{200,0,220}
\definecolor{Zgris}{rgb}{0.85,0.85,0.85}

\newcommand{\red}[1]{\textcolor{red}{#1}}
\newcommand{\blue}[1]{\textcolor{blue}{#1}}
\newcommand{\black}[1]{\textcolor{black}{#1}}

\newcommand{\blueul}[1]{\blue{\underline{\black{#1}}}}
\newcommand{\redul}[1]{\red{\underline{\black{#1}}}}

\newcommand{\eqdef}{\mbox{\,\raisebox{0.2ex}{\scriptsize\ensuremath{\mathrm:}}\ensuremath{=}\,}} 

\newcommand{\MM}{\mathbb{M}}
\newcommand{\NN}{\mathbb{N}}

\newcommand{\QQ}{\mathbb{Q}}
\newcommand{\RR}{\mathbb{R}}

\newcommand{\OO}{\mathbb{O}}
\newcommand{\PP}{\mathbb{P}}

\newcommand{\Sn}{\mathfrak{S}_n} 

\newcommand{\PW}{\mathbf{PW}} 

\newcommand{\SA}{\mathcal{A}} 
\newcommand{\SP}{\mathcal{P}} 
\newcommand{\ST}{\mathcal{T}} 
\newcommand{\SI}{\mathcal{I}} 
\newcommand{\Si}{\mathfrak{i}} 

\newcommand{\PForestred}{\mathfrak{F}_{\red{R}}} 
\newcommand{\PTreered}{\mathfrak{T}_{\red{R}}} 
\newcommand{\PParticularred}{\mathfrak{N}_{\red{R}}} 

\newcommand{\PForestredske}{\mathfrak{F}_{\red{R}ske}} 
\newcommand{\PTreeredske}{\mathfrak{T}_{\red{R}ske}} 
\newcommand{\PParticularredske}{\mathfrak{N}_{\red{R}ske}} 

\newcommand{\PForestredn}{\mathfrak{F}_{\red{R}n}} 
\newcommand{\PTreeredn}{\mathfrak{T}_{\red{R}n}} 
\newcommand{\PParticularredn}{\mathfrak{N}_{\red{R}n}} 

\newcommand{\PForestrednplusone}{\mathfrak{F}_{\red{R}n+1}} 
\newcommand{\PTreerednplusone}{\mathfrak{T}_{\red{R}n+1}} 
\newcommand{\PForestrednmp}{\mathfrak{F}_{\red{R}n-p}} 
\newcommand{\PForestredinfn}{\mathfrak{F}_{\red{R}\leq n}} 
\newcommand{\PTreeredinfn}{\mathfrak{T}_{\red{R}\leq n}} 

\newcommand{\PForestredone}{\mathfrak{F}_{\red{R}1}} 
\newcommand{\PTreeredone}{\mathfrak{T}_{\red{R}1}} 
\newcommand{\PParticularredone}{\mathfrak{N}_{\red{R}1}} 

\newcommand{\PForestredtwo}{\mathfrak{F}_{\red{R}2}} 
\newcommand{\PForestredthree}{\mathfrak{F}_{\red{R}3}} 

\newcommand{\PForestblue}{\mathfrak{F}_{\blue{B}}} 
\newcommand{\PTreeblue}{\mathfrak{T}_{\blue{B}}} 
\newcommand{\PParticularblue}{\mathfrak{N}_{\blue{B}}} 

\newcommand{\PForestblueske}{\mathfrak{F}_{\blue{B}ske}} 
\newcommand{\PTreeblueske}{\mathfrak{T}_{\blue{B}ske}} 
\newcommand{\PParticularblueske}{\mathfrak{N}_{\blue{B}ske}} 

\newcommand{\PForestbluen}{\mathfrak{F}_{\blue{B}n}} 
\newcommand{\PTreebluen}{\mathfrak{T}_{\blue{B}n}} 
\newcommand{\PParticularbluen}{\mathfrak{N}_{\blue{B}n}} 

\newcommand{\PForestbluenplusone}{\mathfrak{F}_{\blue{B}n+1}} 
\newcommand{\PTreebluenplusone}{\mathfrak{T}_{\blue{B}n+1}} 
\newcommand{\PForestbluenmp}{\mathfrak{F}_{\blue{B}n-p}} 
\newcommand{\PForestblueinfn}{\mathfrak{F}_{\blue{B}\leq n}} 
\newcommand{\PTreeblueinfn}{\mathfrak{T}_{\blue{B}\leq n}} 

\newcommand{\PForestblueone}{\mathfrak{F}_{\blue{B}1}} 
\newcommand{\PTreeblueone}{\mathfrak{T}_{\blue{B}1}} 
\newcommand{\PParticularblueone}{\mathfrak{N}_{\blue{B}1}} 

\newcommand{\PForestbluetwo}{\mathfrak{F}_{\blue{B}2}} 
\newcommand{\PForestbluethree}{\mathfrak{F}_{\blue{B}3}} 

\newcommand{\PForestrb}{\mathfrak{F}_{\red{R}\blue{B}}}
\newcommand{\PForestbr}{\mathfrak{F}_{\blue{B}\red{R}}}

\newcommand{\FQSym}{\mathbf{FQSym}{}}
\newcommand{\WQSym}{\mathbf{WQSym}{}}
\newcommand{\PQSym}{\mathbf{PQSym}{}}

\newcommand{\Id}{\operatorname{Id}}
\newcommand{\Imm}{\operatorname{Im}}
\newcommand{\Ker}{\operatorname{Ker}}
\newcommand{\Nodered}{\operatorname{Node_{\red{R}}}}
\newcommand{\Nodeblue}{\operatorname{Node_{\blue{B}}}}
\newcommand{\Noderb}{\operatorname{Node_{\red{R}\blue{B}}}}
\newcommand{\Nodebr}{\operatorname{Node_{\blue{B}\red{R}}}}
\newcommand{\Nodexy}{\operatorname{Node_{XY}}}
\newcommand{\Nodeyx}{\operatorname{Node_{YX}}}
\newcommand{\Primof}{\operatorname{Prim}}
\newcommand{\TPrimof}{\operatorname{TPrim}}
\newcommand{\Prim}{\operatorname{\mathbf{Prim}}}
\newcommand{\TPrim}{\operatorname{\mathbf{TPrim}}}

\newcommand{\std}{\operatorname{std}}
\newcommand{\pack}{\operatorname{pack}}
\newcommand{\gcdot}{\slash} 
\newcommand{\dcdot}{\backslash} 
\newcommand{\ins}{\red{\blacktriangleright}} 
\newcommand{\insl}{\blue{\blacktriangle}} 
\newcommand{\bijirr}{\mu} 
\newcommand{\Bijirr}{M_\mu} 
\newcommand{\bijt}{\sigma_\mu} 
\newcommand{\Bij}{\Sigma_\mu} 

\newcommand{\valshuffle}{\underline{\shuffle}} 
\newcommand{\valcshuffle}{\underline{\cshuffle}} 
\newcommand{\valprec}{\preceq} 
\newcommand{\valsucc}{\succeq} 
\newcommand{\valprecM}{\ll} 
\newcommand{\valsuccM}{\gg} 

\newcommand{\Frske}{\operatorname{F_{\red{R}ske}}}
\newcommand{\Trske}{\operatorname{T_{\red{R}ske}}}
\newcommand{\Fbske}{\operatorname{F_{\blue{B}ske}}}
\newcommand{\Tbske}{\operatorname{T_{\blue{B}ske}}}
\newcommand{\Fr}{\operatorname{F_{\red{R}}}}
\newcommand{\Tr}{\operatorname{T_{\red{R}}}}
\newcommand{\Fb}{\operatorname{F_{\blue{B}}}}
\newcommand{\Tb}{\operatorname{T_{\blue{B}}}}
\newcommand{\Frb}{\operatorname{F_{\red{R}\blue{B}}}}
\newcommand{\Trb}{\operatorname{T_{\red{R}\blue{B}}}}
\newcommand{\Fbr}{\operatorname{F_{\blue{B}\red{R}}}}
\newcommand{\Tbr}{\operatorname{T_{\blue{B}\red{R}}}}

\newcommand{\Frskestar}{\operatorname{F_{\red{R}ske}^{~*}}}
\newcommand{\Trskestar}{\operatorname{T_{\red{R}ske}^{~*}}}
\newcommand{\Fbskestar}{\operatorname{F_{\blue{B}ske}^{~*}}}
\newcommand{\Tbskestar}{\operatorname{T_{\blue{B}ske}^{~*}}}
\newcommand{\Frstar}{\operatorname{F_{\red{R}}^{~*}}}
\newcommand{\Trstar}{\operatorname{T_{\red{R}}^{~*}}}
\newcommand{\Fbstar}{\operatorname{F_{\blue{B}}^{~*}}}
\newcommand{\Tbstar}{\operatorname{T_{\blue{B}}^{~*}}}
\newcommand{\Frbstar}{\operatorname{F_{\red{R}\blue{B}}^{~*}}}
\newcommand{\Trbstar}{\operatorname{T_{\red{R}\blue{B}}^{~*}}}
\newcommand{\Fbrstar}{\operatorname{F_{\blue{B}\red{R}}^{~*}}}
\newcommand{\Tbrstar}{\operatorname{T_{\blue{B}\red{R}}^{~*}}}

\newcommand{\Frskeinv}{\operatorname{F_{\red{R}ske}^{~-1}}}
\newcommand{\Trskeinv}{\operatorname{T_{\red{R}ske}^{~-1}}}
\newcommand{\Fbskeinv}{\operatorname{F_{\blue{B}ske}^{~-1}}}
\newcommand{\Tbskeinv}{\operatorname{T_{\blue{B}ske}^{~-1}}}
\newcommand{\Frinv}{\operatorname{F_{\red{R}}^{~-1}}}
\newcommand{\Trinv}{\operatorname{T_{\red{R}}^{~-1}}}
\newcommand{\Fbinv}{\operatorname{F_{\blue{B}}^{~-1}}}
\newcommand{\Tbinv}{\operatorname{T_{\blue{B}}^{~-1}}}
\newcommand{\Frbinv}{\operatorname{F_{\red{R}\blue{B}}^{~-1}}}
\newcommand{\Trbinv}{\operatorname{T_{\red{R}\blue{B}}^{~-1}}}
\newcommand{\Fbrinv}{\operatorname{F_{\blue{B}\red{R}}^{~-1}}}

\newcommand{\tdeux}{\begin{picture}(7,7)(0,-1)
\put(3,0){\circle*{2}}
\put(3,0){\line(0,1){5}}
\put(3,5){\circle*{2}}
\end{picture}}
\newcommand{\ttroisdeux}{\begin{picture}(5,12)(-2,-1)
\put(0,0){\circle*{2}}
\put(0,0){\line(0,1){5}}
\put(0,5){\circle*{2}}
\put(0,5){\line(0,1){5}}
\put(0,10){\circle*{2}}
\end{picture}}

\newcommand{\ttroisun}{\begin{picture}(15,8)(-5,-1)
\put(3,0){\circle*{2}}
\put(-1,0){$\vee$}
\put(6,7){\circle*{2}}
\put(0,7){\circle*{2}}
\end{picture}}

\newcommand{\tdun}[1]{\begin{picture}(10,0)(-2,-1)
\put(0,0){\circle*{2}}
\put(3,-2){\tiny #1}
\end{picture}}
\newcommand{\tddeux}[2]{\begin{picture}(12,5)(0,-1)
\put(3,0){\circle*{2}}
\put(3,0){\line(0,1){5}}
\put(3,5){\circle*{2}}
\put(6,-2){\tiny #1}
\put(6,3){\tiny #2}
\end{picture}}

\crefname{coro}{Corollary}{Corollaries}
\crefname{defi}{Definition}{Definitions}
\crefname{lem}{Lemma}{Lemmas}
\crefname{algocf}{Algorithm}{Algorithms}
\crefname{ex}{Example}{Examples}
\crefname{rem}{Remark}{Remarks}

\newtheorem{theorem}{Theorem}
\newtheorem{lem}[theorem]{Lemma}

\newtheorem{coro}[theorem]{Corollary}
\newtheorem{prop}[theorem]{Proposition}

\theoremstyle{definition}
\newtheorem{defi}[theorem]{Definition}

\theoremstyle{remark}
\newtheorem{ex}[theorem]{Example}
\newtheorem{rem}[theorem]{Remark}
\newtheorem{nota}[theorem]{Notations}

\tikzstyle{alert} = [color=red, line width = 1.5]
\tikzstyle{bluealert} = [color=blue, line width =1.5]
\tikzstyle{big} = [line width = 1.5]
\tikzstyle{Point} = [fill, radius=0.08]
\tikzstyle{RedPoint} = [fill, radius=0.09, color = red]
\tikzstyle{Leaf} = [color = gray]

\tikzstyle{Red} = [color = red]
\tikzstyle{Blue} = [color = blue]
\tikzstyle{Green} = [color = mygreen]
\tikzstyle{Gray} = [color = gray]

\newcommand\Item[1][]{%
  \ifx\relax#1\relax  \item \else \item[#1] \fi
  \abovedisplayskip=0pt\abovedisplayshortskip=0pt~\vspace*{-\baselineskip}}

\title{Decompositions of packed words and self duality of Word
  Quasisymmetric Functions}

\author{Hugo MLODECKI \thanks{\href{mailto:mlodecki@lisn.fr}{mlodecki@lisn.fr}}}

\address{Université Paris-Saclay, CNRS, Laboratoire Interdisciplinaire des
  Sciences du Numérique, 91405, Orsay, France.}

\received{\today}

\abstract{By Foissy's work, the bidendriform structure of the Word
  Quasisymmetric Functions Hopf algebra (WQSym) implies that it is isomorphic to
  its dual. However, the only known explicit isomorphism due to Vargas does not
  respect the bidendriform structure. This structure is entirely determined by
  so-called totally primitive elements (elements such that the two
  half-coproducts vanish). In this paper, we construct two bases indexed by two
  new combinatorial families called red (dual side) and blue (primal side)
  biplane forests in bijection with packed words. In those bases, primitive
  elements are indexed by biplane trees and totally primitive elements by a
  certain subset of trees. We carefully combine red and blue forests to get
  bicolored forests. A simple recoloring of the edges allows us to obtain the
  first explicit bidendriform automorphism of WQSym.}

\resume{Grâce aux travaux de Foissy, on sait que l'algèbre de Hopf WQSym est
  isomorphe à sa duale car bidendriforme. Cependant, le seul isomorphisme
  explicite connu (dû à Vargas) ne respecte pas la structure
  bidendriforme. Cette structure est entièrement déterminée par les éléments
  totalement primitifs (annulés par les demi co-produits). Dans ce papier, nous
  construisons deux bases indexées par deux nouvelles familles combinatoire
  appelées forêts biplanes rouges (coté duale) et bleues (coté primale), en
  bijection avec les mots tassés. Dans ces bases, les éléments primitifs sont
  indexés par les arbres et les totalement primitifs par un certain
  sous-ensemble d'arbres. Nous combinons soigneusement les forêts rouges et
  bleues pour obtenir des forêts bicolores. Une simple recoloration des arêtes
  nous permet d'obtenir le premier automorphisme bidendriforme explicite de
  WQSym.}


\begin{document}


\maketitle

\section*{Introduction}

Combinatorial Hopf algebras are a common meeting point of different
communities. On one hand, the operad theory allows to construct them frequently
since the free algebra on one generator admits a Hopf structure. For many
operads, one can make the structure explicit using combinatorics, one of the
most basic example being the free dendriform algebra on one generator realized
as the Loday-Ronco hopf algebra of binary trees~\cite{LodRon_PBT}.

On the other hand, the theory of symmetric functions often proceeds through
non-commmutative lifting to better understand the identities. Hence, the
community introduced a series of larger and larger Hopf algebras over a large
variety of combinatorial structures. One of the first step was the introduction
of the dual pair of quasi-symmetric functions and non-commutative symmetric
functions~\cite{Gessel,NCSF1} to understand the inner product of characters
through the descent algebra~\cite{MalReu}. It leads to the discovery of the
Malvenuto-Reutenauer algebra $\FQSym$ of permutations. Another example was the
introduction by Poirier~\cite{Poirier_thesis} of an algebra of Young
tableaux. In~\cite{NCSF6}, it was realized that it can lead to a very simple
proof of the Littlewood-Richardson rule.  \smallskip

An early meeting between the symmetric function and the operad communities was
the discovery~\cite{HivNovThi_PBT} that the same procedure allows to construct
both the algebra of tableaux and the algebra of binary trees from the algebra of
permutations. One just has to enforce some simple relation (respectively plactic
and sylvester relation) in the variable of the polynomial realization.

Aside from the algebra of permutations, there is another non-commutative lifting
of quasi-symmetric functions. Indeed Hivert's action on polynomials whose
invariant are quasi-symmetric functions~\cite{NCSF6,Hivert_thesis} can be lifted
to words. Here, its non-commutative invariants spans the Hopf algebra~$\WQSym$
of packed words or, equivalently, surjections or even ordered set
partitions. This algebra has various applications in the theory of free Lie
algebras is closely related to the Solomon-Tits algebra and twisted descents,
the development of which was motivated by the geometry of Coxeter groups, the
study of Markov chains on hyperplane arrangements (see ~\cite{NovPatThi} and the
reference therein).

\bigskip

To better understand these algebra, one has to investigate their structure.  For
the binary tree algebra it was shown in~\cite{HivNovThi_PBT} that it is free as
an algebra and isomorphic to its dual. Though those properties are quite obvious
for the algebra~$\FQSym$ of permutation, the situation of~$\WQSym$ is much more
difficult. Its first study is due to Bergeron-Zabrocki~\cite{BerZab}. They
showed that it is free and co-free. However, it was only conjectured
in~\cite{NCSF6} that its primitive Lie algebra is free and that it is self
dual. It is only by a deep theorem of Foissy~\cite{Foissy_2007} that one can
show that the second one is too. In particular, until Vargas's
work~\cite{Vargas_thesis}, no concrete isomorphism was known.

%

Independently, Novelli-Thibon worked on parking functions which is a super-set
of packed words. They endowed the hopf algebra of parking functions $\PQSym$
with a bidendriform bialgebra structure~\cite{NovThi_PQSym_bidend}. Then they
describe $\WQSym$ as a sub-bidendriform bialgebra of
$\PQSym$~\cite{NovThi_WQSym_bidend}. Recall that a dendriform algebra is an
abstraction of a shuffle algebra where the product is split in two
half-products. If the coproduct is also split, and certain compatibilities hold,
one gets the notion of bidendriform bialgebra~\cite{Foissy_2007}.

Building on the work of Chapoton and Ronco~\cite{Ronco_2000,Chapoton_2002},
Foissy~\cite{Foissy_2007} showed that the structure of a bidendriform bialgebra
is very rigid. In particular, he defined a specific subspace called the space of
totally primitive elements, and showed that it characterizes the whole
structure. This does not only re-prove the freeness and co-freeness, as well as
the freeness of the primitive lie algebra, but also shows that the structure of
a bidendriform bialgebra depends only on its Hilbert series (the series of
dimensions of its homogeneous components). In particular, any such algebra is
isomorphic to its dual. However, Foissy's isomorphism is not fully explicit and
depends on a choice of a basis of the totally primitive elements. To this end,
one needs an explicit basis of the totally primitive elements. Foissy described
such a construction for $\FQSym$~\cite{Foissy_2011}. In this paper, we construct
a far reaching generalization for packed words and $\WQSym$ so that the basis
discribed in~\cite{Foissy_2011} is simply a restriction to permutations and
$\FQSym$ is a sub-bidendriform bialgebra of $\WQSym$. We provide two explicit
bases of totally primitive elements, for $\WQSym$ and its dual, using a
bijection with certain families of trees called biplane.

\bigskip

We begin with a background section presenting Foissy's two rigidity structure
theorems that prove, among other things, the self-duality of any bidendriform
bialgebra (\cref{brace,prim_tot}). We then define the notion of packed words as
well as the two specific bases ($\QQ$ and $\RR$) of~$\WQSym$ and its dual, which
will be the starting point of our combinatorial analysis.

\cref{sect2_pw_trees} is devoted to the combinatorial construction of biplane
forests (\cref{red_packed_forest,blue_packed_forest}) which are our first key
ingredient. They record a recursive decomposition of packed words according to
their global descents (\cref{gd_fact}) and positions of the maximum letter
(\cref{red_fact}) or the value of the last letter (\cref{blue_fact}). We show
that the cardinalities of some specific sets of biplane trees match the
dimensions of primitive and totally primitive elements
(\cref{thm:bij_primal,thm:bij_dual}).

In \cref{sect3_bases} we construct two new bases ($\OO$ and $\PP$ \cref{O,P}) of
$\WQSym$ and its dual which each contain as a subset a basis for the primitive
and totally primitve elements (see~\cref{thm:P,thm:O}). To do so we decompose
the space of totally primitive elements as a certain direct sum which matchs the
combinatorial decomposition of packed words
(\cref{tot_stable_phi,tot_stable_psi}).

Finally in \cref{sect4_bij} we make explicit how bases $\OO$ and $\PP$ are
sufficient to have an infinite number of bidendriform automorphism of
$\WQSym$. Then we give an explicit isomorphism based on an involution on packed
words. The definition of the bijection require a new kind of forest mixing red
and blue, namely bicolored-packed forests.



\section{Background}

\subsection{Cartier-Milnor-Moore theorems for Bidendriform bialgebras}

The goal of this section is to recall the elements of the definition of
bidendriform bialgebras which are useful for the comprehension of this paper. We
refer to~\cite{Foissy_2007} for the full list of axioms.

A bialgebra is a vector space over a field $K$, endowed with an unitary
associative product $\cdot$ and a counitary coassociative coproduct $\Delta$ satisfying a
compatibility relation called the Hopf relation
$\Delta(a\cdot b) = \Delta(a)\cdot\Delta(b)$.  In this paper all bialgebras are assumed to be graded
and connected (\textit{i.e.} the homogeneous component of degree $0$ is $K$). They
are therefore Hopf algebras, as the existence of the antipode is
implied.

A \textbf{dendriform algebra}
(see~\cite{Loday_dialgebras,LodRon_PBT,Ronco_2000,Ronco_2002}) 
$A$ is a $K$-vector space, endowed with two binary bilinear operations $\prec$,
$\succ$ satisfying the following axioms, for all $a, b, c\in A$:
\begin{align}
  \label{E1} (a \prec b) \prec c &= a \prec (b \prec c + b \succ c),\\
  \label{E2} (a \succ b) \prec c &= a \succ (b \prec c),\\
  \label{E3} (a \prec b + a \succ b) \succ c &= a \succ (b \succ c).
\end{align}
Adding together~\cref{E1,E2,E3} show that the product
$a \cdot b \eqdef a \prec b + a \succ b$ is associative. Adding a subspace of scalars,
this defines a unitary algebra structure on $K \oplus A$.  In this paper,
all the dendriform algebras are graded and have null $0$-degree component so
that the associated algebra is connected.

Dualizing, one gets a notion of \textbf{co-dendriform co-algebra}
(see~\cite{Foissy_2007}) which is a $K$-vector space with two binary
co-operations (\textit{i.e.}, linear maps $A \to A\otimes A$) denoted by
$\Delta_\prec$, $\Delta_\succ$ satisfying the dual axioms of~\cref{E1,E2,E3}:
\begin{align}
  \label{E4} (\Delta_\prec \otimes \Id) \circ \Delta_\prec(a) &= (\Id \otimes \Delta_\prec + \Id \otimes \Delta_\succ) \circ \Delta_\prec(a),\\
  \label{E5} (\Delta_\succ \otimes \Id) \circ \Delta_\prec(a) &= (\Id \otimes \Delta_\prec) \circ \Delta_\succ(a),\\
  \label{E6} (\Delta_\prec \otimes \Id + \Delta_\succ \otimes \Id) \circ \Delta_\succ(a) &= (\Id \otimes \Delta_\succ) \circ \Delta_\succ(a).
\end{align}
Adding together~\cref{E4,E5,E6} show that the reduced
coproduct~$\tilde{\Delta}(a) \eqdef \Delta_\prec(a) + \Delta_\succ(a)$ is co-associative. On
$K \oplus A$, setting
$\Delta(a) \eqdef 1 \otimes a + a \otimes 1 + \tilde{\Delta}(a)$ defines a co-associative and
co-unitary coproduct.

A \textbf{bidendriform bialgebra} is a $K$-vector space which is both a dendriform
algebra and a co-dendriform co-algebra satisfying a set of four relations
relating respectively $\prec$ and $\succ$ with $\Delta_\prec$, $\Delta_\succ$ (see~\cite{Foissy_2007} for
more details). In these equations, we use a kind of Einstein notation where
$\tilde{\Delta}(a) = a' \otimes a''$ and $\Delta_\alpha(b) = b'_\alpha \otimes b''_\alpha$ with $\alpha \in \{\prec, \succ\}$.
\begin{align}
  \label{E7} \Delta_\succ (a\succ b) &=a'b'_\succ\otimes a''\succ b''_\succ +b'_\succ\otimes a\succ b''_\succ+ab'_\succ \otimes b''_\succ+a'\otimes a''\succ b+a\otimes b,\\
  \label{E8} \Delta_\succ (a\prec b) &=a'b'_\succ\otimes a''\prec b''_\succ +b'_\succ\otimes a\prec b''_\succ+a'\otimes a''\prec b,\\
  \label{E9} \Delta_\prec (a\succ b) &=a'b'_\prec\otimes a''\succ b''_\prec +b'_\prec\otimes a\succ b''_\prec+ab'_\prec \otimes b''_\prec,\\
  \label{E10} \Delta_\prec (a\prec b)&=a'b'_\prec\otimes a''\prec b''_\prec +b'_\prec\otimes a\prec b''_\prec+a'b \otimes a''+b\otimes a.
\end{align}
Adding those four relations shows that $\cdot$ and $\Delta$ as defined above defines a
proper bi-algebra.

\bigskip

We recall here the relevant results of Foissy~\cite{Foissy_2007} on the
rigidity of bidendriform bialgebras based on the works of Chapoton and
Ronco~\cite{Ronco_2000,Chapoton_2002}.

Let $A$ be a bidendriform bialgebra. We
denote~$\Primof(A)\eqdef \Ker(\tilde{\Delta})$ the set of \textbf{primitive} elements of
$A$. We also denote by~$\SA(z)$ and~$\SP(z)$ the Hilbert series of $A$ and
$\Primof(A)$ defined as $\SA(z) \eqdef \sum_{n=1}^{+\infty}\dim(A_n)z^n$ and
$\SP(z) \eqdef \sum_{n=1}^{+\infty}\dim(\Primof(A_n))z^n$.  The present work is based on
two analogues of the Cartier-Milnor-Moore theorems~\cite{Foissy_2007} which we
present now.  The first one is extracted from the proof
of~{\cite[Proposition~6]{Foissy_2011}}:
\begin{prop}\label{left_prod}
  Let $A$ be a bidendriform bialgebra and let $p_1\dots p_n \in \Primof(A)$. Then
  the map
  \begin{equation}
    p_1 \otimes p_2 \otimes \ldots \otimes p_n
    \mapsto p_1 \prec (p_2 \prec (\ldots \prec p_n)\ldots).
  \end{equation}
  is an isomorphism of co-algebras from $T^+(\Primof(A))$ (the non trivial part
  of the tensor algebra with deconcatenation as
  coproduct) to $A$.  As a consequence, taking a basis $(p_i)_{i\in I}$ of
  $\Primof(A)$, the
  family~$(p_{w_1} \prec (p_{w_2} \prec (\dots \prec p_{w_n})\dots))_w$
  where $w=w_1\dots w_n$ is a non empty word on $I$ defines a basis of $A$.  This
  implies the equality of Hilbert series $\SA=\SP/(1-\SP)$.
\end{prop}

One can further analyze $\Primof(A)$ using the so-called \textbf{totally
  primitive} elements of $A$ defined as
$\TPrimof(A) \eqdef \Ker(\Delta_\prec) \cap \Ker(\Delta_\succ)$. The associated Hilbert serie is
defined as $\ST(z) \eqdef \sum_{n=1}^{+\infty}\dim(\TPrimof(A_n))z^n$.  Recall that a
brace algebra is a $K$-vector space $A$ together with an $n$-multilinear
operation denoted as $\langle\ldots\rangle$ for all $n ≥ 2$ which satisfies certain relations
(see \cite{Ronco_2000} for details).
\begin{theorem}[{\cite[Theorem~4~and~5]{Foissy_2011}}]\label{brace}
  Let $A$ be a bidendriform bialgebra. Then $\Primof(A)$ is freely generated as
  a brace algebra by $\TPrimof(A)$ with brackets given by
  \begin{multline*}
    \langle p_1, \dots, p_{n-1}; p_n\rangle \eqdef \sum_{i=0}^{n-1}\ (-1)^{n-1-i}\\
    (p_1 \prec (p_2 \prec (\cdots \prec p_i)\cdots)) \succ p_n \prec ((\cdots(p_{i+1} \succ p_{i+2}) \succ \cdots) \succ p_{n-1}).
  \end{multline*}
\end{theorem}

A basis of $\Primof(A)$ is described by ordered trees that are decorated with
elements of $\TPrimof(A)$ where $p_n$ is the root and $p_1, \ldots, p_{n-1}$ are the
children (see \cite{Ronco_2000,Chapoton_2002,Foissy_2011}). This is reflected on
their Hilbert series as~\cite[{Corollary~37}]{Foissy_2007}:
$\ST = \SA / (1+\SA)^2$ or equivalently $\SP = \ST (1 + \SA)$.

Using \cref{left_prod,brace} together with a dimension argument, one can show the
two following corollaries:
\begin{coro}[{\cite[Theorem~2]{Foissy_2011}}]\label{prim_tot}
  Let $A$ be a bidendriform bialgebra. Then $A$ is freely generated as a
  dendriform algebra by $\TPrimof(A)$.
\end{coro}

\begin{coro}[{\cite{Foissy_tree1,Foissy_tree2}}]\label{coro:prim_tot}
  A basis of $A$ is described by ordered forests of ordered trees that are
  decorated with a basis of $\TPrimof(A)$.
\end{coro}

On this basis, the product can be described using grafting (see Proposition~28
in~\cite{Foissy_tree2}) and the coproduct as the deconcatenation of forests that
are word of trees (see Theorem~35 equation 7.(c) in~\cite{Foissy_tree1}).

\subsection{Packed words}

The algebra $\WQSym$ is a Hopf algebra whose bases are indexed by ordered set
partitions or equivalently surjections or even packed words. In this paper, we
use the latter which we define now.
\medskip

In this paper we will deal with words over the alphabet of positive integers
$\NN_{>0}$. We start with basic notations: First, $\max(w)$ is the maximum
letter of the word $w$ with the convention that $\max(\epsilon)=0$. Then $|w|$ is the
length (or size) of the word $w$. The concatenation of the two words $u$ and $v$
is denoted as $u\cdot v$. The shift of a word $w$ of a value $i$ is denoted by
$w^{[i]}$. Once that said, $u\gcdot v \eqdef u^{[\max(v)]}\cdot v$ (resp.
$u\dcdot v \eqdef u\cdot v^{[\max(u)]}$) is the left-shifted (resp. right-shifted)
concatenation of the two words where all the letters of the left (resp. right)
word are shifted by the maximum of the right (resp. left) word:
$1121\gcdot 3112 = 44543112$ and $1121\dcdot 3112 = 11215334$. We also use the
notation $u|_{\leq i}$ (resp. $u|_{>i}$) for the subword containing all letters
smaller (resp. strictly greater) than a value $i$.

\begin{defi}\label{packed_word}
  A word over the alphabet~$\NN_{>0}$ is \textbf{packed} if all the letters from
  $1$ to its maximum $m$ appears at least once. By convention, the empty
  word~$\epsilon$ is packed.  For $n\in\NN$, we denote by $\PW_n$ the set of all packed
  words of length (also called size) $n$ and $\PW = \bigsqcup_{n\in\NN}\PW_n$ the
  set of all packed words.
\end{defi}

\begin{table}[!h]
  \[
    \begin{tabular}{|c|c|c|c|c|c|c|c|c|c||l|}
      \hline
      $n$ & 1 & 2 & 3 & 4 & 5 & 6 & 7 & 8 & 9 & OEIS \\
      \hline
      $\PW_n$ & 1 & 3 & 13 & 75 & 541 & 4 683 & 47 293 & 545 835 & 7 087 261 & A000670 \\
      \hline
    \end{tabular} 
  \]
  \caption{Number of packed words of size smaller than $9$.}
  \label{oeis_PW}
\end{table}

\begin{defi}\label{pack}
  The packed word~$u\eqdef \pack(w)$ associated with a word over the
  alphabet~$\NN_{>0}$ is obtained by the following process: if
  $b_1< b_2< \dots < b_r$ are the distinct letters occurring in $w$, then $u$ is the
  image of $w$ by the homomorphism $b_i\mapsto i$.
\end{defi}
A word $u$ is packed if and only if $\pack(u) = u$.
\begin{ex}
  The word $4152142$ is not packed because the letter $3$ does not appear while
  the maximum letter is $5 > 3$. Meanwhile $\pack(4152142) = 3142132$ is a packed
  word. Here are all packed words of size $1$, $2$ and $3$ in lexicographic
  order:
  \[
    1,\quad 11\ 12\ 21,\quad
    111\ 112\ 121\ 122\ 123\ 132\ 211\ 212\ 213\ 221\ 231\ 312\ 321
  \]
\end{ex}

The fonction $\pack(w)$ is the analogue of the \textbf{standardization}
$\std(w)$ that returns a permutation.

\begin{defi}
  The standardized word $\std(w)$ associated with a word over the
  alphabet~$\NN_{>0}$ is obtained by iteratively scanning $w$ from left to
  right, and labelling the occurrences of its smallest letter, then labelling
  the occurrences of the next one, and so on.
\end{defi}

\begin{ex}
  For example, $\std(4152142) = \std(3142132) = 5173264$.
\end{ex}
\bigskip

For the reader familiar with ordered set partitions, there is a classical
bijection between packed words and ordered set partitions. The one corresponding
to a packed word $w_1\cdot w_2\cdots w_n$ is obtained by placing the index $i$ into the
$w_i$-th block.

\begin{ex}\label{OSP_PW}
  The word $121$ is associated with $\{\{13\},\{2\}\}$ and the word $113223$
  with $\{\{12\},\{45\},\{36\}\}$.
\end{ex}
\bigskip

To depict some definitions or lemmas, we will use box diagrams with Cartesian
coordonates for packed words. On these diagrams, positions are from left to
right (as reading direction) and values are from bottom to top. These diagrams
will also be used to represent different decompositions with different
colors. Transparency will order the decompositions.

\begin{ex}
  Here we have three examples: the representation of the packed word
  $214313$. Then the word $3415251$ decomposed with red-factorization (see
  \cref{red_fact}). Finally the general case of the red-blue-factorization (see
  \cref{rb_br_fact}) where it can be seen clearly that the blue-factorization is
  done after the red-factorization thanks to transparency.
  $$\scalebox{0.5}{\input{figures/box_d/pw}} \qquad
  \scalebox{0.5}{\input{figures/box_d/pw_ins}} \qquad
  \scalebox{0.9}{\input{figures/box_d/ins_of_insl}}$$
\end{ex}

Global descent are defined in \cite{AguSot_2005} on permutations, here we
generalise the definition on packed words.

\begin{defi} \label{global_descent} A \textbf{global descent} of a packed word
  $w$ is a position $c$ such that all the letters before or at position $c$ are
  strictly greater than all letters after position $c$.
\end{defi}

\begin{ex}
  The global descents of $w = 54664312$ are the positions 5 and 6. Indeed, all
  letters of $54664$ are greater than the letters of $312$ and this is also true
  for $546643$ and $12$.
\end{ex}

\begin{defi}\label{irreducible}
  A packed word $w$ is \textbf{irreducible} if it is non empty and it has no
  global descent.
\end{defi}
\begin{table}[!h]
  \[
    \begin{tabular}{|c|c|c|c|c|c|c|c|c|c||l|}
      \hline
      $n$ & 1 & 2 & 3 & 4 & 5 & 6 & 7 & 8 & 9 & OEIS \\
      \hline
      $p_n$ & 1 & 2 & 8 & 48 & 368 & 3 376 & 35 824 & 430 512 & 5 773 936 & A095989 \\
      \hline
    \end{tabular} 
  \]
  \caption{Number of irreducible packed words of size smaller than $9$.}
  \label{oeis_irr}
\end{table}

\begin{ex}
  The word $w' = 21331$ is irreducible.
\end{ex}

\begin{lem} \label{gd_fact} Each word~$w$ admits a unique factorization as
  $w = w_1\gcdot w_2\gcdot\dots\gcdot w_k$ such that $w_i$ is irreducible for
  all $i$.
\end{lem}

\begin{ex}
  The global descent decomposition of $54664312$ is $21331\gcdot 1\gcdot 12$.
  The word $n\cdot n-1\cdot...\cdot 1$ has $1\gcdot 1\gcdot ...\gcdot 1$ as global descent
  decomposition.
  \begin{figure}[!h]
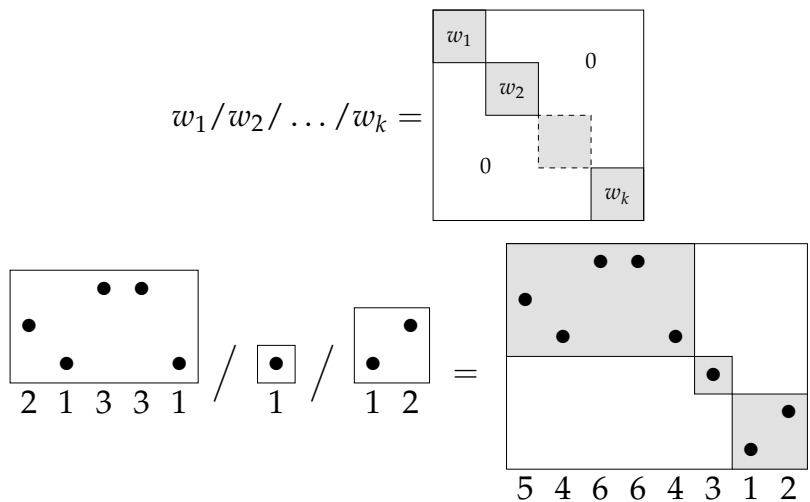

    $$w_1\gcdot w_2 \gcdot \dots \gcdot w_k = \scalebox{0.7}{\input{figures/box_d/gd_fact}}$$
    
    $$\scalebox{0.5}{\input{figures/box_d/21331}}~\Big{\gcdot}~\
    \scalebox{0.5}{\input{figures/box_d/1}}~\Big{\gcdot}~\
    \scalebox{0.5}{\input{figures/box_d/12}}~=~\
    \scalebox{0.5}{\input{figures/box_d/pw_gd_fact}}$$\vspace{-1em}
    \caption{Box diagrams: global descent decomposition.}
    \label{gd}
  \end{figure}
\end{ex}


\begin{defi}\label{def:shuffle}
  $u\shuffle v$ denotes the shuffle product of the two words. It is
  recursively defined by $u \shuffle \epsilon \eqdef \epsilon \shuffle u \eqdef u$ and
  \begin{equation}\label{eq:shuffle}
    ua \shuffle vb \eqdef (u \shuffle vb)\cdot a + (ua \shuffle v)\cdot b
  \end{equation}
  where $u$ and $v$ are words and $a$ and $b$ are letters. Analogously to the
  shifted concatenation, one can define the right shifted-shuffle
  $u\cshuffle v \eqdef u \shuffle v^{[\max(u)]}$ where all the letters of the
  right word $v$ are shifted by the maximum of the left word $u$.
\end{defi}

\begin{ex}
  $12\cshuffle \red{1}\red{1} =
  12\shuffle \red{3}\red{3} =
  12\red{3}\red{3} + 1\red{3}2\red{3} +
  1\red{3}\red{3}2 + \red{3}12\red{3} + \red{3}1\red{3}2 + \red{3}\red{3}12$.
  $$\scalebox{0.4}{\begin{tikzpicture}[baseline=1]
  \draw (0,0) rectangle ++(2,2);
  
  \draw[anchor = center] (0.5,0.5) node (w1) {\huge{$\bullet$}};
  \draw[anchor = center] (1.5,1.5) node (w1) {\huge{$\bullet$}};

  
  \draw[anchor = mid] (0.5,-0.7) node {\Huge{$1$}};
  \draw[anchor = mid] (1.5,-0.7) node {\Huge{$2$}};
\end{tikzpicture}} \;\scalebox{1.2}{$\cshuffle$}\;
  \scalebox{0.4}{\begin{tikzpicture}[auto]
  \draw (0,0) rectangle ++(2,1);
  
  \draw[anchor = center] (0.5,0.5) node (w1) {\huge{$\red{\bullet}$}};
  \draw[anchor = center] (1.5,0.5) node (w1) {\huge{$\red{\bullet}$}};

  
  \draw[anchor = mid] (0.5,-0.7) node {\Huge{$1$}};
  \draw[anchor = mid] (1.5,-0.7) node {\Huge{$1$}};
\end{tikzpicture}} \;\scalebox{1.2}{=}\;
  \scalebox{0.4}{\input{figures/box_d/1233}} \;\scalebox{1.2}{+}\;
  \scalebox{0.4}{\input{figures/box_d/1323}} \;\scalebox{1.2}{+}\;
  \scalebox{0.4}{\input{figures/box_d/1332}} \;\scalebox{1.2}{+}\;
  \scalebox{0.4}{\input{figures/box_d/3123}} \;\scalebox{1.2}{+}\;
  \scalebox{0.4}{\input{figures/box_d/3132}} \;\scalebox{1.2}{+}\;
  \scalebox{0.4}{\input{figures/box_d/3312}}$$
\end{ex}

\begin{defi}[{\cite[Example~5.4.(a)]{Loday_dialgebras}}]\label{def:half_shuffle}
  The recursive definition of the shuffle product \cref{eq:shuffle} contains two
  summands. The two half shuffle products on words $\prec$ and $\succ$ are defined respectively by:
  \begin{equation}\label{eq:half_shuffle}
    ua \prec vb \eqdef (u \shuffle vb)\cdot a
    \qquad\text{and}\qquad
    ua \succ vb \eqdef (ua \shuffle v)\cdot b.
  \end{equation}
\end{defi}

\begin{ex}
  $12 \prec \red{3}\red{3} = 1\red{3}\red{3}2 + \red{3}1\red{3}2 + \red{3}\red{3}12
  \qquad\text{and}\qquad 12 \succ \red{33} = 12\red{3}\red{3} + 1\red{3}2\red{3} +
  \red{3}12\red{3}$.
\end{ex}

\begin{defi}\label{def:valshuffle}
  $u\valshuffle v$ denote the dualisation of the deconcatenation
  using the function $\pack(w)$ of \cref{pack}.
  \begin{equation}\label{eq:valshuffle}
    u \valshuffle v \eqdef \sum_{\substack{u = \pack(u') \\ v = \pack(v')}}u'\cdot v'
  \end{equation}
  where $u, v$ and $u'\cdot v'$ are packed words.
  We also use the non-overlapping shuffle product on values by adding the
  constraint that letters of the two parts are distinct:
  \begin{equation}\label{eq:valcshuffle}
    u \valcshuffle v \eqdef \sum_{\substack{u = \pack(u') \\ v = \pack(v')\\\forall i,j, u'_i \neq v'_j}}u'\cdot v'.
  \end{equation}
\end{defi}

\begin{ex}
  $12\valshuffle \blue{11} = 12\blue{11} + 12\blue{22} + 12\blue{33} +
  13\blue{22} + 23\blue{11}$, \quad
  $12\valcshuffle \blue{11} = 12\blue{33} + 13\blue{22} + 23\blue{11}$.
  $$\scalebox{0.4}{\begin{tikzpicture}[baseline=1]
  \draw (0,0) rectangle ++(2,2);
  
  \draw[anchor = center] (0.5,0.5) node (w1) {\huge{$\bullet$}};
  \draw[anchor = center] (1.5,1.5) node (w1) {\huge{$\bullet$}};

  
  \draw[anchor = mid] (0.5,-0.7) node {\Huge{$1$}};
  \draw[anchor = mid] (1.5,-0.7) node {\Huge{$2$}};
\end{tikzpicture}} \;\scalebox{1.2}{$\valshuffle$}\;
  \scalebox{0.4}{\begin{tikzpicture}[auto]
  \draw (0,0) rectangle ++(2,1);
  
  \draw[anchor = center] (0.5,0.5) node (w1) {\huge{$\blue{\bullet}$}};
  \draw[anchor = center] (1.5,0.5) node (w1) {\huge{$\blue{\bullet}$}};

  
  \draw[anchor = mid] (0.5,-0.7) node {\Huge{$1$}};
  \draw[anchor = mid] (1.5,-0.7) node {\Huge{$1$}};
\end{tikzpicture}} \;\scalebox{1.2}{=}\;
  \scalebox{0.4}{\input{figures/box_d/1211}} \;\scalebox{1.2}{+}\;
  \scalebox{0.4}{\input{figures/box_d/1222}} \;\scalebox{1.2}{+}\;
  \scalebox{0.4}{\input{figures/box_d/1233b}} \;\scalebox{1.2}{+}\;
  \scalebox{0.4}{\input{figures/box_d/1322}} \;\scalebox{1.2}{+}\;
  \scalebox{0.4}{\input{figures/box_d/2311}}$$
  $$\scalebox{0.4}{\begin{tikzpicture}[baseline=1]
  \draw (0,0) rectangle ++(2,2);
  
  \draw[anchor = center] (0.5,0.5) node (w1) {\huge{$\bullet$}};
  \draw[anchor = center] (1.5,1.5) node (w1) {\huge{$\bullet$}};

  
  \draw[anchor = mid] (0.5,-0.7) node {\Huge{$1$}};
  \draw[anchor = mid] (1.5,-0.7) node {\Huge{$2$}};
\end{tikzpicture}} \;\scalebox{1.2}{$\valcshuffle$}\;
  \scalebox{0.4}{\begin{tikzpicture}[auto]
  \draw (0,0) rectangle ++(2,1);
  
  \draw[anchor = center] (0.5,0.5) node (w1) {\huge{$\blue{\bullet}$}};
  \draw[anchor = center] (1.5,0.5) node (w1) {\huge{$\blue{\bullet}$}};

  
  \draw[anchor = mid] (0.5,-0.7) node {\Huge{$1$}};
  \draw[anchor = mid] (1.5,-0.7) node {\Huge{$1$}};
\end{tikzpicture}} \;\scalebox{1.2}{=}\;
  \scalebox{0.4}{\input{figures/box_d/1233b}} \;\scalebox{1.2}{+}\;
  \scalebox{0.4}{\input{figures/box_d/1322}} \;\scalebox{1.2}{+}\;
  \scalebox{0.4}{\input{figures/box_d/2311}}$$
\end{ex}

\begin{rem}
  Using the classical bijection between ordered set partitions and packed words
  (see~\cref{OSP_PW}), the product $\valcshuffle$ is equivalent to the shifted
  shuffle on ordered set partitions defined in~\cite{BerZab}.
\end{rem}

\begin{defi}
  Analogously to the two half shuffle product of \cref{def:half_shuffle}, we
  split the two products $\valshuffle$ and $\valcshuffle$ in two parts.
  \begin{equation} \label{half_valshuffle}
    u \valprecM v \eqdef \sum_{\substack{u = \pack(u') \\ v = \pack(v')\\\max(u')>\max(v')}}u'\cdot v'
    \qquad\text{and}\qquad
    u \valsuccM v \eqdef \sum_{\substack{u = \pack(u') \\ v = \pack(v')\\\max(u')\leq\max(v')}}u'\cdot v'.
  \end{equation}
  \begin{equation} \label{half_valshuffle}
    u \valprec v \eqdef \sum_{\substack{u = \pack(u') \\ v = \pack(v')\\\forall i,j, u'_i \neq v'_j\\\max(u')>\max(v')}}u'\cdot v'
    \qquad\text{and}\qquad
    u \valsucc v \eqdef \sum_{\substack{u = \pack(u') \\ v = \pack(v')\\\forall i,j, u'_i \neq v'_j\\\max(u')<\max(v')}}u'\cdot v'.
  \end{equation}
\end{defi}

\begin{ex}
  ~\\
  $12 \valprecM \blue{11} = 12\blue{11} + 13\blue{22} + 23\blue{11} \qquad
  \text{and}\qquad 12 \valsuccM \blue{11} = 12\blue{22} + 12\blue{33}.$\\
  $12 \valprec \blue{11} = 13\blue{22} + 23\blue{11} \qquad\qquad\quad\,\,\,
  \text{and}\qquad 12 \valsucc \blue{11} = 12\blue{33}.$
\end{ex}

To sum up in a few words, in $u \prec v$ the last letter is comming from $u$ and
the rest is shuffled, in $u \valprec v$ the maximum value is comming from $u$
and the rest is shuffled.

\subsection{The Hopf algebra of word-quasisymmetric
  functions~\texorpdfstring{$\WQSym$}{WQSym}}

We are now in position to define the Hopf algebra of word-quasisymmetric
functions $\WQSym$. It was first defined as a Hopf algebra in
\cite{Hivert_thesis}. Novelli-Thibon proved later that $\WQSym$ and its dual are
bidendriform bialgebras~{\cite[{Theorems~2.5~and~2.6}]{NovThi_WQSym_bidend}}. Their
products and coproducts in the monomial basis $(\MM_w)_{w\in\PW}$ involve
overlapping-shuffle. However, to deal with the bidendriform structure, it will
be easier for us to chose, among the various bases known in the literature
\cite{Hivert_thesis,BerZab,NovThi_WQSym_bidend,Vargas_thesis} a basis where the shuffles
are non-overlapping. Therefore, for $\WQSym^*$, we take the basis denoted
$(\QQ_w)_{w\in\PW}$ of \cite[Equation~23]{BerZab} using the classical bijection
between ordered set partitions and packed words (see~\cref{OSP_PW}). For the
primal $\WQSym$, we define the dual basis denoted $(\RR_w)_{w\in\PW}$. In this
section we transfer the bidendriform structure on the bases $(\QQ_w)_{w\in\PW}$
and $(\RR_w)_{w\in\PW}$. \smallskip

Following Novelli-Thibon, we start from the basis $(\MM_w)_{w\in\PW}$ and compute
expressions of half product \cref{QQ_half_prod} and half coproduct
\cref{QQ_left_coprod,QQ_right_coprod} in the basis $(\QQ_w)_{w\in\PW}$. Then we
dualise these operations \cref{RR_half_prod,RR_left_coprod,RR_right_coprod} to
define the basis $(\RR_w)_{w\in\PW}$. The Hopf algebra product and reduced
coproduct are respectively recovered as the sum of the half products and half
coproducts. \bigskip

The monomial word-quasisymmetric function of a totally ordered alphabet
$\mathcal{A}$ associated to the packed word $u$ is the linear combination of
words defined by
$$\MM_u \eqdef \sum_{\substack{w \in \mathcal{A}^*,\\\pack(w) = u}}w.$$
It turns out that the concatenation of two such elements is a sum of
$(\MM_u)_{u\in\PW}$ so that $\WQSym \eqdef Vect(\MM_u \mid u \in \PW)$ is an
algebra. This can be refined to a bidendriform bialgebra structure. The
operations $\valprecM$, $\valsuccM$, $\Delta_{\valprecM}$ and
$\Delta_{\valsuccM}$ on
$$(\WQSym)_+ \eqdef Vect(\MM_u \mid u \in \PW_n, n \geq 1)$$ are defined in the following
way: for all $u = u_1\cdots u_n\in \PW_{n\geq1}$ and $v \in \PW_{m\geq1}$,
\begin{equation} \label{MM_half_prod}
  \MM_u \valprecM \MM_v \eqdef \sum_{w \in u \valprecM v} \MM_w,
  \qquad\text{and}\qquad
  \MM_u \valsuccM \MM_v \eqdef \sum_{w \in u \valsuccM v} \MM_w.
\end{equation}

\begin{align}
  \label{MM_left_coprod}  \Delta_{\valprecM}(\MM_u) &\eqdef \sum_{i = u_n}^{\max(u) - 1}
  \MM_{u|_{\leq i}} \otimes \MM_{\pack(u|_{> i})},\\
  \label{MM_right_coprod}  \Delta_{\valsuccM}(\MM_u) &\eqdef \sum_{i = 1}^{u_n - 1}
  \MM_{u|_{\leq i}} \otimes \MM_{\pack(u|_{> i})}.
\end{align}
\begin{ex}
  \begin{align*}
    \MM_{112} \valprecM \MM_{\blue{12}} &= \MM_{113\blue{12}} + \MM_{114\blue{23}} + \MM_{223\blue{12}} + \MM_{224\blue{13}} + \MM_{334\blue{12}},\\
    \MM_{112} \valsuccM \MM_{\blue{12}} &= \MM_{112\blue{12}} + \MM_{112\blue{13}} + \MM_{112\blue{23}} + \MM_{112\blue{34}} + \MM_{113\blue{23}} + \MM_{113\blue{24}} + \MM_{223\blue{13}} + \MM_{223\blue{14}},\\
    \Delta_{\valprecM}(\MM_{212536434}) &= \MM_{2123434} \otimes \MM_{\pack(56)} + \MM_{21253434} \otimes \MM_{\pack(6)},\\
                             &= \MM_{2123434} \otimes \MM_{12} + \MM_{21253434} \otimes \MM_{1},\\
    \Delta_{\valsuccM}(\MM_{212536434}) &= \MM_{1} \otimes \MM_{\pack(22536434)} + \MM_{212} \otimes \MM_{\pack(536434)} + \MM_{21233} \otimes \MM_{\pack(5644)}\\
                             &= \MM_{1} \otimes \MM_{11425323} + \MM_{212} \otimes \MM_{314212} + \MM_{21233} \otimes \MM_{2311}.
  \end{align*}
\end{ex}

\begin{theorem}{\cite[{Theorem~2.5}]{NovThi_WQSym_bidend}}
  $((\WQSym)_+, \valprecM, \valsuccM, \Delta_\valprecM, \Delta_\valsuccM)$ is a bidendriform bialgebra.
\end{theorem}

As we said, we want a basis without overlapping-shuffle for the
product. Following~\cite{BerZab}, we first define a partial order on packed
words then we define the new basis.

\begin{defi}{\cite{BerZab}}\label{BerZab_order}
  We say that the packed word $u$ is smaller than $v$ for the relation $\leq_*$ if
  $u$ and $v$ have the same standardization and if $u_i = u_j$ implies
  $v_i = v_j$ for all $i$ and $j$.
  $$ u \leq_* v \iff \std(u) = \std(v) \text{ and } (u_i = u_j \implies v_i = v_j).$$

\end{defi}
\vspace{-5mm}
\begin{figure}[!h]
  $$\scalebox{1}{\input{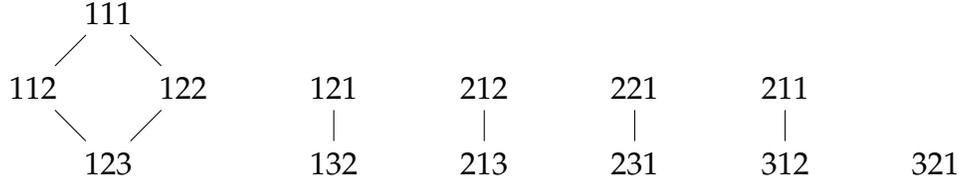}}$$
  \caption{The Hasse diagram of $(\PW_3, \leq_*)$.}
  \label{leq_star}
\end{figure}

We give two immediate lemmas on this order that are useful.

\begin{lem}\label{max-include}
  For $u \leq_* v$, let $m_u$ (resp. $m_v$) be the set of positions of occurrences
  of the maximum value letters in $u$ (resp. $v$). Then $m_u$ is included in
  $m_v$ and all positions in $m_v$ that are not in $m_u$ are smaller to the
  minimum of $m_u$.
\end{lem}

\begin{proof}
  It is immediate with the definition of $\leq_*$.
\end{proof}

\begin{lem}\label{order-restricted}
  For $u \leq_* v$, let $i$ and $i'$ such that $u|_{\leq i}$ and $v|_{\leq i'}$ are of
  the same size then $u|_{\leq i} \leq_* v|_{\leq i'}$ and $u|_{> i} \leq_* v|_{> i'}$.
\end{lem}

\begin{proof}
  It is immediate with the definition of $\leq_*$.
\end{proof}

Now we can recall ({\cite[{Equation~23}]{BerZab}}) the definition of the basis
$(\QQ_w)_{w\in\PW}$
\begin{equation}\label{QQ_MM}
  \QQ_u \eqdef \sum_{u\leq_*v}\MM_v.
\end{equation}

\begin{ex}
  \begin{align*}
    \QQ_{123} &= \MM_{123} + \MM_{122} + \MM_{112} + \MM_{111}  &\QQ_{43132} &= \MM_{43132} + \MM_{32121} \\
    \QQ_{412234} &= \MM_{412234} + \MM_{312223} + \MM_{311123} + \MM_{211112} &\QQ_{2131} &= \MM_{2131} + \MM_{2121} \\
  \end{align*}
\end{ex}

It is proved in {\cite[{Theorem~17}]{BerZab}} that the product in basis
$(\QQ_w)_{w\in\PW}$ is
$$\QQ_u\QQ_v = \sum_{w\in u \valcshuffle v}\QQ_w.$$
Thanks to \cref{max-include}, we have the two expressions for the two half
products.
\begin{equation} \label{QQ_half_prod}
  \QQ_u \valprec \QQ_v \eqdef \sum_{w \in u \valprec v} \QQ_w,
  \qquad\text{and}\qquad
  \QQ_u \valsucc \QQ_v \eqdef \sum_{w \in u \valsucc v} \QQ_w.
\end{equation}
For the coproduct, we start with the definition of the coproduct in basis $\MM$.
\begin{align*}
  \Delta(\QQ_u) &= \sum_{v\geq_*u}\Delta(\MM_v)\\
           &= \sum_{v\geq_*u} \left(\sum_{i=0}^{\max(v)}\MM_{v|_{\leq i}} \otimes \MM_{v|_{> i}}\right)\\
           &= \sum_{i=0}^{\max(u)} \left(\sum_{v \geq_*u|_{\leq i}}\MM_v \otimes \sum_{v' \geq_*u|_{> i}}\MM_{v'}\right) \text{\qquad (by \cref{order-restricted})}\\
           &= \sum_{i=0}^{\max(u)} \QQ_{u|_{\leq i}} \otimes \QQ_{u|_{> i}}.
\end{align*}

Then, with \cref{max-include} we have the two expressions for the two half
coproducts.

\begin{align}
  \label{QQ_left_coprod}  \Delta_{\valprec}(\QQ_u) &\eqdef \sum_{i = u_n}^{\max(u) - 1}
  \QQ_{u|_{\leq i}} \otimes \QQ_{\pack(u|_{> i})},\\
  \label{QQ_right_coprod}  \Delta_{\valsucc}(\QQ_u) &\eqdef \sum_{i = 1}^{u_n - 1}
  \QQ_{u|_{\leq i}} \otimes \QQ_{\pack(u|_{> i})}.
\end{align}
\begin{ex}
  \begin{align*}
    \QQ_{1312} \valprec \QQ_{\blue{12}} &= \QQ_{1512\blue{34}} + \QQ_{1513\blue{24}} + \QQ_{1514\blue{23}} + \QQ_{2523\blue{14}} + \QQ_{2524\blue{13}} + \QQ_{3534\blue{12}},\\
    \QQ_{1312} \valsucc \QQ_{\blue{12}} &= \QQ_{1312\blue{45}} + \QQ_{1412\blue{35}} + \QQ_{1413\blue{25}} + \QQ_{2423\blue{15}},\\
    \Delta_{\valprec}(\QQ_{212536434}) &= \QQ_{2123434} \otimes \QQ_{\pack(56)} + \QQ_{21253434} \otimes \QQ_{\pack(6)},\\
                             &= \QQ_{2123434} \otimes \QQ_{12} + \QQ_{21253434} \otimes \QQ_{1},\\
    \Delta_{\valsucc}(\QQ_{212536434}) &= \QQ_{1} \otimes \QQ_{\pack(22536434)} + \QQ_{212} \otimes \QQ_{\pack(536434)} + \QQ_{21233} \otimes \QQ_{\pack(5644)}\\
                             &= \QQ_{1} \otimes \QQ_{11425323} + \QQ_{212} \otimes \QQ_{314212} + \QQ_{21233} \otimes \QQ_{2311}.
  \end{align*}
\end{ex}




Finally we define $\prec$, $\succ$, $\Delta_\prec$ and $\Delta_\succ$ on
$(\WQSym^*)_+ \eqdef Vect(\RR_u \mid u \in \PW_n, n \geq 1)$ by dualizing half products and
half coproducts of the basis $(\QQ_w)_{w\in\PW}$ in the following way: for all
$u = u_1\cdots u_n\in \PW_{n\geq1}$ and $v \in \PW_{m\geq1}$,
\begin{equation}\label{RR_half_prod}
  \RR_u \prec \RR_v \eqdef \sum_{w \in u \prec v^{[\max(u)]}} \RR_w, 
  \qquad\text{and}\qquad
  \RR_u \succ \RR_v \eqdef \sum_{w \in u \succ v^{[\max(u)]}} \RR_w.
\end{equation}
\begin{align}
  \label{RR_left_coprod}  \Delta_\prec(\RR_u) &\eqdef \sum_{\substack{i=k\\\{u_1, \dots, u_i\} \cap \{u_{i+1}, \dots, u_n\} = \emptyset\\u_k = \max(u)}}^{n-1}
                                    \RR_{\pack(u_1\cdots u_i)} \otimes \RR_{\pack(u_{i+1}\cdots u_n)},\\
  \label{RR_right_coprod}  \Delta_\succ(\RR_u) &\eqdef \sum_{\substack{i=1\\\{u_1, \dots, u_i\} \cap \{u_{i+1}, \dots, u_n\} = \emptyset\\u_k = \max(u)}}^{k-1}
               \RR_{\pack(u_1\cdots u_i)} \otimes \RR_{\pack(u_{i+1}\cdots u_n)}.
\end{align}
\begin{ex}
  \begin{align*}
    \RR_{211} \prec \RR_{\red{12}} &= \RR_{21\red{3}\red{4}1} + \RR_{2\red{3}1\red{4}1} + \RR_{2\red{3}\red{4}11} + 
                                 \RR_{\red{3}21\red{4}1} + \RR_{\red{3}2\red{4}11} + \RR_{\red{3}\red{4}211},\\
    \RR_{221} \succ \RR_{\red{12}} &= \RR_{211\red{3}\red{4}} + \RR_{21\red{3}1\red{4}} + \RR_{2\red{3}11\red{4}}
                                 + \RR_{\red{3}211\red{4}},\\
    \Delta_\prec(\RR_{2125334}) &= \RR_{2123} \otimes \RR_{112} + \RR_{212433} \otimes \RR_{1},\\
    \Delta_\succ(\RR_{2125334}) &= \RR_{212} \otimes \RR_{3112}. 
  \end{align*}
\end{ex}

\begin{theorem}
  $((\WQSym)_+, \valprec, \valsucc, \Delta_{\valprec}, \Delta_{\valsucc})$ and
  $((\WQSym^*)_+, \prec, \succ, \Delta_\prec, \Delta_\succ)$ are two dual bidendriform bialgebras.
\end{theorem}

From now on $\Primof(\WQSym)$ and $\TPrimof(\WQSym)$ are respectively abbreviated to
$\Prim$ and $\TPrim$. Moreover, we denote homogeneous components using indices
and dualization using a $*$ in exponent as in $\Prim^*_n$. We give the first
values of the dimensions~$a_n\eqdef\dim(\WQSym_n)$, $p_n\eqdef\dim(\Prim_n)$ and
$t_n\eqdef\dim(\TPrim_n)$:
\begin{table}[!h]
  \[
    \begin{tabular}{|c|c|c|c|c|c|c|c|c|c||l|}
      \hline
      $n$ & 1 & 2 & 3 & 4 & 5 & 6 & 7 & 8 & 9 & OEIS \\
      \hline
      $a_n$ & 1 & 3 & 13 & 75 & 541 & 4 683 & 47 293 & 545 835 & 7 087 261 & A000670 \\
      \hline
      $p_n$ & 1 & 2 & 8 & 48 & 368 & 3 376 & 35 824 & 430 512 & 5 773 936 & A095989 \\
      \hline
      $t_n$ & 1 & 1 & 4 & 28 & 240 & 2 384 & 26 832 & 337 168 & 4 680 272 &  \\
      \hline
    \end{tabular} 
  \]
  \caption{Dimensions of homogeneous components for $\WQSym$, $\Prim$ and $\TPrim$.}
  \label{oeis_WQSym}
\end{table}

Though the numbers $(t_n)_n$ are easy to obtain thanks to the relation of
\cref{brace}: ${\ST = \SA / (1+\SA)^2}$, no combinatorial interpretation
existed. The first results of this paper are two different subsets of packed
words that are counted by these dimensions (red-irreducible and blue-irreducible).


\section{Decorated forests}\label{sect2_pw_trees}

In this paper we will generalize twice the construction of \cite{Foissy_2011},
one for $\WQSym$ and one for its dual. This section is devoted to the
combinatorial ingredient, that is a notion of biplane forests suitable for
indexing the various bases of primitive elements. Each time, we start by
decomposing packed words through global descents and removal of specific
letters. We then perform those decompositions recursively, encoding the result
in a forests. We hence obtain so called biplane forests, which are in bijection
with packed words. Later, the recursive structure of forests will be understood
as a chaining of brace and dendriform operations generating some elements of
$\WQSym$ or its dual. This will allow us to construct two bases of respectively
$\TPrim$ and its dual by characterizing a subfamily of biplane trees.  \bigskip

From now on, we associate the color \blue{blue} to the primal ($\WQSym$) and the
color \red{red} to the dual ($\WQSym^*$). We start by explaining the
construction on $\WQSym^*$ (\red{red}) then we dualize the construction to the
primal $\WQSym$ (\blue{blue}).

\subsection{Dual (\red{Red})}\label{Primal_Red}

For the red side $\WQSym^*$, the decomposition of packed words is made through
global descents and removal of maximum values. One step of this decomposition is
called the red-factorization.

\subsubsection{Decomposition of packed words through maximums}

In this section, we define two combinatorial operations on packed words
($\phi_I$ and $\ins$) and the red-factorization that uses them. The unary operation
$\phi_I$ inserts new maximums in a word in positions $I$. A word that cannot be
factorized $u \ins v$ in a non trivial way is called
red-irreducible. Red-irreducible words will index our basis of $\TPrim^*$.

\begin{defi}\label{phi}
  Fix $n \in \NN$ and $w \in \PW_n$. We write $m' \eqdef \max(w) + 1$. For any
  $p>0$ and any subset~$I \subseteq [1,\dots, n+p]$ of cardinality~$p$, we define
  $\phi_I(w) \eqdef u_1\dots u_{n+p}$ as the packed word of length $n+p$ obtained
  by inserting $p$ occurrences of the letter $m'$ in $w$ so that they end up in
  positions $i \in I$. In other words $u_i = m'$ if $i\in I$ and $w$ is obtained
  from $\phi_I(w)$ by removing all occurrences of $m'$. Notice that $\phi_I(w)$ is
  only defined if $n + p \geq i_p$.
\end{defi}

\begin{ex}
  $\phi_{2,4,7}(1232) = 1424324$ and $\phi_{1,2,3}(\epsilon)=111$.
\end{ex}

\begin{nota}
  For the rest of this paper, $I = [i_1, \dots, i_p]$ will always denote a
  non-empty ($p > 0$) list of increasing non-zero integers. For any integer $k$,
  $I' = I + k$ denote the list $I' = [i_1+k, \dots, i_p+k]$. Let $\PW_n^I$
  denote the set of packed words of size $n$ whose maximums are in positions
  $i\in I$. This way $\phi_I(w)\in \PW_{n+p}^I$ for any $w\in \PW_{n}$.
\end{nota}


\begin{lem}\label{phi_bij}
  Let $n \in \NN$ and $p>0$, for any
  $I = [i_1, \dots, i_p] \subseteq [1,\dots, n+p]$ of size $p$, $\phi_I$ is a bijection
  from~$\PW_{n}$ to the~$\PW_{n+p}^I$.
  
  Moreover, for any $W\in\PW_\ell$ where~$\ell>0$ there exists a unique pair
  $(I,w)$ where $I\subseteq [1\dots \ell]$ and $w$ is packed, such that $W = \phi_I(w)$.
\end{lem}
The box diagram that pictures this lemma is $W = \input{figures/box_d/phi}$.
\begin{proof}
  Let $W\in\PW_\ell$ with~$\ell>0$ and $m$ the value of the maximum letter of
  $W$. Let $I = [i_1, \dots, i_p] \subseteq [1,\dots, \ell]$ be the list of the positions
  of $m$ in $W$ and let $w$ be the word obtain by removing all occurrences of
  $m$ in $W$, then $W = \phi_I(w)$. If $\phi_I(u) = \phi_J(v)$ then positions of maximum
  values are the same so $I = J$ and words obtain by removing these maximum
  values are also the same so $u = v$.
\end{proof}

\begin{defi}\label{ins}
  Let $u, v \in \PW$ with $v \neq \epsilon$. By~\cref{phi_bij}, there is a unique pair
  $(I, v')$ such that $v=\phi_I(v')$. Let $I' = I + |u|$, we define
  $u \ins v \eqdef \phi_{I'}(u\gcdot v')$. In other words, we remove the maximum
  letter of the right word, perform a left shifted concatenation and reinsert
  the removed letters as new maximums.
\end{defi}

\begin{ex}
  $2123\ins \red{3}\blue{22}\red{3}\blue{12} ~=~ 2123 \ins \phi_{1, 4}(\blue{2212}) ~=~
  \phi_{1+4, 4+4}(4345\blue{2212}) ~=~ 4345\red{6}\blue{22}\red{6}\blue{12}$.
\end{ex}
\begin{figure}[!h]
  $$u\ins v = \scalebox{1}{\input{figures/box_d/ins}} \qquad
  \scalebox{0.4}{\input{figures/box_d/2123}}~\ins~\scalebox{0.4}{\input{figures/box_d/322312}}~=~\scalebox{0.4}{\input{figures/box_d/pw_ins_ex}}$$\vspace{-1em}
  \caption{Box digrams: the operation $\ins$}
  \label{boxd_ins}
\end{figure}

\begin{lem}\label{red_fact}
  Let~$w$ be an irreducible packed word. There exists a unique factorization of
  the form $w=u \ins v$ which maximizes the size of $u$. In this factorization,
  let $v'$ and $I$ be such that $v = \phi_I(v')$, then
  \begin{itemize}
  \item either $v' = \epsilon$ and $I = [1, \dots, p]$ for some $p$,
  \item or the global descent decomposition $v'=v_1\gcdot\dots\gcdot v_r$ of
    $v'$ satisfies the inequalities $1 \leq i_1 \leq |v_1|$, and
    $1 \leq (|I| + |v'|) + 1 - i_p \leq |v_r|$ with $I = [i_1, \dots, i_p]$.
  \end{itemize}
  We call it the \textbf{red-factorization} of a word.
\end{lem}

\begin{ex}
  Here is a first detailed example of a red-factorization of an irreducible
  packed word:

  Consider the irreducible packed word $w = 543462161$.
  \begin{itemize}
  \item The first step is to remove all the occurrences of the maximal value but
    keep in memory the positions in the initial word. We get $w' = 5434\_21\_1$
    which is a packed word, but is not irreducible.
  \item The second step is to decompose the new word $w'$ in irreducible
    factors $w' = 1\gcdot212\gcdot\_1\gcdot1\_1$. We still keep in memory the
    positions of the removed value. (when we have the choice, we cut to the left
    of the removed value.)
  \item We can distinguish two groups of factors, those strictly before the
    first maximum withdrawn and the others
    $w' = 1\gcdot212\quad\gcdot\quad\_1\gcdot1\_1$.
  \item Finally, by numbering the positions of the maximum removed value in the
    right factor (positions $1$ and $4$), we get the following decomposition of
    $w$ (see \cref{phi} for $\phi$ and \cref{ins} for $\ins$): \\
    $w = 543462161 = (1 \gcdot 212) \ins \phi_{1,4}(1 \gcdot 11) = (3212) \ins
    \phi_{1,4}(211) = 3212 \ins 32131$.
  \end{itemize}
\end{ex}

\begin{ex} Here are some other red-factorizations:
  \begin{alignat*}{2}
    21331 &= 1\ins\phi_{2,3}(11)  = 1\ins 1221 &\qquad 1231 &= \epsilon \ins \phi_3(121) = \epsilon \ins 1231\\
    1233 &= 12 \ins \phi_{1,2}(\epsilon) = 12 \ins 11 &\qquad 111 &= \epsilon \ins \phi_{1,2,3}(\epsilon) = \epsilon \ins 111\\
    56434126 &= 1 \ins \phi_{1,7}(212 \gcdot 12) \rlap{$ = 1 \ins \phi_{1,7}(43412) = 1 \ins 5434125$}
  \end{alignat*}
\end{ex}

\begin{proof} Let $w$ be irreducible and let $(I, w')$ be the unique pair such
  that $w=\phi_I(w')$ according to~\cref{phi_bij}. By~\cref{gd_fact}, we write
  $w' = w_1'\gcdot w_2' \gcdot\dots\gcdot w_k'$, the unique decomposition into
  irreducibles. Let $\ell$ be such that $w_{\ell}'$ is the last factor which is
  entirely before the first removed maximum, it is the only choice to maximize
  the size of $u$. Then with $r=k-\ell$ we can rewrite $w'$ as
  $(u_1\gcdot \dots\gcdot u_\ell)\gcdot(v_1\gcdot \dots\gcdot v_{r})$. Now we get
  $I'$ by subtracting $|u_1\gcdot \dots\gcdot u_\ell|$ to all parts of $I$
  ($I' = I - |u_1\gcdot \dots\gcdot u_\ell|$) and we obtain
  $$w=u_1\gcdot \dots\gcdot u_\ell \ins \phi_{I'}(v_1\gcdot \dots\gcdot v_r)$$ with
  $i_1' \leq |v_1|$ or $r = 0$.

  In the case of $v'\neq\epsilon$, the inequality $(|v'| + |I|) + 1 - i_p \leq |v_r|$
  is always true otherwise $w$ would not be irreducible.
\end{proof}




\begin{defi}\label{red_irreducible}
  A packed word $w$ is said to be \textbf{red-irreducible} if $w$ is irreducible
  and the equality $w = u \ins v$ implies that $u = \epsilon (\text{and } w = v)$.
\end{defi}

Here are all red-irreducible packed words of size $1$, $2$, $3$ and $4$ in
lexicographic order:
\[
  1,\quad 11,\quad111\ 121\ 132\ 212,
\]
\[
  1111\ 1121\ 1132\ 1211\ 1212\ 1221\ 1231\ 1232\ 1243\ 1312\ 1321\ 1322\
  1323\ 1332\
\]
\[
  1342\ 1423\ 1432\ 2112\ 2121\ 2122\ 2132\ 2143\ 2212\ 2312\ 2413\ 3123\
  3132\ 3213.
\]
\bigskip

Here are some useful lemmas on the operation $\ins$.

\begin{lem} \label{ins_pseudo_assoc_left}
  For any $u, v, w \in \PW$ with $w \neq \epsilon$, we have $u \ins (v \ins w) = (u \gcdot v) \ins w$.
  $$\input{figures/box_d/ins_of_ins} \qquad\scalebox{2}{=}\qquad\quad \input{figures/box_d/gd_in_ins}$$
\end{lem}

\begin{proof} 
  Let $u, v, w \in \PW$ with $w \neq \epsilon$ and let $w'$ and $I_w$ such that $w = \phi_{I_w}(w')$.
  \begin{align*}
    u \ins (v \ins w) &= u \ins (\phi_{I_w+|v|}(v \gcdot w'))\\
                      &= \phi_{I_w+|v|+|u|}(u \gcdot (v \gcdot w'))\\
                      &= \phi_{I_w+|v|+|u|}((u \gcdot v) \gcdot w')\\
                      &= (u \gcdot v) \ins \phi_{I_w}(w')\\
                      &= (u \gcdot v) \ins w. \qedhere
  \end{align*}
\end{proof}

\begin{lem} \label{ins_pseudo_assoc_right}
  For any $u, v, w \in \PW$ with $v \neq \epsilon$, we have $u \ins (v \gcdot w) = (u \ins v) \gcdot w$.
  $$\scalebox{1}{\input{figures/box_d/ins_of_gd} \qquad\scalebox{2}{=}\qquad\quad \input{figures/box_d/ins_in_gd}}$$
\end{lem}

\begin{proof}
  Let $u, v, w \in \PW$ with $v \neq \epsilon$ and let $v'$ and $I_v$ such that $v = \phi_{I_v}(v')$.
  \begin{align*}
    u \ins (v \gcdot w) &= u \ins \phi_{I_v}(v' \gcdot w)\\
                        &=  \phi_{I_v+|u|}(u \gcdot (v' \gcdot w))\\
                        &=  \phi_{I_v+|u|}((u \gcdot v') \gcdot w)\\
                        &=  \phi_{I_v+|u|}(u \gcdot v') \gcdot w\\
                        &= (u \ins v) \gcdot w. \qedhere
  \end{align*}
\end{proof}

\begin{rem}\label{rem:skew_dup_red}
  Adding the associativity of shifted concatenation
  $u \gcdot (v \gcdot w) = (u \gcdot v) \gcdot w$, the two operations $\ins$ and
  $\gcdot$ verify relations of the \textit{skew-duplicial operad}
  \cite{BurDel_dup}.
\end{rem}

\begin{coro} \label{ins_irr}
  For any $u, v \in \PW$, we have that $u \ins v$ is irreducible if and only if
  $v$ is irreducible.
\end{coro}

\begin{proof}
  By contradiction, if $v = v_1 \gcdot v_2$ then by
  \cref{ins_pseudo_assoc_right} $u \ins v = (u \ins v_1) \gcdot v_2$. Now if
  $u \ins v = w_1 \gcdot w_2$, as the position of the first maximum of $u \ins v$
  is greater than $|u|$ we have that $w_1 = w_1' \cdot w_1''$ such that
  $\pack(w_1') = u$. We also have that $\pack(w_1'') \gcdot w_2 = v$.
\end{proof}



\begin{prop} \label{red_fact_irr}
  For any word $w$, $w = u \ins v$ is the red-factorization of $w$ if and only
  if $v$ is red-irreducible.
\end{prop}

\begin{proof}
  Let $w \in \PW$ and let $u \ins v$ be the red-factorization of $w$. Let $v_1$ and
  $v_2$ such that $v = v_1 \ins v_2$, then $(u \gcdot v_1) \ins v_2 = w$ by
  \cref{ins_pseudo_assoc_left}, but in the red-factorization the size of $u$ is maximized
  so $|(u \gcdot v_1)| \leq |u|$ and then we have that $v_1 = \epsilon$ so $v$ is
  red-irreducible.

  Let $w \in \PW$ and let $u$ and $v$ such that $w = u \ins v$ and $v$ is
  red-irreducible. By contradiction, suppose that there exists $u', v'$ such
  that $w = u' \ins v'$ with $|u| < |u'|$ and $v' \neq \epsilon$. Then necessarily
  $u$ is a prefix of $u'$. Let $u''$ such that $u' = u \cdot u''$, then
  $pack(u'') \ins v' = v$. But $v$ is red-irreducible. So the size of $u$ is
  maximal if $v$ is red-irreducible.
\end{proof}

For the reader who is familiar with ordered set partitions, all the definitions
in \cref{Primal_Red} can be easily written with these. However in
\cref{Dual_Blue} it is easier to do all the definitions on packed words and in
\cref{sect4_bij} we must have the same object on both sides to explicit the
isomorphism. So we decided to stick to packed words.

\subsubsection{Red-forests from decomposed packed words using~\texorpdfstring{$\phi$}{phi}}
\label{red_forests}

We now apply recursively the red-factorization of the previous section to construct a
bijection between packed words and a certain kind of trees that we now define.

\begin{defi}
  An unlabeled \textbf{biplane tree} is an ordered tree (sometimes also called
  planar) whose children are organized in a pair of two (possibly empty) ordered
  forests, which we call the left and right forests, a forest being an ordered
  list of trees.
\end{defi}

In the picture, we naturally draw the children of the left (resp. right)
forest on the left (resp. right) of their father.
\begin{ex}\label{ex:biplan_tree}
  The biplane trees $\scalebox{0.5}{{ \newcommand{\nodea}{\node[draw,ellipse] (a) {}
;}\newcommand{\nodeb}{\node[draw,ellipse] (b) {}
;}\newcommand{\nodec}{\node[draw,ellipse] (c) {}
;}\begin{tikzpicture}[auto]
\matrix[column sep=.3cm, row sep=.3cm,ampersand replacement=\&]{
         \&         \& \nodea  \\ 
 \nodeb  \& \nodec  \&         \\
};

\path[ultra thick, red] (a) edge (b) edge (c);
\end{tikzpicture}}}$,
  $\scalebox{0.5}{{ \newcommand{\nodea}{\node[draw,ellipse] (a) {}
;}\newcommand{\nodeb}{\node[draw,ellipse] (b) {}
;}\newcommand{\nodec}{\node[draw,ellipse] (c) {}
;}\begin{tikzpicture}[auto]
\matrix[column sep=.3cm, row sep=.3cm,ampersand replacement=\&]{
         \& \nodea  \&         \\ 
 \nodeb  \&         \& \nodec  \\
};

\path[ultra thick, red] (a) edge (b) edge (c);
\end{tikzpicture}}}$ and
  $\scalebox{0.5}{{ \newcommand{\nodea}{\node[draw,ellipse] (a) {}
;}\newcommand{\nodeb}{\node[draw,ellipse] (b) {}
;}\newcommand{\nodec}{\node[draw,ellipse] (c) {}
;}\begin{tikzpicture}[auto]
\matrix[column sep=.3cm, row sep=.3cm,ampersand replacement=\&]{
 \nodea  \&         \&         \\ 
         \& \nodeb  \& \nodec  \\
};

\path[ultra thick, red] (a) edge (b) edge (c);
\end{tikzpicture}}}$ are different.
  Indeed in the first case, the left forest contains two trees and the right
  forest is empty, in the second case both forests contain exactly one tree
  while in the third case the left forest is empty and the right contains two
  trees. Here is an example of a bigger biplane tree where the root has two
  trees in both left and right forests
  \scalebox{0.5}{{ \newcommand{\nodea}{\node[draw,ellipse] (a) {}
;}\newcommand{\nodeb}{\node[draw,ellipse] (b) {}
;}\newcommand{\nodec}{\node[draw,ellipse] (c) {}
;}\newcommand{\noded}{\node[draw,ellipse] (d) {}
;}\newcommand{\nodee}{\node[draw,ellipse] (e) {}
;}\newcommand{\nodef}{\node[draw,ellipse] (f) {}
;}\newcommand{\nodeg}{\node[draw,ellipse] (g) {}
;}\newcommand{\nodeh}{\node[draw,ellipse] (h) {}
;}\newcommand{\nodei}{\node[draw,ellipse] (i) {}
;}\newcommand{\nodej}{\node[draw,ellipse] (j) {}
;}
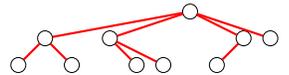
\begin{tikzpicture}[auto]
\matrix[column sep=.3cm, row sep=.3cm,ampersand replacement=\&]{
         \&         \&         \&         \&         \&         \&         \& \nodea  \&         \&         \&         \\ 
         \& \nodeb  \&         \&         \& \nodee  \&         \&         \&         \&         \& \nodeh  \& \nodej  \\ 
 \nodec  \&         \& \noded  \&         \&         \& \nodef  \& \nodeg  \&         \& \nodei  \&         \&         \\
};

\path[ultra thick, red] (b) edge (c) edge (d)
	(e) edge (f) edge (g)
	(h) edge (i)
	(a) edge (b) edge (e) edge (h) edge (j);
\end{tikzpicture}}}.
\end{ex}
  
\begin{defi}\label{skeleton}
  A \textbf{skeleton biplane tree} is a biplane tree where no node has a right
  forest.
\end{defi}

These skeleton biplane trees can also be seen as planar trees. In
\cite{Foissy_2011} we have planar trees recursively labeled by planar
trees. Skeleton biplane trees are similar to these planar trees, we prefer to
see them as biplane tree with no right forest in order to keep some
constistency.

\begin{defi}
  The \textbf{size} of a biplane tree is the number of node in the tree.
\end{defi}

\begin{rem}
  Biplane forests $\mathfrak{F}$ (\textit{i.e.} ordered list of biplane trees
  $\mathfrak{T}$) are counted by the sequence $A001764$ in OEIS~\cite{oeis} whose
  explicit formula is $a(n) = \binom{3n}{n}/(2n+1)$. Biplane forests are in
  bijection with ternary trees. The bijection is the following, in a biplane
  forest a node has a first left child and a first right child and a right
  brother. A consequence is that unlabeled biplane trees are counted by the
  sequence $A006013$ in OEIS~\cite{oeis} whose explicit formula is
  $a(n) = \binom{3n+1}{n}/(n+1)$. Indeed, biplane trees are in bijection with
  pair of ternary trees. Here is an example of the ternary tree in bijection
  with the biplane forest constituted of one tree, the big biplane tree in
  \cref{ex:biplan_tree}
  \scalebox{0.5}{{ \newcommand{\nodea}{\node[draw,ellipse] (a) {}
;}\newcommand{\nodeb}{\node[draw,ellipse] (b) {}
;}\newcommand{\nodec}{\node[draw,ellipse] (c) {}
;}\newcommand{\noded}{\node[draw,ellipse] (d) {}
;}\newcommand{\nodee}{\node[draw,ellipse] (e) {}
;}\newcommand{\nodef}{\node[draw,ellipse] (f) {}
;}\newcommand{\nodeg}{\node[draw,ellipse] (g) {}
;}\newcommand{\nodeh}{\node[draw,ellipse] (h) {}
;}\newcommand{\nodei}{\node[draw,ellipse] (i) {}
;}\newcommand{\nodej}{\node[draw,ellipse] (j) {}
;}\newcommand{\nodear}{\node (ar) {}
;}\newcommand{\nodecl}{\node (cl) {}
;}\newcommand{\nodecc}{\node (cc) {}
;}\newcommand{\nodecr}{\node (cr) {}
;}\newcommand{\nodedl}{\node (dl) {}
;}\newcommand{\nodedc}{\node (dc) {}
;}\newcommand{\nodedr}{\node (dr) {}
;}\newcommand{\nodeel}{\node (el) {}
;}\newcommand{\nodeer}{\node (er) {}
;}\newcommand{\nodefl}{\node (fl) {}
;}\newcommand{\nodefc}{\node (fc) {}
;}\newcommand{\nodegl}{\node (gl) {}
;}\newcommand{\nodegc}{\node (gc) {}
;}\newcommand{\nodegr}{\node (gr) {}
;}\newcommand{\nodehc}{\node (hc) {}
;}\newcommand{\nodeil}{\node (il) {}
;}\newcommand{\nodeic}{\node (ic) {}
;}\newcommand{\nodeir}{\node (ir) {}
;}\newcommand{\nodejl}{\node (jl) {}
;}\newcommand{\nodejc}{\node (jc) {}
;}\newcommand{\nodejr}{\node (jr) {}
;}
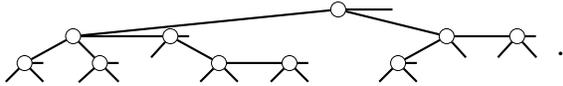
\begin{tikzpicture}[auto]
\matrix[column sep=.3cm, row sep=.3cm,ampersand replacement=\&]{
         \&         \&         \&         \&         \&         \&         \&         \&         \&         \&         \&         \&         \&         \& \nodea  \&         \&         \& \nodear \&         \&         \&         \&         \&         \&  \\ 
         \&         \&         \& \nodeb  \&         \&         \&         \& \nodee  \& \nodeer \&         \&         \&         \&         \&         \&         \&         \&         \&         \&         \& \nodeh  \&         \&         \& \nodej  \& \nodejr \\ 
         \& \nodec  \& \nodecr \&         \& \noded  \& \nodedr \& \nodeel \&         \&         \& \nodef  \&         \&         \& \nodeg  \& \nodegr \&         \&         \&         \& \nodei  \& \nodeir \&         \& \nodehc \& \nodejl \&         \& \nodejc \\ 
 \nodecl \&         \& \nodecc \& \nodedl \&         \& \nodedc \&         \&         \& \nodefl \&         \& \nodefc \& \nodegl \&         \& \nodegc \&         \&         \& \nodeil \&         \& \nodeic \&         \&         \&         \&         \&  \\ 
};

\path[ultra thick, black]
	(a) edge (b) edge (h) edge (ar)
        (b) edge (c) edge (d) edge (e)
        (c) edge (cl) edge (cc) edge (cr)
        (d) edge (dl) edge (dc) edge (dr)
        (e) edge (el) edge (f) edge (er)
        (f) edge (fl) edge (fc) edge (g)
        (g) edge (gl) edge (gc) edge (gr)
        (h) edge (i) edge (hc) edge (j)
        (i) edge (il) edge (ic) edge (ir)
        (j) edge (jl) edge (jc) edge (jr);
      \end{tikzpicture}}
}.
\end{rem}
\begin{table}[!h]
  \[
    \begin{tabular}{|c|c|c|c|c|c|c|c|c|c||l|}
      \hline
      $n$ & 0 & 1 & 2 & 3 & 4 & 5 & 6 & 7 & 8 & OEIS \\
      \hline
      $\mathfrak{F}_n$ & 1 & 1 & 3 & 12 & 55 & 273 & 1 428 & 7 752 & 43 263 & A001764 \\
      \hline
      $\mathfrak{T}_n$ & 0 & 1 & 2 & 7 & 30 & 143 & 728 & 3 876 & 21 318 & A006013 \\
      \hline
    \end{tabular} 
  \]
  \caption{Number of biplane forests and biplane trees.}
  \label{oeis_biplan}
\end{table}

\begin{rem}\label{rem:L}
  As we can see on OEIS~\cite{oeis}, sequence $A006013$ which counts unlabeled
  biplane trees is the dimensions of the free L-algebra on one generator
  (see~\cite{Leroux}). It would be interesting to investigate the link between
  $L$-algebras and bidendriform bialgebras using biplane trees.
\end{rem}

In our construction we will deal with labeled biplane trees with colored
edges. For a labeled biplane tree, we denote by $\Nodered(x, f_\ell, f_r)$ the tree
whose edges are colored in red, root is labeled by $x$ and whose left
(resp. right) forest is given by $f_\ell$ (resp. $f_r$). We also denote
by~$[t_1, \dots t_k]$ a forest of $k$ trees. The edge color (for now, only red)
will play a role later in the paper.



\begin{ex}
  $\Nodered((1),\ [],[]) =$ \scalebox{0.5}{{ \newcommand{\nodea}{\node[draw,ellipse] (a) {$\,1\,$}
;}\begin{tikzpicture}[baseline={([yshift=-1ex]current bounding box.center)}]
\matrix[column sep=.3cm,row sep=.3cm,ampersand replacement=\&]{
 \nodea  \\
};
\path[ultra thick,red] ;
\end{tikzpicture}}},
  and $\Nodered((1,3),\ [],[\Nodered((1),\ [],[])]) =$
  \scalebox{0.5}{{ \newcommand{\nodea}{\node[draw,ellipse] (a) {$1, 3$}
;}\newcommand{\nodeb}{\node[draw,ellipse] (b) {$\,1\,$}
;}\begin{tikzpicture}[baseline={([yshift=-1ex]current bounding box.center)}]
\matrix[column sep=.3cm,row sep=.3cm,ampersand replacement=\&]{
         \& \nodea  \&         \\ 
         \&         \& \nodeb  \\
};
\path[ultra thick,red] (a) edge (b);
\end{tikzpicture}}}.
\end{ex}

We now apply recursively the global descent decomposition and the
red-factorization of~\cref{gd_fact,red_fact}. We obtain an algorithm which takes
a packed word and returns a biplane forest where nodes are decorated by
red-irreducible packed words:

\begin{defi}\label{def:construction_red_skeleton}
  We now define two functions $\Frske$ and $\Trske$. These functions transform
  respectively a packed word and an irreducible packed word into respectively a
  skeleton biplane forest and a skeleton biplane tree labeled by red-irreducible
  words. These functions are defined in a mutual recursive way as follows:
  \begin{itemize}
  \item $\Frske(\epsilon) = []$ (empty forest),
  \item for any packed word $w$, let $w_1\gcdot w_2\gcdot\dots\gcdot w_k$ be the
    global descent decomposition of $w$, then
    $\Frske(w) \eqdef [\Trske(w_1), \Trske(w_2), \dots, \Trske(w_k)]$.
    $$\scalebox{0.7}{\input{figures/box_d/gd_fact} \qquad\scalebox{2}{$\to$}\qquad\quad $\overbrace{\input{figures/arbres_pw/arbrew1_ske}}^{\Trske(w_1)}$
      $\overbrace{\input{figures/arbres_pw/arbrew2_ske}}^{\Trske(w_2)}$
      $\cdots$ $\overbrace{\input{figures/arbres_pw/arbrew1_ske}}^{\Trske(w_k)}$}.$$
  \item for any irreducible packed word $w$, let $w=u \ins v$ be the
    red-factorization of $w$. We define
    $\Trske(w) \eqdef \Nodered(v, \Frske(u), [])$.
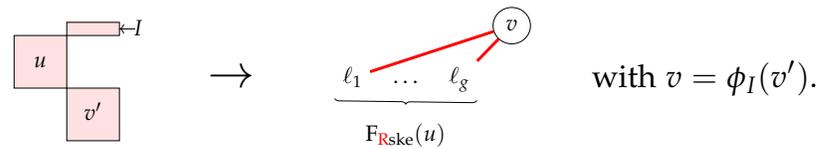
    $$\scalebox{0.7}{\input{figures/box_d/ins} \qquad\scalebox{2}{$\to$}\qquad\quad { \newcommand{\nodea}{\node[draw,ellipse] (a) {$v$}
;}\newcommand{\nodeb}{\node (b) {$\ell_1$}
;}\newcommand{\nodebc}{\node (bc) {$\ldots$}
;}\newcommand{\nodec}{\node (c) {$\ell_g$}
;}\begin{tikzpicture}[auto]
\matrix[column sep=.3cm, row sep=.3cm,ampersand replacement=\&]{
      \&  \&  \& \nodea  \& \&    \&    \\ 
 \nodeb  \& \nodebc \& \nodec \&  \& \& \& \\ 
};

\path[ultra thick, red] (a) edge (b) edge (c); 
\draw[snake=brace,raise snake = 4mm,mirror snake] (b.west) -- (c.east);
\node[below of = bc] (F) {$\Frske(u)$};

\end{tikzpicture}}} \text{with $v = \phi_I(v')$}.$$
  \end{itemize}
\end{defi}

\begin{ex}\label{ex:red_skeleton}
  Let $w = 876795343912$, the global descent decomposition of~\cref{gd_fact}
  gives $w = w_1\gcdot w_2$ with $w_1 = 6545731217$ and $w_2 = 12$. Now, we have
  the red-factorization of $w_1$ and $w_2$ using~\cref{red_fact} as
  \[
    w_1 = 3212 \ins \phi_{1,6}( 3121 ) = (1 \gcdot 212) \ins 431214
    \quad\text{and}\quad
    w_2 = 1 \ins \phi_1(\epsilon) = 1 \ins 1.
  \]
  $$\scalebox{0.3}{\input{figures/box_d/pw_red_ex}}$$
  It gives the following forest:
  \begin{align*}
    \Frske(876795343912) & = [\Trske(6545731217),\ \Trske(12)] \\
                         & =
                           \scalebox{0.7}{{ \newcommand{\nodea}{\node[draw,ellipse] (a) {$431214$}
;}\newcommand{\nodeb}{\node (b) {$\Frske(3212)$}
;}\begin{tikzpicture}[auto]
\matrix[column sep=.3cm, row sep=.3cm,ampersand replacement=\&]{
  \& \nodea \\
  \nodeb \& \\
};

\path[ultra thick, red] (a) edge (b);
\end{tikzpicture}}
}
                           \scalebox{0.7}{{ \newcommand{\nodea}{\node[draw,ellipse] (a) {$\,1\,$}
;}\newcommand{\nodeb}{\node (b) {$\Frske(1)$}
;}\begin{tikzpicture}[auto]
\matrix[column sep=.3cm, row sep=.3cm,ampersand replacement=\&]{
         \& \nodea  \\
 \nodeb  \&         \\
};

\path[ultra thick, red] (a) edge (b);
\end{tikzpicture}}
} \\ 
                         & =
                           \scalebox{0.7}{\input{figures/arbres_pw/arbre876795343912ske}}.
  \end{align*}
\end{ex}

\begin{defi}\label{def:red_skeleton}
  A labeled biplane forest (resp. tree) is a \textbf{red-skeleton forest}
  (resp. \textbf{tree}) if it is labeled by red-irreducible words and no node
  has a right child.
\end{defi}

We want to prove that the functions $\Frske$ and $\Trske$ are bijections. To do
that we first define two functions that are the inverses.

\begin{defi}\label{def:construction_red_skeleton_star}
  We now define two functions $\Frskestar$ and $\Trskestar$ that
  transform respectively a red-skeleton forest and tree into packed words. These
  functions are defined in a mutual recursive way as follows:
  \begin{itemize}
  \item $\Frskestar([]) = \epsilon$,
  \item for any red-skeleton forest $f = [t_1, \dots, t_k]$, we define \newline
    $\Frskestar(f) \eqdef \Trskestar(t_1)\gcdot \dots \gcdot\Trskestar(t_k)$.
    $$\scalebox{0.7}{$\overbrace{\input{figures/arbres_pw/arbrew1_ske}}^{t_1}$
      $\overbrace{\input{figures/arbres_pw/arbrew2_ske}}^{t_2}$
      $\cdots$ $\overbrace{\input{figures/arbres_pw/arbrew1_ske}}^{t_k}$
      \qquad\scalebox{2}{$\to$}\qquad\quad \input{figures/box_d/gd_fact}} \text{with $w_i = \Trskestar(t_i)$}.$$
  \item for any red-skeleton tree $t = \Nodered(v, f_\ell, [])$, we define \newline
    $\Trskestar(t) \eqdef \Frskestar(f_\ell) \ins v$.
    $$\scalebox{0.7}{{ \newcommand{\nodea}{\node[draw,ellipse] (a) {$v$}
;}\newcommand{\nodeb}{\node (b) {$\ell_1$}
;}\newcommand{\nodebc}{\node (bc) {$\ldots$}
;}\newcommand{\nodec}{\node (c) {$\ell_g$}
;}\begin{tikzpicture}[auto]
\matrix[column sep=.3cm, row sep=.3cm,ampersand replacement=\&]{
      \&  \&  \& \nodea  \& \&    \&    \\ 
 \nodeb  \& \nodebc \& \nodec \&  \& \& \& \\ 
};

\path[ultra thick, red] (a) edge (b) edge (c); 
\draw[snake=brace,raise snake = 4mm,mirror snake] (b.west) -- (c.east);
\node[below of = bc] (F) {$f_\ell$};

\end{tikzpicture}} \qquad\scalebox{2}{$\to$}\qquad\quad \input{figures/box_d/ins}} \text{with $v = \phi_I(v')$ and $u = \Frskestar(f_\ell)$}.$$
  \end{itemize}
\end{defi}
\begin{lem}\label{lem:red_skeleton_bij}
  The functions $\Frske$ and $\Frskestar$ (resp. $\Trske$ and $\Trskestar$) are
  two converse bijections between packed words and red-skeleton forests
  (resp. irreducible packed words and red-skeleton trees). That is to
  say $\Frskeinv = \Frskestar$ and $\Trskeinv = \Trskestar$.
\end{lem}

\begin{proof}
  We start to prove that domain and codomain are as announced (see Items (a) and
  (b) bellow), then we prove that the functions $\Frske$ and $\Frskestar$
  (resp. $\Trske$ and $\Trskestar$) are inverse to each other (see Items (c) and
  (d) bellow).

  (a) By \cref{def:construction_red_skeleton}, a forest (resp. tree) obtain by
  $\Frske$ (resp. $\Trske$) is a red-skeleton forest (resp. tree). Indeed,
  thanks to \cref{red_fact_irr} nodes are labeled by red-irreducible words
  because a red-factorization is done and nodes have no right children.
  
  (b) We prove by a mutual induction that $\Frskestar$ returns a packed word and
  that $\Trskestar$ returnes an irreducible packed word. Indeed, we do an
  induction on the size of the forest or the tree. Here is our induction
  hypothesis for $n \in \NN$:
  \begin{equation}\label{HR_rskestar}
    \begin{array}{l}
      \forall t, \text{red-skeleton tree of size } \leq n, \Trskestar(t) \text{ is an irreducible packed word,} \\
      \forall f, \text{red-skeleton forest of size } \leq n, \Frskestar(f) \text{ is a packed word.}
    \end{array}
  \end{equation}
  
  The base case ($n = 0$) is given by the first item of
  \cref{def:construction_red_skeleton_star} (\textit{i.e.}
  $\Frskestar([]) = \epsilon$ the empty packed word).
  
  \noindent Now, let us fix $n \geq 1$ and suppose that the hypothesis
  (\ref{HR_rskestar}) holds. Let $f = [t_1, \dots, t_k]$ be a red-skeleton forest of
  size $n + 1$.
  
  \noindent $\bullet$ If $k = 1$, then $f$ is reduced to a single tree $t$. We need
  to prove that $\Trskestar(t)$ is an irreducible packed word (which also
  gives that $\Frskestar(f)$ is a packed word as in this case
  $\Frskestar(f) = \Trskestar(t)$). Let $t = \Nodered(v, f_\ell, [])$ be a
  red-skeleton tree of size $n + 1$ (notice that the word $v$ can be of any
  size). The size of $f_\ell$ is $n$, so by induction $\Frskestar(f_\ell)$ is a
  packed word, as $v$ is red-irreducible it is by definition irreducible so by
  \cref{ins_irr} $\Trskestar(t) = \Frskestar(f_\ell) \ins v$ is irreducible.
  
  \noindent $\bullet$ If $k \geq 2$, \textit{i.e.}, the forest contains at least two
  trees, since all trees are of size at least one, $t_1, \dots, t_k$ are at
  must of size $n$, so we have by induction that
  $\Trskestar(t_1), \dots \Trskestar(t_k)$ are irreducible packed
  words. $\Frskestar(f)$ is the shifted concatenation of
  $\Trskestar(t_1), \dots \Trskestar(t_k)$ and thus it is a packed word.

  (c) We now prove by a mutual induction on the size of the forest or the tree
  that, for $n \in \NN$:
  \begin{equation}
    \begin{array}{l}
      \forall t, \text{red-skeleton tree of size } \leq n, \Trske(\Trskestar(t)) = t,\\
      \forall f, \text{red-skeleton forest of size } \leq n, \Frske(\Frskestar(f)) = f.
    \end{array}
    \label{H_n}
  \end{equation}
  
  The base case ($n = 0$) is given by the first item of
  \cref{def:construction_red_skeleton,def:construction_red_skeleton_star} as
  $\Frske(\Frskestar([])) = \Frske(\epsilon) = []$.

  \noindent Now let us fix $n \geq 1$ and suppose that the hypothesis (\ref{H_n})
  holds. Let $f = [t_1, \dots, t_k]$ be a red-skeleton forest of size $n + 1$.

  \noindent $\bullet$ If $k = 1$, then the forest $f$ is reduced to a single tree
  $t$, then it is sufficient to prove $\Trske(\Trskestar(t)) = t$ (as in this
  case $\Frskestar(f) = \Trskestar(t)$). Let $t = \Nodered(v, f_\ell, [])$ a
  red-skeleton tree of size $n + 1$. As the label $v$ is a red-irreducible
  packed word, with the induction hypothesis on $\Frskestar(f_\ell)$ and with
  \cref{red_fact_irr}, $\Frskestar(f_\ell) \ins v$ is the red-factorization so:
  \begin{align*}
    \Trske(\Trskestar(t)) &= \Trske(\Frskestar(f_\ell) \ins v)\\
                          &= \Nodered(v, \Frske(\Frskestar(f_\ell), []))\\
                          &= \Nodered(v, f_\ell, []) = t.
  \end{align*}
  
  \noindent $\bullet$ If $k \geq 2$, since all trees are of size at least one, they are
  at most of size $n$, so we have by induction that:
  \begin{align*}
    \Frske(\Frskestar(f)) &= \Frske(\Trskestar(t_1)\gcdot \dots \gcdot\Trskestar(t_k))\\
                          & \qquad \text{as $\Trskestar(t_i)$ are irreducible packed words}\\
                          &= [\Trske(\Trskestar(t_1)), \dots ,\Trske(\Trskestar(t_k))]\\
                          &= [t_1, \dots, t_k] = f.
  \end{align*}
  
  (d) Finally we prove by a mutual induction on the size of the word $w$ that,
  for $n \in \NN$:
  \begin{equation}
    \begin{array}{l}
      \forall v \in \PW, \text{irreducible packed word of size } \leq n, \Trskestar(\Trske(v)) = v,\\
      \forall w \in \PW, \text{packed word of size } \leq n, \Frskestar(\Frske(w)) = w.
    \end{array}
    \label{H^*_n}
  \end{equation}

  The base case ($n = 0$) is given by the first item of
  \cref{def:construction_red_skeleton,def:construction_red_skeleton_star} as
  $\Frskestar(\Frske(\epsilon)) = \Frskestar([]) = \epsilon$.

  \noindent Now let us fix $n \geq 1$ and suppose that the hypothesis \eqref{H^*_n}
  holds.  Let $w \in \PW_{n+1}$ a packed word of size $n + 1$. Let
  $w = w_1\gcdot w_2\gcdot\dots\gcdot w_k$ be the global descent decomposition
  of $w$.

  \noindent $\bullet$ If $k = 1$, the packed word $w$ is irreducible then
  $\Frske(w) = [\Trske(w)]$ so we need to prove that
  $\Trskestar(\Trske(w)) = w$. Let $w=u \ins v$ be the red-factorization of
  $w$, then we can use the induction hypothesis on $u$, indeed as $v$ is not
  empty the size of $u$ is smaller than $n$:
  \begin{align*}
    \Trskestar(\Trske(w)) &= \Trskestar(\Nodered(v, \Frske(u), []))\\
                          &= \Frskestar(\Frske(u)) \ins v\\
                          &= u \ins v = w.
  \end{align*}

  \noindent $\bullet$ If $k \geq 2$, then we use the induction hypothesis on each factors, so
  we have:
  \begin{align*}
    \Frskestar(\Frske(w)) &= \Frskestar([\Trske(w_1), \dots, \Trske(w_k)])\\
                          &= \Trskestar(\Trske(w_1))\gcdot \dots \gcdot \Trskestar(\Trske(w_k))\\
                          &= w_1 \gcdot \dots \gcdot w_k = w. \qedhere
  \end{align*}
\end{proof}

Now that we have the red-skeleton, we will add right forests to every nodes to
obtain biplane trees. For every node, if $v$ is the red-irreducible word in
label, with $\phi_I(v') = v$, then $I$ is the new label and $\Fr(v')$ is the new
right forest.

Here is the formal definition of $\Fr(w)$ and $\Tr(w)$ which are very similar to
\cref{def:construction_red_skeleton}, only the third item is different. The
labels are now lists of integers.

\begin{defi}\label{def:construction_red_packed}
  The forest~$\Fr(w)$~(resp. tree~$\Tr(w)$) associated to a packed word (resp.
  irreducible packed word) $w$ are defined in a mutual recursive way as follows:
  \begin{itemize}
  \item $\Fr(\epsilon) = []$ (empty forest),
  \item for any packed word $w$, let
    $w_1\gcdot w_2\gcdot\dots\gcdot w_k$ be the global descent decomposition
    of $w$, then $\Fr(w) \eqdef [\Tr(w_1), \Tr(w_2), \dots, \Tr(w_k)]$.
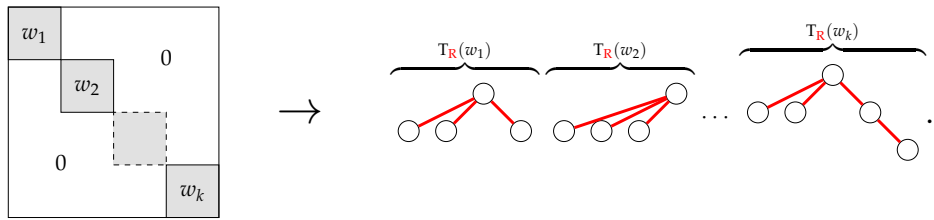
    $$\scalebox{0.7}{\input{figures/box_d/gd_fact} \qquad\scalebox{2}{$\to$}\qquad\quad $\overbrace{{ \newcommand{\nodea}{\node[draw,ellipse] (a) {}
;}\newcommand{\nodeb}{\node[draw,ellipse] (b) {}
;}\newcommand{\nodec}{\node[draw,ellipse] (c) {}
;}\newcommand{\noded}{\node[draw,ellipse] (d) {}
;}\begin{tikzpicture}[baseline={([yshift=-1ex]current bounding box.center)}]
\matrix[column sep=.3cm,row sep=.3cm,ampersand replacement=\&]{
         \&         \& \nodea \& \\ 
 \nodeb  \& \nodec  \&        \& \noded  \\
};
\path[ultra thick,red] (a) edge (b) edge (c) edge (d);
\end{tikzpicture}}}^{\Tr(w_1)}$
      $\overbrace{\input{figures/arbres_pw/arbrew2}}^{\Tr(w_2)}$ $\cdots$ $\overbrace{{ \newcommand{\nodea}{\node[draw,ellipse] (a) {}
;}\newcommand{\nodeb}{\node[draw,ellipse] (b) {}
;}\newcommand{\nodec}{\node[draw,ellipse] (c) {}
;}\newcommand{\noded}{\node[draw,ellipse] (d) {}
;}\newcommand{\nodee}{\node[draw,ellipse] (e) {}
;}\begin{tikzpicture}[baseline={([yshift=-1ex]current bounding box.center)}]
\matrix[column sep=.3cm,row sep=.3cm,ampersand replacement=\&]{
         \&         \& \nodea \& \& \\ 
 \nodeb  \& \nodec  \&        \& \noded \& \\
 \&  \&        \& \& \nodee  \\
};
\path[ultra thick,red] (a) edge (b) edge (c) edge (d)
(d) edge (e);
\end{tikzpicture}}}^{\Tr(w_k)}$}.$$
  \item for any irreducible packed word $w$, we define
    $\Tr(w) \eqdef \Nodered(I, \Fr(u), \Fr(v'))$ where $w=u\ins\phi_I(v')$ is the
    red-factorization of $w$.
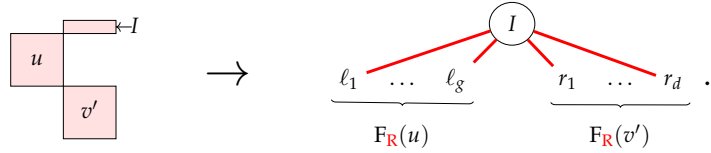
    $$\scalebox{0.7}{\input{figures/box_d/ins} \qquad\scalebox{2}{$\to$}\qquad\quad { \newcommand{\nodea}{\node[draw,ellipse] (a) {$\,I\,$}
;}\newcommand{\nodeb}{\node (b) {$\ell_1$}
;}\newcommand{\nodebc}{\node (bc) {$\ldots$}
;}\newcommand{\nodec}{\node (c) {$\ell_g$}
;}\newcommand{\nodee}{\node (e) {$r_1$}
;}\newcommand{\nodeef}{\node (ef) {$\ldots$}
;}\newcommand{\nodef}{\node (f) {$r_d$}
;}\begin{tikzpicture}[auto]
\matrix[column sep=.3cm, row sep=.3cm,ampersand replacement=\&]{
      \&  \&  \& \nodea  \& \&    \&    \\ 
 \nodeb  \& \nodebc \& \nodec \&  \& \nodee \& \nodeef \& \nodef  \\
};

\path[ultra thick, red] (a) edge (b) edge (c) edge (e) edge (f);
\draw[snake=brace,raise snake = 4mm,mirror snake] (b.west) -- (c.east);
\node[below of = bc] (F) {$\Fr(u)$};
\draw[snake=brace,raise snake = 4mm,mirror snake] (e.west) -- (f.east);
\node[below of = ef] (F) {$\Fr(v')$};

\end{tikzpicture}}}.$$
  \end{itemize}
\end{defi}



\begin{ex}
  Consider again $w = 876795343912$. We start from the red-skeleton forest from
  \cref{ex:red_skeleton}.
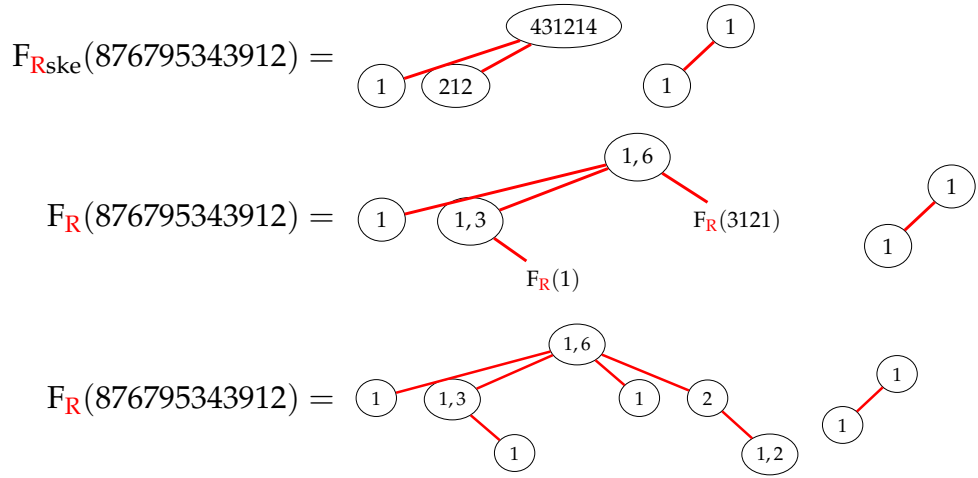
  \begin{align*}
    \Frske(876795343912) & = \scalebox{0.7}{\input{figures/arbres_pw/arbre876795343912ske}}\\
    \Fr(876795343912) & = \scalebox{0.7}{{ \newcommand{\nodea}{\node[draw,ellipse] (a) {$1, 6$}
;}\newcommand{\nodeb}{\node[draw,ellipse] (b) {$\,1\,$}
;}\newcommand{\nodec}{\node[draw,ellipse] (c) {$1, 3$}
;}\newcommand{\noded}{\node (d) {$\Fr(1)$}
;}\newcommand{\nodee}{\node (e) {$\Fr(3121)$}
;}\begin{tikzpicture}[auto]
\matrix[column sep=.3cm, row sep=.3cm,ampersand replacement=\&]{
         \&         \&         \&         \& \nodea  \&         \&         \&         \&         \\ 
 \nodeb  \&         \& \nodec  \&         \&         \& \nodee  \&         \&   \&         \\ 
         \&         \&         \& \noded  \&         \&         \&         \&         \&   \\
};

\path[ultra thick, red] (c) edge (d)
	(a) edge (b) edge (c) edge (e);
\end{tikzpicture}}
{  \newcommand{\nodea}{\node[draw,ellipse] (a) {$\,1\,$}
;}\newcommand{\nodeb}{\node[draw,ellipse] (b) {$\,1\,$}
;}\begin{tikzpicture}[auto]
\matrix[column sep=.3cm, row sep=.3cm,ampersand replacement=\&]{
         \& \nodea  \\ 
 \nodeb  \&         \\
};

\path[ultra thick, red] (a) edge (b);
\end{tikzpicture}}}\\
    \Fr(876795343912) & = \scalebox{0.6}{{ \newcommand{\nodea}{\node[draw,ellipse] (a) {$1, 6$}
;}\newcommand{\nodeb}{\node[draw,ellipse] (b) {$\,1\,$}
;}\newcommand{\nodec}{\node[draw,ellipse] (c) {$1, 3$}
;}\newcommand{\noded}{\node[draw,ellipse] (d) {$\,1\,$}
;}\newcommand{\nodee}{\node[draw,ellipse] (e) {$\,1\,$}
;}\newcommand{\nodef}{\node[draw,ellipse] (f) {$\,2\,$}
;}\newcommand{\nodeg}{\node[draw,ellipse] (g) {$1, 2$}
;}\begin{tikzpicture}[auto]
\matrix[column sep=.3cm, row sep=.3cm,ampersand replacement=\&]{
         \&         \&         \&         \& \nodea  \&         \&         \&         \&         \\ 
 \nodeb  \&         \& \nodec  \&         \&         \& \nodee  \&         \& \nodef  \&         \\ 
         \&         \&         \& \noded  \&         \&         \&         \&         \& \nodeg  \\
};

\path[ultra thick, red] (c) edge (d)
	(f) edge (g)
	(a) edge (b) edge (c) edge (e) edge (f);
\end{tikzpicture}}
{  \newcommand{\nodea}{\node[draw,ellipse] (a) {$\,1\,$}
;}\newcommand{\nodeb}{\node[draw,ellipse] (b) {$\,1\,$}
;}\begin{tikzpicture}[auto]
\matrix[column sep=.3cm, row sep=.3cm,ampersand replacement=\&]{
         \& \nodea  \\ 
 \nodeb  \&         \\
};

\path[ultra thick, red] (a) edge (b);
\end{tikzpicture}}}
  \end{align*}
\end{ex}

\begin{defi} \label{red_weight}
  Let $t$ be a labeled biplane tree. We write $t = \Nodered(I, f_\ell, f_r)$ where
  ${I = [i_1, \ldots, i_p]}$, $p > 0$, $(1 \leq i_1 < \dots < i_p)$,
  $f_\ell = [\ell_1, \ldots, \ell_g]$ and $f_r = [r_1, \ldots, r_d]$, which is depicted as follows:
  \begin{equation*}
    t = \scalebox{0.8}{{ \newcommand{\nodea}{\node[draw,ellipse] (a) {$\,I\,$}
;}\newcommand{\nodeb}{\node (b) {$\ell_1$}
;}\newcommand{\nodebc}{\node (bc) {$\ldots$}
;}\newcommand{\nodec}{\node (c) {$\ell_g$}
;}\newcommand{\nodee}{\node (e) {$r_1$}
;}\newcommand{\nodeef}{\node (ef) {$\ldots$}
;}\newcommand{\nodef}{\node (f) {$r_d$}
;}
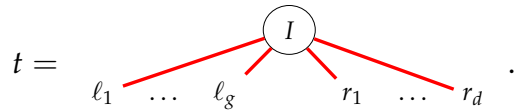
\begin{tikzpicture}[auto]
\matrix[column sep=.3cm, row sep=.3cm,ampersand replacement=\&]{
      \&  \&  \& \nodea  \& \&    \&    \\ 
 \nodeb  \& \nodebc \& \nodec \&  \& \nodee \& \nodeef \& \nodef  \\
};

\path[ultra thick, red] (a) edge (b) edge (c) edge (e) edge (f);

\end{tikzpicture}}}.
  \end{equation*}
  
  The \textbf{weight} of $t$ is recursively defined by
  $\omega(t) = p + \sum_{i=0}^{g}\omega(\ell_i) + \sum_{j=0}^{d}\omega(r_j)$. In particular, if
  $t$ is a single node then $\omega(t) = p$. 
  By extension, the \textbf{weight} of a forest is the sum of the weight of its
  trees.
\end{defi}

\begin{lem}\label{Frweight}
  The weight of a forest (resp. a tree) obtained by the functions $\Fr$
  (resp. $\Tr$) is equal to the size of the word, \textit{i.e.}  For
  all$w \in \PW$ then $\omega(\Fr(w)) = |w|$ and for all $w \in \PW$ with $w$ irreducible
  then $\omega(\Tr(w)) = |w|$.
\end{lem}

\begin{proof}
  We prove by induction with the following hypothesis, for $n \in \NN$:
  \begin{equation}
    \begin{array}{rl}
      \forall w \in \PW_n, & \omega(\Fr(w)) = |w|,\\
      \forall w \in \PW_n \text{ with $w$ irreducible}, & \omega(\Tr(w)) = |w|.
    \end{array}
    \label{H_nsize}
  \end{equation}
  The base case is given by the first item of
  \cref{def:construction_red_packed} as $\Fr(\epsilon) = []$ and $\omega([]) = |\epsilon| = 0$.

  \noindent Let us fix $n \geq 1$ and suppose that the hypothesis
  (\ref{H_nsize}) holds. Let $w \in \PW_{n+1}$ and
  ${w = w_1\gcdot w_2\gcdot\dots\gcdot w_k}$ be the global descent decomposition
  of $w$.

  \noindent $\bullet$ If $k = 1$, we have $\Fr(w) = [\Tr(w)]$. Let $w = u \ins \phi_I(v)$
  with $I = [i_1, \dots, i_p], p > 0$ be the red-factorization of $w$, then
  $$\omega(\Tr(w)) = \omega(\Nodered(I, \Fr(u), \Fr(v))) = p + \omega(\Fr(u)) + \omega(\Fr(v)).$$ As
  $p > 0$, the sizes of $u$ and $v$ are at most $n$, by induction
  ${\omega(\Tr(w)) = p + |u| + |v| = |w|}$.

  \noindent $\bullet$ If $k \geq 2$, by induction on each factors, we have that
  \[\omega(\Fr(w)) = \omega(\Tr(w_1)) + \dots + \omega(\Tr(w_k)) = |w_1| + \dots + |w_k| = |w|.\qedhere\]
\end{proof}

\begin{defi} \label{red_packed_forest}
  Using the same notations as in previous \cref{red_weight}, we say that $t$ is a
  \textbf{red-packed tree} if it satisfies:
  \begin{equation*}
    \left\{
      \begin{array}{l}
        d = 0,\\
        i_k = k \text{ for all }k \leq p,\\
        \ell_1, \ldots, \ell_g \text{ are red-packed trees.}
      \end{array}
    \right.
    \text{ or }
    \left\{
      \begin{array}{l}
        d \geq 1,\\
        \begin{array}{lcl}
          \!\!\!1 \leq &i_1 &\leq \omega(r_1),\\
          \!\!\!1 \leq &p + \omega(f_r) + 1 - i_p &\le \omega(r_d),
        \end{array}\\
        \ell_1, \ldots, \ell_g \text{ and } r_1, \ldots, r_d \text{ are red-packed trees.}
      \end{array}
    \right.
  \end{equation*}
  An ordered list of red-packed trees is a \textbf{red-packed forest}.
\end{defi}

\begin{rem}\label{red_ske_flat}
  Red-skeleton trees can be interpreted as flattened representations of
  red-packed trees. Symmetrically, red-packed trees can be interpreted as
  unfolded representations of red-skeleton trees. We use the operation $\phi_I$ to
  change between red-packed and red-skeleton trees.
\end{rem}

\begin{nota}\label{notation_red} 
  From now on, we use these notations:
  \begin{itemize}
  \item $\PForestredn$ the set of red-packed forests of weight $n$,
    ($\PForestredn \eqdef \{\Fr(w)\}_{w\in\PW_n}$),
  \item $\PTreeredn$ the set of red-packed trees of weight $n$,
    ($\PTreeredn \eqdef \{\Tr(w)\}_{w\in\PW_n}$ with $w$ irreducible),
  \item $\PParticularredn$ the set of red-packed trees of weight $n$ such that
    the left forest of the root is empty, ($\PParticularredn \eqdef \{\Tr(w)\}_{w\in\PW_n}$ with $w$
    red-irreducible). In particular, the red-skeleton of a tree of
    $\PParticularredn$ consist of a single node labeled by a red-irreducible word.
  \end{itemize}
\end{nota}

\begin{rem} \label{rem:PParticularredn}
  The set $\PParticularredn$ can be described as a disjointed union of sets
  depending on $I= [i_1, \dots, i_p]$. Let $\PForestredn^I$ denote the set of
  red-packed forests of weight $n$ that can be right children of a node labeled
  by $I$ (see~\cref{red_packed_forest} for conditions), we have the following description:
  \begin{align} \label{eq:PParticularredn}
    \PParticularredn = \bigsqcup_I~\{\Nodered(I,[],f_r)~|~f_r \in \PForestrednmp^I\}.
  \end{align}
\end{rem}

Analogously, we use $\PForestredske$, $\PTreeredske$ and $\PParticularredske$
for red-skeleton forests, trees and trees with only one node.

\noindent We can remark that for $n = 1$ we have
$\PParticularredone = \PTreeredone = \PForestredone$ and
$\forall n > 1, \PParticularredn \subsetneq \PTreeredn \subsetneq \PForestredn$.

As with \cref{def:construction_red_skeleton_star}, we want to prove that the
functions $\Fr$ and $\Tr$ are bijections. To do that we first define the two
inverse functions.

\begin{defi}\label{def:Fr-1}
  We define here the functions $\Frstar$~(resp. $\Trstar$) that transform a
  red-packed forest $f$ (resp. tree $t$) into a packed word. We reverse all
  instructions of \cref{def:construction_red_packed} as follows:
  \begin{itemize}
  \item $\Frstar([]) = \epsilon$,
  \item for any non empty red-packed forest $f = [t_1, t_2\dots, t_k]$, then \newline
    $\Frstar(f) = \Trstar(t_1)\gcdot \Trstar(t_2)\gcdot\dots\gcdot \Trstar(t_k)$.
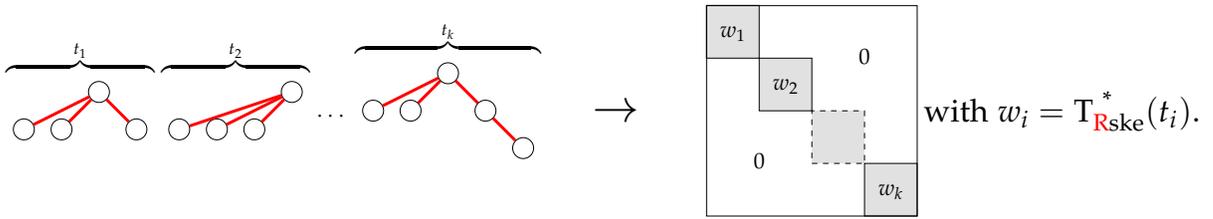
    $$\scalebox{0.7}{$\overbrace{{ \newcommand{\nodea}{\node[draw,ellipse] (a) {}
;}\newcommand{\nodeb}{\node[draw,ellipse] (b) {}
;}\newcommand{\nodec}{\node[draw,ellipse] (c) {}
;}\newcommand{\noded}{\node[draw,ellipse] (d) {}
;}\begin{tikzpicture}[baseline={([yshift=-1ex]current bounding box.center)}]
\matrix[column sep=.3cm,row sep=.3cm,ampersand replacement=\&]{
         \&         \& \nodea \& \\ 
 \nodeb  \& \nodec  \&        \& \noded  \\
};
\path[ultra thick,red] (a) edge (b) edge (c) edge (d);
\end{tikzpicture}}}^{t_1}$
      $\overbrace{\input{figures/arbres_pw/arbrew2}}^{t_2}$
      $\cdots$ $\overbrace{{ \newcommand{\nodea}{\node[draw,ellipse] (a) {}
;}\newcommand{\nodeb}{\node[draw,ellipse] (b) {}
;}\newcommand{\nodec}{\node[draw,ellipse] (c) {}
;}\newcommand{\noded}{\node[draw,ellipse] (d) {}
;}\newcommand{\nodee}{\node[draw,ellipse] (e) {}
;}\begin{tikzpicture}[baseline={([yshift=-1ex]current bounding box.center)}]
\matrix[column sep=.3cm,row sep=.3cm,ampersand replacement=\&]{
         \&         \& \nodea \& \& \\ 
 \nodeb  \& \nodec  \&        \& \noded \& \\
 \&  \&        \& \& \nodee  \\
};
\path[ultra thick,red] (a) edge (b) edge (c) edge (d)
(d) edge (e);
\end{tikzpicture}}}^{t_k}$
      \qquad\scalebox{2}{$\to$}\qquad\quad \input{figures/box_d/gd_fact}} \text{with $w_i = \Trskestar(t_i)$}.$$
  \item for any non empty red-packed tree $t=\Nodered(I,f_\ell, f_r)$, then \newline
    $\Trstar(t) = \Frstar(f_\ell)\ins\phi_I(\Frstar(f_r))$.
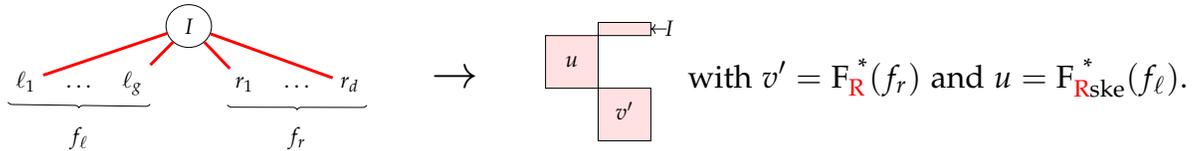
    $$\scalebox{0.7}{{ \newcommand{\nodea}{\node[draw,ellipse] (a) {$\,I\,$}
;}\newcommand{\nodeb}{\node (b) {$\ell_1$}
;}\newcommand{\nodebc}{\node (bc) {$\ldots$}
;}\newcommand{\nodec}{\node (c) {$\ell_g$}
;}\newcommand{\noded}{\node (d) {$\bullet$}
;}\newcommand{\nodee}{\node (e) {$r_1$}
;}\newcommand{\nodeef}{\node (ef) {$\ldots$}
;}\newcommand{\nodef}{\node (f) {$r_d$}
;}\begin{tikzpicture}[auto]
\matrix[column sep=.3cm, row sep=.3cm,ampersand replacement=\&]{
      \&  \&  \& \nodea  \& \&    \&    \\ 
 \nodeb  \& \nodebc \& \nodec \&  \& \nodee \& \nodeef \& \nodef  \\
};

\path[ultra thick, red] (a) edge (b) edge (c) edge (e) edge (f);
\draw[snake=brace,raise snake = 4mm,mirror snake] (b.west) -- (c.east);
\node[below of = bc] (F) {$f_\ell$};
\draw[snake=brace,raise snake = 4mm,mirror snake] (e.west) -- (f.east);
\node[below of = ef] (F) {$f_r$};

\end{tikzpicture}} \qquad\scalebox{2}{$\to$}\qquad\quad \input{figures/box_d/ins}} \text{with $v' = \Frstar(f_r)$ and $u = \Frskestar(f_\ell)$}.$$
  \end{itemize}
  There might be a problem with this definition since $\phi_I(\Frstar(f_r))$ is
  only defined if $i_p \leq |\Frstar(f_r)| + p$~(see~\cref{phi}). We
  prove in the following \cref{Fr-1welldef} that the inequality holds if
  $t \in \PTreered$.
\end{defi}

\begin{lem}\label{Fr-1welldef}
  For any red-packed forest $f$, $\Frstar(f)$ is a well defined word of size $\omega(f)$.
  For any red-packed tree $t$, $\Trstar(t)$ is a well defined word of size $\omega(t)$.
\end{lem}

\begin{proof}
  We prove by induction with the following hypothesis, for $n \in \NN$:
  \begin{equation}
    \begin{array}{rl}
      \forall f \in \PForestredinfn,
      & \Frstar(f) \text{ is well defined and } |\Frstar(f)| = \omega(f),\\
      \forall t \in \PTreeredinfn,
      & \Trstar(t) \text{ is well defined and } |\Trstar(t)| = \omega(t).
    \end{array}
    \label{H_nwelldef}
  \end{equation}
  The base case is given by the first item of
  \cref{def:Fr-1} as $\Frstar([]) = \epsilon$ and $\omega([]) = |\epsilon| = 0$.

  \noindent Let us fix $n \geq 1$ and suppose that the hypothesis
  (\ref{H_nwelldef}) holds. Let $f = [t_1, \dots, t_k] \in \PForestrednplusone$.

  \noindent $\bullet$ If $k = 1$, it is sufficient to prove the second item of
  (\ref{H_nwelldef}). Let $t = \Nodered(I, f_\ell, f_r) \in
  \PTreerednplusone$. According to \cref{red_packed_forest} with notations of
  \cref{red_weight} there are two cases:

  $\centerdot$ $d = 0$ and $I = [1, \dots, p]$. We have that $i_p = p$ and
  $|\Frstar(f_r)| = |\epsilon| = 0$ so $i_p \leq 0 + p$ and
  $\phi_I(\epsilon) = 11\dots11$ of size $p$. Now by induction on $f_\ell$, we have that
  $\Frstar(f_\ell)$ is a well defined word of size $\omega(f_\ell)$. Finally
  $\Trstar(t) = \Frstar(f_\ell) \ins \phi_I(\epsilon) = \Frstar(f_\ell) \dcdot \phi_I(\epsilon) $ is a
  well defined word of size $|\Frstar(f_\ell)| + p = \omega(t)$.

  $\centerdot$ $d \geq 1$. As $p > 0$ we can apply the hypothesis (\ref{H_nwelldef}) on
  $f_r$ and $f_\ell$. According to \cref{red_packed_forest} we have that
  \begin{align*}
    1 &\leq p + \omega(f_r) + 1 - i_p,\\
    i_p &\leq p + |\Frstar(f_r)|.
  \end{align*}
  So $\Trstar(t)$ is well defined. Moreover
  \begin{align*}
    |\Trstar(t)| &= |\Frstar(f_\ell)\ins\phi_I\Frstar(f_r)|\\
                 &= |\Frstar(f_\ell)| + p + |\Frstar(f_r)|\\
                 &= \omega(f_\ell) + p + \omega(f_r) = \omega(t).
  \end{align*}

  \noindent $\bullet$ If $k \geq 2$, the weight of trees are at least $1$ so we can apply
  (\ref{H_nwelldef}) on trees of $f$.
\end{proof}


\begin{theorem}
  \label{bij_red}
  The functions $\Fr$ and $\Frstar$ (resp. $\Tr$ and $\Trstar$) are two converse
  bijections between packed words of size $n$ and red-packed forests
  (resp. irreducible packed words and red-packed trees) of weight $n$. That is
  to say $\Frinv = \Frstar$ and $\Trinv = \Trstar$.
\end{theorem}

\begin{proof}
  The proof is very similar to the one of \cref{lem:red_skeleton_bij}. Indeed,
  we start to prove that domain and codomain are as announced (see Items (a) and
  (b) bellow), then we prove that the functions $\Fr$ and $\Frstar$ (resp. $\Tr$
  and $\Trstar$) are inverse to each other (see Items (c) and (d) bellow).

  We now give the differences with the proof of \cref{lem:red_skeleton_bij} and
  we advise the reader to read the two proofs in parallel. While for Item (a) in
  \cref{lem:red_skeleton_bij} it was simple, we need to do an induction here to
  prove that conditions on labels are respected. For Items (b), (c) and (d), the
  same inductions are done with one additional argument, so only the different
  argument of the induction is explicited here.

  (a) We prove by a mutual induction that $\Fr$ returns a red-packed forest and
  that $\Tr$ returns a red-packed tree. Indeed, we do an induction on the size
  of the word $w$. Here is our induction hypothesis for $n \in \NN$:
  \begin{equation}
    \begin{array}{rl}
      \forall v\in\PW, \text{irreducible packed word of size} \leq n,& \Tr(v) \text{ is a red-packed tree},\\
      \forall w, \text{packed word of size} \leq n,& \Fr(w) \text{ is a red-packed forest}.
    \end{array}
    \label{HR_r}
  \end{equation}
  
  The base case $(n=0)$ is given by the first item of \cref{def:construction_red_packed}.

  \noindent Now let us fix $n \geq 1$ and suppose that the hypothesis (\ref{HR_r})
  holds. Let $w \in \PW$ be a packed word of size $n+1$ and let
  $w_1 \gcdot \cdots \gcdot w_k$ be the global descent decomposition of $w$.

  \noindent $\bullet$ If $k = 1$ ($w$ is irreducible), then $\Fr(w)$ is reduced to a
  single tree $\Tr(w)$. We need to prove that $\Tr(w)$ is a red-packed tree
  (which also gives that $\Fr(w)$ is a red-packed forest). Let $w = u \ins v$ be
  the red-factorization of $w$. With $\phi_I(v') = v$, $f_\ell = \Fr(u)$ and
  $f_r = \Fr(v')$ we have that $\Tr(w) = \Nodered(I, f_\ell, f_r)$. The
  inequalities on $I$ and $v'$ in \cref{red_fact} are the same as the
  inequalities on $I$ and $f_r$ in \cref{red_packed_forest}. Therefore by
  \cref{Frweight} and (\ref{HR_r}) on $f_r$, we have that $\Tr(w)$ belongs to $\PTreered$.

  \noindent $\bullet$ If $k \geq 2$, the hypothesis (\ref{HR_r}) can be applied to each
  factors.
  
  (b) Compared to the proof of \cref{lem:red_skeleton_bij} we use the same
  general arguments to prove that $\Frstar$ and $\Frskestar$ return a packed
  word. First of all the base case and the case were the size of the forest is
  $k \geq 2$ are dealt with by a similar argumentation. It remains to prove that
  $\Trstar$ returns an irreducible packed word. We thus suppose that the
  induction hypothesis (\ref{HR_rstar}) holds for a given $n \in \NN$:
  \begin{equation}
    \begin{array}{rl}
      \forall t, \text{red-packed tree of size} \leq n,& \Trstar(t) \text{ is an irreducible packed word,} \\
      \forall f, \text{red-packed forest of size} \leq n,& \Frstar(f) \text{ is a packed word.}
    \end{array}
    \label{HR_rstar}
  \end{equation}

  \noindent Let $t = \Nodered(I, f_\ell, f_r)$ be a red-packed tree of size
  $n + 1$. By induction we have that $\Frstar(f_\ell)$ and $\Frstar(f_r)$ are
  packed words. Moreover either $f_r = \epsilon$ and $I = [1, \dots, p]$ or
  $1 \leq p + \omega(f_r) + 1 - i_p \leq \omega(r_d)$ with
  $I = [i_1, \dots, i_p]$ and $f_r = [r_1, \dots, r_d]$. In both cases,
  $\phi_I(\Frstar(f_r))$ is an red-irreducible packed word. Indeed, we recognize the
  two cases of \cref{red_irreducible} with the same inequalities. Finally
  $\Trstar(t) = \Frstar(f_\ell) \ins \phi_I(\Frstar(f_r))$ is an irreducible packed
  word according to \cref{ins_irr}.

  (c) We now want to prove that for any forest $f$ (resp. $t$),
  $\Fr(\Frstar(f))=f$ (resp. $\Tr(\Trstar(t))=t$). As in Item (b), the arguments
  are the same as is the proof of \cref{lem:red_skeleton_bij} for $\Fr$ and
  $\Frske$. In the case of $\Tr$ the new arguments are the same as in Item (b)
  (\textit{i.e.}  $\phi_I(\Frstar(f_r))$ is a red-irreducible packed word).

  (d) Finally, we want to prove that for any packed word $w$ (resp. irreducible
  packed word $v$), $\Frstar(\Fr(w))=w$ (resp. $\Trstar(\Tr(v))=v$). Once again,
  the only difference with \cref{lem:red_skeleton_bij} is the former second
  case. It remains to prove that point and we thus suppose that the induction
  hypothesis (\ref{H^{-1}_n}) holds for $n \in \NN$:
  \begin{equation}
    \begin{array}{rl}
      \forall v \in \PW, \text{irreducible packed word of size } \leq n,& \Trstar(\Tr(v)) = v,\\
      \forall w \in \PW, \text{packed word of size } \leq n,& \Frstar(\Fr(w)) = w.
    \end{array}
    \label{H^{-1}_n}
  \end{equation}

  \noindent Let $v$ be an irreducible packed word of size $n + 1$. Let
  $v = v' \ins \phi_I(v'')$ be the red-factorization of $v$. We have by
  \cref{def:construction_red_packed} that
  $\Tr(v) = \Nodered(I, \Fr(v'), \Fr(v''))$. As $|I| > 0$, the sizes of $v'$ and
  $v''$ are smaller than $n$ so we can apply (\ref{H^{-1}_n}). We have:
  \begin{align*}
    \Trstar(\Tr(v)) &= \Trstar(\Nodered(I, \Fr(v'), \Fr(v'')))\\
                    &= \Frstar(\Fr(v')) \ins \phi_I(\Frstar(\Fr(v'')))\\
                    &= v' \ins \phi_I(v'') = v. \qedhere
  \end{align*}

\end{proof}

\begin{ex}
  There is a unique forest in $\PForestredone$, namely
  \scalebox{0.5}{{ \newcommand{\nodea}{\node[draw,ellipse] (a) {$\,1\,$}
;}\begin{tikzpicture}[baseline={([yshift=-1ex]current bounding box.center)}]
\matrix[column sep=.3cm,row sep=.3cm,ampersand replacement=\&]{
 \nodea  \\
};
\path[ultra thick,red] ;
\end{tikzpicture}}}, here are the 3 forests of
  $\PForestredtwo$ with the associated packed
  word:$\Fr(12)=$\scalebox{0.5}{{ \newcommand{\nodea}{\node[draw,ellipse] (a) {$\,1\,$}
;}\newcommand{\nodeb}{\node[draw,ellipse] (b) {$\,1\,$}
;}\begin{tikzpicture}[baseline={([yshift=-1ex]current bounding box.center)}]
\matrix[column sep=.3cm,row sep=.3cm,ampersand replacement=\&]{
         \& \nodea  \\ 
 \nodeb  \&         \\
};
\path[ultra thick,red] (a) edge (b);
\end{tikzpicture}}},
  $\Fr(21)=$\scalebox{0.5}{\input{figures/arbres_pw/arbre21}},
  $\Fr(11)=$\scalebox{0.5}{{ \newcommand{\nodea}{\node[draw,ellipse] (a) {$1, 2$}
;}\begin{tikzpicture}[baseline={([yshift=-1ex]current bounding box.center)}]
\matrix[column sep=.3cm,row sep=.3cm,ampersand replacement=\&]{
 \nodea  \\
};
\path[ultra thick,red] ;
\end{tikzpicture}}}. We show below the
  13 forests of $\PForestredthree$ with the corresponding packed word:

  \noindent
  $\Fr(123)=$ \scalebox{0.5}{\input{figures/arbres_pw/arbre123}},
  $\Fr(132)=$ \scalebox{0.5}{\input{figures/arbres_pw/arbre132}},
  $\Fr(213)=$ \scalebox{0.5}{\input{figures/arbres_pw/arbre213}},
  $\Fr(231)=$ \scalebox{0.5}{\input{figures/arbres_pw/arbre231}},
  $\Fr(312)=$ \scalebox{0.5}{\input{figures/arbres_pw/arbre312}},
  $\Fr(321)=$ \scalebox{0.5}{\input{figures/arbres_pw/arbre321}},
  $\Fr(122)=$ \scalebox{0.5}{{ \newcommand{\nodea}{\node[draw,ellipse] (a) {$1, 2$}
;}\newcommand{\nodeb}{\node[draw,ellipse] (b) {$\,1\,$}
;}\begin{tikzpicture}[baseline={([yshift=-1ex]current bounding box.center)}]
\matrix[column sep=.3cm,row sep=.3cm,ampersand replacement=\&]{
         \& \nodea  \\ 
 \nodeb  \&         \\
};
\path[ultra thick,red] (a) edge (b);
\end{tikzpicture}}},
  $\Fr(212)=$ \scalebox{0.5}{{ \newcommand{\nodea}{\node[draw,ellipse] (a) {$1, 3$}
;}\newcommand{\nodeb}{\node[draw,ellipse] (b) {$\,1\,$}
;}\begin{tikzpicture}[baseline={([yshift=-1ex]current bounding box.center)}]
\matrix[column sep=.3cm,row sep=.3cm,ampersand replacement=\&]{
         \& \nodea  \&         \\ 
         \&         \& \nodeb  \\
};
\path[ultra thick,red] (a) edge (b);
\end{tikzpicture}}},
  $\Fr(221)=$ \scalebox{0.5}{\input{figures/arbres_pw/arbre221}},
  $\Fr(112)=$ \scalebox{0.5}{{ \newcommand{\nodea}{\node[draw,ellipse] (a) {$\,1\,$}
;}\newcommand{\nodeb}{\node[draw,ellipse] (b) {$1, 2$}
;}\begin{tikzpicture}[baseline={([yshift=-1ex]current bounding box.center)}]
\matrix[column sep=.3cm,row sep=.3cm,ampersand replacement=\&]{
         \& \nodea  \\ 
 \nodeb  \&         \\
};
\path[ultra thick,red] (a) edge (b);
\end{tikzpicture}}},
  $\Fr(121)=$ \scalebox{0.5}{{ \newcommand{\nodea}{\node[draw,ellipse] (a) {$\,2\,$}
;}\newcommand{\nodeb}{\node[draw,ellipse] (b) {$1, 2$}
;}\begin{tikzpicture}[baseline={([yshift=-1ex]current bounding box.center)}]
\matrix[column sep=.3cm,row sep=.3cm,ampersand replacement=\&]{
         \& \nodea  \&         \\ 
         \&         \& \nodeb  \\
};
\path[ultra thick,red] (a) edge (b);
\end{tikzpicture}}},
  $\Fr(211)=$ \scalebox{0.5}{\input{figures/arbres_pw/arbre211}},
  $\Fr(111)=$ \scalebox{0.5}{{ \newcommand{\nodea}{\node[draw,ellipse] (a) {$1, 2, 3$}
;}\begin{tikzpicture}[baseline={([yshift=-1ex]current bounding box.center)}]
\matrix[column sep=.3cm,row sep=.3cm,ampersand replacement=\&]{
 \nodea  \\
};
\path[ultra thick,red] ;
\end{tikzpicture}}}.
  
  More examples can be found in the annexes section with
  \cref{F(w)123,F(1234),F(1233),F(1223),F(1123),F(1111)}.
\end{ex}

We conclude by the main theorem of this subsection. It is a generalization of
the construction of \cite{Foissy_2011} for $\FQSym$ and permutations to $\WQSym$
and packed words. Indeed, if we restrict the construction on permutations and we
consider right children of a node as label of this node, we have the same
construction as in \cite{Foissy_2011} with a shift of 1 for labels. Here are
some examples of trees in \cite{Foissy_2011} and the equivalent red-packed tree.

\begin{table}[h]
  \noindent\scalebox{0.9}{\parbox{5cm}{
\begin{center}
  \begin{tabular}{cc}
    \begin{tabular}{|c|c|}
      \hline
      $\tdun{T}$ with $T = (\ttroisun, 1)$ & \scalebox{0.5}{\input{figures/arbres_pw/arbre2413}}\\
      \hline
      $\tdun{T}$ with $T = (\ttroisun, 2)$ & \scalebox{0.5}{\input{figures/arbres_pw/arbre2143}}\\
      \hline
      $\tdun{T}$ with $T = (\ttroisdeux, 1)$ & \scalebox{0.5}{\input{figures/arbres_pw/arbre1423}}\\
      \hline
      $\tdun{T}$ with $T = (\ttroisdeux, 2)$ & \scalebox{0.5}{\input{figures/arbres_pw/arbre1243}}\\
      \hline
    \end{tabular}
  &
    \begin{tabular}{|c|c|}
      \hline
      $\tdun{T}$ with $T = (\tdun{T'}, 1)$ with $T' = (\tdeux, 1)$ & \scalebox{0.5}{\input{figures/arbres_pw/arbre1432}}\\
      \hline
      $\tdun{T}$ with $T = (\tdun{T'}, 2)$ with $T' = (\tdeux, 1)$ & \scalebox{0.5}{\input{figures/arbres_pw/arbre1342}}\\
      \hline
      $\tddeux{}{T}$ with $T = (\tdeux, 1)$ & \scalebox{0.5}{\input{figures/arbres_pw/arbre3142}}\\
      \hline
      $\tddeux{T}{}$ with $T = (\tdeux, 1)$ & \scalebox{0.5}{\input{figures/arbres_pw/arbre1324}}\\
      \hline
    \end{tabular}
  \end{tabular}
\end{center}
}}
\caption{Equivalence between trees of \cite{Foissy_2011} and red-packed trees.}
\end{table}

\bigskip

All the constructions with red-packed forests have been done in order to have
this theorem.

\begin{theorem}\label{thm:bij_primal}
  For all $n\in\NN$ we have the three following equalities :
  \begin{equation*}
    \dim(\WQSym^*_n) = \#\PForestredn
    \quad\text{and}\quad
    \dim(\Prim^*_n) = \#\PTreeredn
    \quad\text{and}\quad
    \dim(\TPrim^*_n) = \#\PParticularredn
  \end{equation*}
\end{theorem}

\begin{proof}
  \cref{bij_red} proves the first equality. It also gives a relation between
  $\#\PForestred$ and $\#\PTreered$. Indeed a red-packed forest of weight $n$ is
  an ordered sequence of red-packed trees of weight $(n_k)$ such that $\sum_k(n_k)=n$. This
  relation is the same between $\dim(\WQSym^*_n)$ and $\dim(\Prim^*_n)$ according to
  \cref{left_prod} (\textit{i.e.} $\SA=\SP/(1-\SP)$).

  \noindent Red-skeleton trees are equivalent to ordered trees decorated by
  red-irreducible words as said in \cref{red_ske_flat}. Recall that a basis of
  primitive elements is given by~\cref{brace} as ordered trees decorated by
  totally primitive elements. Elements of $\PParticularred$ are by definition in
  bijections with red-irreducible words, labels of red-skeleton trees.
\end{proof}


\subsection{Primal (\blue{Blue})}\label{Dual_Blue}

Now we do the same work for the primal side: $\WQSym$. This subsection follows
the same structure of statements as the previous one. Recall that in
\cref{Primal_Red} we constructed a bijection between packed words and
red-forests by recursively decomposing packed words using global descent and
removal of maximums. In this section we follow the same path: incerting the last
letter using $\psi_{i^\alpha}$ (lowercase $i$ designates the integer value) instead of
new maximums using $\phi_I$ (uppercase $I$ designates the list of their
positions). We define a blue-factorization of packed words. When used
recursively, blue-factorization and global descent decomposition construct a
bijection between PW and so-called blue-packed forest.  Since the general
structure of proofs are the same as in the previous section, we will mostly
focus on the differences between combinatorials arguments.

\subsubsection{Decomposition of packed words through last letter}

In this section, we define two combinatorial operations on packed words
($\psi_{i^\alpha}$ and $\insl$) and the blue-factorization that use them. The unary
operation $\psi_{i^\alpha}$ insert the new value $i$ at the end of a given word. A word
that cannot be factorized $u \insl v$ in a non trivial way is called
blue-irreducible. Blue-irreducible words will index a new basis of $\TPrim$.

\begin{defi}\label{psi}
  Fix $n \in \NN$ and $w \in \PW_n$. For any $1 \leq i \le \max(w) + 1$ (with the
  convention $\max(\epsilon) = 0$), we denote by
  $\psi_{i^\circ}(w) = u_1 \cdots u_n \cdot i$ the packed word defined by
  $u_k = w_k$ if $w_k < i$ and $u_k = w_k+1$ otherwise. We also define
  $\psi_{i^\bullet}(w) = w \cdot i$ for any $1 \le i \le \max(w)$.
\end{defi}

\begin{ex}
  $\psi_{\blue{2}^\circ}(1232) = 1343\blue{2}$, $\psi_{\blue{2}^\bullet}(1232) = 1232\blue{2}$,
  $\psi_{\blue{4}^\circ}(1232) = 1232\blue{4}$ and $\psi_{\blue{1}^\circ}(\epsilon) = \blue{1}$.
\end{ex}



\begin{lem} \label{psi_bij} For any $W\in\PW_\ell$ where~$\ell>0$ there exists a unique
  triplet $(i, \alpha, w)$ where ${i \in [1\dots \ell+1]}$,
  $\alpha \in \{\circ, \bullet\}$ and $w$ is a packed word, such that
  $W = \psi_{i^\alpha}(w)$.
\end{lem}
Depending on $\alpha$, the box diagram can be represented as \newline
$W = \input{figures/box_d/psi_circ}$ or $W =
\input{figures/box_d/psi_bullet}$. In the general case, we will note
$W = \input{figures/box_d/psi}$.

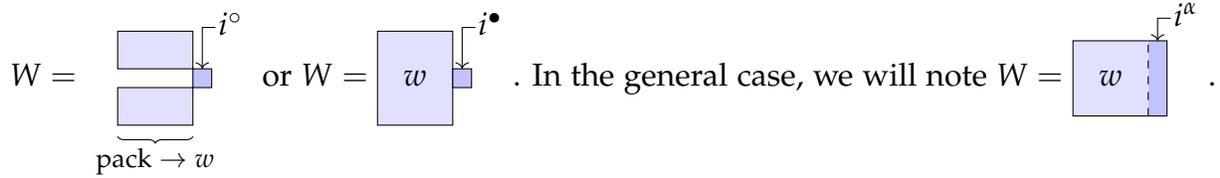
\begin{proof}
  Let $W\in\PW_\ell$ with~$\ell>0$ and $i$ the value of the last letter of $W$.
  \begin{itemize}
  \item If $i$ appears multiple times in $W$, then let $w = W_1\dots W_{\ell-1}$, we
    only remove the last letter $i$ of $W$. We have $W = \psi_{i^\bullet}(w)$.
  \item Otherwise, $i$ appears only as the last letter, then let
    $w = \pack(W_1\dots W_{\ell-1})$, we remove the last letter $i$ of $W$ and pack
    the word. We have $W = \psi_{i^\circ}(w)$.
  \end{itemize}
  If $\psi_{i^\alpha}(u) = \psi_{j^\beta}(v)$ then the last letter is the same so
  $i = j$, the multiplicity of this letter is the same so $\alpha = \beta$ and the prefix
  are the same $u = v$.
\end{proof}

\begin{defi}\label{insl}
  Let $u, v \in \PW$ with $v \neq \epsilon$. By~\cref{psi_bij}, there is a unique triplet
  $(i, \alpha, v')$ such that $v=\psi_{i^\alpha}(v')$. Let $i' = i + \max(u)$, we define
  $u\insl v \eqdef \psi_{i'^\alpha}(v'\gcdot u)$. In other words, we remove the last
  letter of the right word, perform a reversed left shifted concatenation and
  adding back the last letter also shifted.
\end{defi}

\begin{ex}
  $2123\insl \blue{31231}\red{2} = 2123 \insl \psi_{(2)^\bullet}(\blue{31231}) =
  \psi_{(2+3)^\bullet}(\blue{64564}2123) = \blue{64564}2123\red{5}$.
\end{ex}
\begin{figure}[!h]
$$u\insl v = \scalebox{1}{\input{figures/box_d/insl}} \qquad
\scalebox{0.4}{\input{figures/box_d/2123}}~\insl~\scalebox{0.4}{\input{figures/box_d/312312}}~=~\scalebox{0.4}{\input{figures/box_d/pw_insl_ex}}$$\vspace{-1em}
  \caption{Box digrams: the operation $\insl$}
  \label{boxd_insl}
\end{figure}

\begin{lem}\label{blue_fact}
  Let~$w$ be an irreducible packed word. There exists a unique factorization of
  the form $w=u\insl v$ which maximizes the size of $u$. In this factorization,
  let $v'$ and $i^\alpha$ such that $v = \psi_{i^\alpha}(v')$,
  \begin{itemize}
  \item either $v' = \epsilon$ and $i^\alpha = 1^\circ$,
  \item or $v'$ is irreducible and $1 \le i \leq \max(v')$.
  \end{itemize}
  We call it the \textbf{blue-factorization} of a word.
\end{lem}

\begin{ex}
  Here is a first detailled example of a blue-factorization of an irreducible
  packed word:

  Consider the irreducible packed word $w = 654623314$.
  \begin{itemize}
  \item The first step is to remove the last letter $i=4$. Here there are
    multiple occurences of the last letter in $w$, then $\alpha = \bullet$, we get
    $w' = 65462331$ which is a packed word, but is not irreducible.
  \item The second step is to set $w'_1$ as the first irreducible factor of
    $w'$ and $u$ the rest of $w'$. This way $w'=w'_1 \gcdot u$ and the size of
    $u$ is maximized. Here $w'_1 = 3213$ and $u = 2331$. Let
    $i' = i - \max(u) = 4 - 3 = 1$.
  \item Finally, we get the following decomposition of $w$ (see \cref{psi} for
    $\psi$ and \cref{insl} for $\insl$): \\
    $w = 654623314 = u \insl \psi_{i'^\alpha}(w'_1) = (2331) \insl \psi_{1^\bullet}(3213) = 2331
    \insl 32131$.
  \end{itemize}
\end{ex}

\begin{ex} Here are some other blue-factorizations:
  \begin{alignat*}{2}
    234313 &= 1\insl \psi_{2^\bullet}(1232) = 1 \insl 12322 &\qquad
    245413 &= 1\insl \psi_{2^\circ}(1232) = 1 \insl 13432 \\
    11 &= \epsilon \insl \psi_{1^\bullet}(1) = \epsilon \insl 11 &\qquad
    112 &= 11 \insl \psi_{1^\circ}(\epsilon) = 11 \insl 1
  \end{alignat*}
\end{ex}

\begin{proof}
  Let $w$ be irreducible and let $(i, \alpha, w')$ be the unique triplet such that
  $w=\psi_{i^\alpha}(w')$ according to~\cref{psi_bij}.
  
  If $i = \max(w)$ and it appears only one time (\textit{i.e.} $\alpha = \circ$ and
  $i = \max(w') + 1$) then the blue-factorizations is $w = w' \insl \psi_{1^\circ}(\epsilon)$.

  In any other case, we write $w' = w_1'\gcdot w_2' \gcdot\dots\gcdot w_k'$, the
  decomposition into irreducibles. Let $u = w_2' \gcdot\dots\gcdot w_k'$ and
  $i' = i - \max(u)$. We have that $i' \leq \max(w_1')$ otherwise $w$ wouldn't be
  packed and $1 \leq i'$ otherwise $w$ wouldn't be irreducible. If
  $\alpha = \circ$ then $1 \neq i'$ otherwise $w$ wouldn't be irreducible. Then we have
  $w = u \insl \psi_{i'^\alpha}(w_1')$ where the size of $u$ is maximized.
\end{proof}

\begin{rem}\label{rem:inv_permu}
  When restricted to permutations, blue-factorization is equal to a
  red-factorization applies to the inverse. Let $\sigma$ be a permutation and
  $\sigma = \mu \ins \nu$ be the red-factorization of $\sigma$, then
  $\sigma^{-1}=\mu^{-1} \insl \nu^{-1}$ is the blue-factorization of $\sigma^{-1}$.
\end{rem}

\begin{defi} \label{blue_irreducible}
  A packed word $w$ is \textbf{blue-irreducible} if $w$ is irreducible and 
  $w = u \insl v$ implies that $u = \epsilon (\text{and } w = v)$.
\end{defi}

Here are some useful lemmas on the operation $\insl$. There are some
similarities with
\cref{ins_pseudo_assoc_left,ins_pseudo_assoc_right,ins_irr,red_fact_irr}.

\begin{lem} \label{insl_pseudo_assoc_left}
  For any $u, v, w \in \PW$ with $w \neq \epsilon$, we have $u \insl (v \insl w) = (v \gcdot u) \insl w$.
  $$\scalebox{1}{\input{figures/box_d/insl_of_insl} \qquad\scalebox{2}{=}\qquad\quad \input{figures/box_d/gd_in_insl}}$$
\end{lem}

\begin{proof}
  Let $u, v, w \in \PW$ with $w \neq \epsilon$, and let $w'$ and $i^\alpha$ such that $w = \psi_{i^\alpha}(w')$.
  \begin{align*}
    u \insl (v \insl w) &= u \insl (\psi_{(i+\max(v))^\alpha}(w' \gcdot v))\\
                        &= \psi_{(i+\max(v)+\max(u))^\alpha}((w' \gcdot v) \gcdot u)\\
                        &= \psi_{(i+\max(v)+\max(u))^\alpha}(w' \gcdot (v \gcdot u))\\
                        &= (v \gcdot u) \insl \psi_{i^\alpha}(w')\\
                        &= (v \gcdot u) \insl w. \qedhere
  \end{align*}
\end{proof}

\begin{lem} \label{insl_pseudo_assoc_right}
  For any $u, v, w \in \PW$ with $v \neq \epsilon$, we have $u \insl (w \gcdot v) = w \gcdot (u \insl v)$.
  $$\scalebox{1}{\input{figures/box_d/insl_of_gd} \qquad\scalebox{2}{=}\qquad\quad \input{figures/box_d/insl_in_gd}}$$
\end{lem}

\begin{proof}
  Let $u, v, w \in \PW$ with $v \neq \epsilon$ and let $v'$ and $i^\alpha$ such that $v = \psi_{i^\alpha}(v')$.
  \begin{align*}
    u \insl (w \gcdot v) &= u \insl (w \gcdot \psi_{i^\alpha}(v'))\\
                         &= u \insl \psi_{i^\alpha}(w \gcdot v')\\
                         &=  \psi_{(i+\max(u))^\alpha}((w \gcdot v') \gcdot u)\\
                         &=  \psi_{(i+\max(u))^\alpha}(w \gcdot (v' \gcdot u))\\
                         &=  w \gcdot \psi_{(i+\max(u))^\alpha}(v' \gcdot u)\\
                         &= w \gcdot (u \insl v). \qedhere
  \end{align*}
\end{proof}

\begin{rem}\label{rem:skew_dup_blue}
  Theses relations are the same up to symmetry as the one with $\ins$
  (\cref{ins_pseudo_assoc_left,ins_pseudo_assoc_right}). So adding the
  associativity of shifted concatenation
  $u \gcdot (v \gcdot w) = (u \gcdot v) \gcdot w$, the two operations $\insl$
  and $\gcdot$ verify relations of the \textit{skew-duplicial operad}
  \cite{BurDel_dup}.
\end{rem}

\begin{coro}\label{insl_irr}
  For any $u, v \in \PW$, we have that $u \insl v$ is irreducible if and only if
  $v$ is irreducible.
\end{coro}

\begin{proof}
  By contradiction, if $v = v_1 \gcdot v_2$ then by
  \cref{insl_pseudo_assoc_right} $u \insl v = v_1 \gcdot (u \insl v_2)$. Now if
  $u \insl v = w_1 \gcdot w_2$ as the value of the last letter of $u \insl v$ is
  greater than $\max(u)$ we have that $w_2 = w_2' \cdot w_2'' \cdot i$ such that
  $\pack(w_2'') = u$. We also have that $w_1 \gcdot \pack(w_2' \cdot i) = v$.
\end{proof}



\begin{prop}\label{blue_fact_irr}
  For any word $w$, $w = u \insl v$ is the blue-factorization of $w$ if and only
  if $v$ is blue-irreducible.
\end{prop}

\begin{proof}
  Let $w \in \PW$ and let $u \insl v$ be the blue-factorization of $w$. Let $v_1$ and
  $v_2$ such that $v = v_1 \insl v_2$, then $(v_1 \gcdot u) \insl v_2 = w$ by
  \cref{insl_pseudo_assoc_left}, but in the blue-factorization the size of $u$ is maximized
  so $|(v_1 \gcdot u)| \leq |u|$ and then we have that $v_1 = \epsilon$ so $v$ is
  blue-irreducible.

  Let $w \in \PW$ and let $u$ and $v$ such that $w = u \insl v$ and $v$ is
  blue-irreducible. By contradiction, suppose that there exists $u', v'$ such
  that $w = u' \insl v'$ with $|u| < |u'|$ and $v' \neq \epsilon$. Then necessarily
  $u$ is a suffix of $u'$. Let $u''$ such that $u' = u'' \cdot u$, then
  $pack(u'') \insl v' = v$. But $v$ is blue-irreducible. So the size of $u$ is
  maximal if $v$ is blue-irreducible.
\end{proof}

Thanks \cref{rem:inv_permu} the following proposition is immediate.

\begin{prop}
  A permutation $\sigma$ is blue-irreducible if and only if $\sigma^{-1}$ is red-irreducible.
\end{prop}

\subsubsection{Blue-forests from decomposed packed words using~\texorpdfstring{$\psi$}{psi}}
\label{blue_forest}

As in \cref{red_forests} we will apply recursively the blue-factorization of the
former section to construct a bijection between packed words and a certain kind
of labeled biplane trees.

In this construction, the labels can be a blue-irreducible word for skeleton,
or an integer with a sign $\alpha \in \{\circ, \bullet\}$. In order to differentiate the trees
from the one of the previous section, we will draw them in \blue{blue}. As
before, for a labeled biplane tree, we denote the trees by $\Nodeblue(x, f_\ell, f_r)$.

\begin{ex}
  $\Nodeblue(1^\circ,\ [],[]) =$ \scalebox{0.5}{{ \newcommand{\nodea}{\node[draw,ellipse] (a) {$1^\circ$}
;}\begin{tikzpicture}[auto]
\matrix[column sep=.3cm, row sep=.3cm,ampersand replacement=\&]{
 \nodea  \\
};
\end{tikzpicture}}}, and
  $\Nodeblue(1^\bullet,[\Nodeblue(1^\circ,[],[])],[\Nodeblue(1^\circ,[],[])]) =$ \scalebox{0.5}{\input{figures/arbres_pw/arbre212left_lbl}}.
\end{ex}

We apply recursively the global descent decomposition and the blue-factorization
of~\cref{blue_fact}. We obtain an algorithm which takes a packed word and
returns a biplane forest where nodes are decorated by blue-irreducible words:
\begin{defi}
  \label{def:construction_blue_skeleton}
  Exactly as \cref{def:construction_red_skeleton} of $\Frske$ and $\Trske$, we
  now define two functions $\Fbske$ and $\Tbske$. These functions transform
  respectively a packed word and an irreducible packed word into
  respectively a biplane forest and a biplane tree. These functions are defined
  in a mutual recursive way as follow:
  \begin{itemize}
  \item $\Fbske(\epsilon) = []$ (empty forest),
  \item for any packed word $w$, let $w_1\gcdot w_2\gcdot\dots\gcdot w_k$ be the
    global descent decomposition of $w$, then
    $\Fbske(w) \eqdef [\Tbske(w_k), \Tbske(w_{k-1}), \dots, \Tbske(w_1)]$
    (notice the inversion compared to~\cref{def:construction_red_skeleton}).
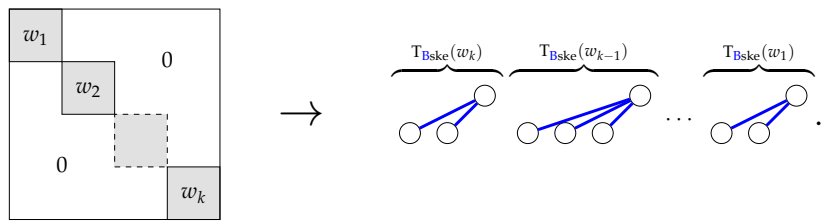
    $$\scalebox{0.7}{\input{figures/box_d/gd_fact} \qquad\scalebox{2}{$\to$}\qquad\quad $\overbrace{\input{figures/arbres_pw/arbrew1_ske_left}}^{\Tbske(w_k)}$
      $\overbrace{\input{figures/arbres_pw/arbrew2_ske_left}}^{\Tbske(w_{k-1})}$
      $\cdots$ $\overbrace{\input{figures/arbres_pw/arbrew1_ske_left}}^{\Tbske(w_1)}$}.$$
  \item for any irreducible packed word $w$, we define
    $\Tbske(w) \eqdef \Nodeblue(v, \Fbske(u), [])$ where
    $w=u \insl v$ is the blue-factorization of $w$.
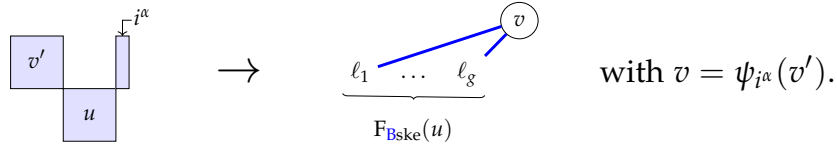
    $$\scalebox{0.7}{\input{figures/box_d/insl} \qquad\scalebox{2}{$\to$}\qquad\quad { \newcommand{\nodea}{\node[draw,ellipse] (a) {$v$}
;}\newcommand{\nodeb}{\node (b) {$\ell_1$}
;}\newcommand{\nodebc}{\node (bc) {$\ldots$}
;}\newcommand{\nodec}{\node (c) {$\ell_g$}
;}\begin{tikzpicture}[auto]
\matrix[column sep=.3cm, row sep=.3cm,ampersand replacement=\&]{
      \&  \&  \& \nodea  \& \&    \&    \\ 
 \nodeb  \& \nodebc \& \nodec \&  \& \& \& \\ 
};

\path[ultra thick, blue] (a) edge (b) edge (c); 
\draw[snake=brace,raise snake = 4mm,mirror snake] (b.west) -- (c.east);
\node[below of = bc] (F) {$\Fbske(u)$};

\end{tikzpicture}}} \text{with $v = \psi_{i^\alpha}(v')$}.$$
  \end{itemize}
\end{defi}


\begin{ex}\label{ex:blue_skeleton}
  Let $w = 8967647523314$, here is the global descent decomposition $w = w_1\gcdot w_2$
  with $w_1 = 12$ and $w_1 = 67647523314$. Now, we have the blue-factorization
  of $w_1$ and $w_2$ using~\cref{blue_fact} as
  \[
    w_1 = 1 \insl \psi_{1^\circ}(\epsilon).  \quad\text{and}\quad
    w_2 = 2331 \insl \psi_{1^\bullet}(343142) = (122
    \gcdot 1) \insl 3431421,
  \]
  $$\scalebox{0.3}{\input{figures/box_d/pw_blue_ex}}$$
  It gives the following forest:
  \begin{align*}
    \Fbske(8967647523314) & = [\Tbske(12),\ \Tbske(67647523314)] \\
                           & =
                             \scalebox{0.7}{{ \newcommand{\nodea}{\node[draw,ellipse] (a) {$3431421$}
;}\newcommand{\nodeb}{\node (b) {$\Fbske(2331)$}
;}
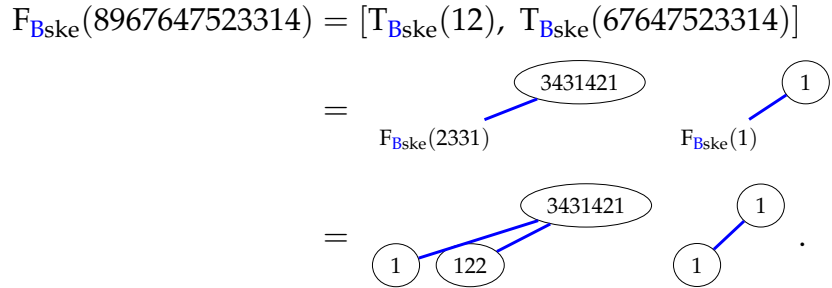
\begin{tikzpicture}[baseline={([yshift=-1ex]current bounding box.center)}]
\matrix[column sep=.3cm, row sep=.3cm,ampersand replacement=\&]{
         \& \nodea   \\
 \nodeb  \&          \\
};

\path[ultra thick, blue] (a) edge (b);
\end{tikzpicture}}
}
                             \scalebox{0.7}{{ \newcommand{\nodea}{\node[draw,ellipse] (a) {$\,1\,$}
;}\newcommand{\nodeb}{\node (b) {$\Fbske(1)$}
;}\begin{tikzpicture}[baseline={([yshift=-1ex]current bounding box.center)}]
\matrix[column sep=.3cm, row sep=.3cm,ampersand replacement=\&]{
         \& \nodea   \\
 \nodeb  \&          \\
};

\path[ultra thick, blue] (a) edge (b);
\end{tikzpicture}}
}\\ 
                           & =
                             \scalebox{0.7}{\input{figures/arbres_pw/arbre8967647523314leftske}}.
  \end{align*}
\end{ex}

\begin{defi} \label{def:blue_skeleton}
  A labeled biplane forest (resp. tree) is a \textbf{blue-skeleton forest}
  (resp. \textbf{tree}) if and only if it is labeled by blue-irreducible words
  and no node has a right child.
\end{defi}

We want to prove that the functions $\Fbske$ and $\Tbske$ are bijections. To do
that, as for $\Frske$ and $\Trske$, we first define the two inverse functions.

\begin{defi}\label{def:construction_blue_skeleton_star}
  We now define two functions $\Fbskestar$ and $\Tbskestar$ that transform
  respectively a blue-skeleton forest and tree into packed words. These
  functions are defined in a mutual recursive way as follow:
  \begin{itemize}
  \item $\Fbskestar([]) = \epsilon$,
  \item for any blue-skeleton forest $f = [t_1, \dots, t_k]$, we define \newline
    ${\Fbskestar(f) \eqdef \Tbskestar(t_k)\gcdot \dots \gcdot\Tbskestar(t_1)}$. \newline
    (notice the inversion compared to~\cref{def:construction_red_skeleton_star})
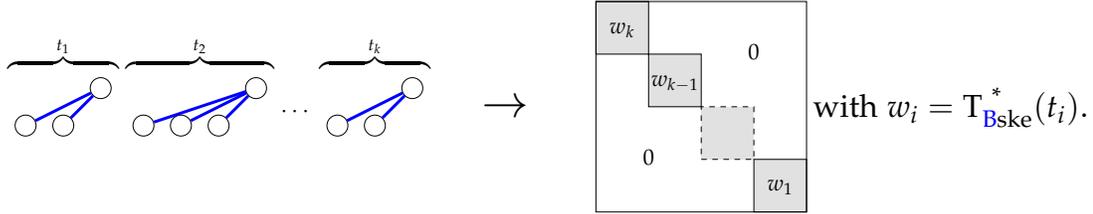
    $$\scalebox{0.7}{$\overbrace{\input{figures/arbres_pw/arbrew1_ske_left}}^{t_1}$
      $\overbrace{\input{figures/arbres_pw/arbrew2_ske_left}}^{t_2}$
      $\cdots$ $\overbrace{\input{figures/arbres_pw/arbrew1_ske_left}}^{t_k}$
      \qquad\scalebox{2}{$\to$}\qquad\quad \input{figures/box_d/gd_fact_left}} \text{with $w_i = \Tbskestar(t_i)$}.$$
  \item for any blue-skeleton tree $t = \Nodeblue(v, f_r, [])$, we define \newline
    $\Tbskestar(t) \eqdef \Fbskestar(f_r) \insl v$.
    $$\scalebox{0.7}{{ \newcommand{\nodea}{\node[draw,ellipse] (a) {$v$}
;}\newcommand{\nodeb}{\node (b) {$\ell_1$}
;}\newcommand{\nodebc}{\node (bc) {$\ldots$}
;}\newcommand{\nodec}{\node (c) {$\ell_g$}
;}\begin{tikzpicture}[auto]
\matrix[column sep=.3cm, row sep=.3cm,ampersand replacement=\&]{
      \&  \&  \& \nodea  \& \&    \&    \\ 
 \nodeb  \& \nodebc \& \nodec \&  \& \& \& \\ 
};

\path[ultra thick, blue] (a) edge (b) edge (c); 
\draw[snake=brace,raise snake = 4mm,mirror snake] (b.west) -- (c.east);
\node[below of = bc] (F) {$f_\ell$};

\end{tikzpicture}} \qquad\scalebox{2}{$\to$}\qquad\quad \input{figures/box_d/insl}} \text{with $v = \psi_{i^\alpha}(v')$ and $u = \Fbskestar(f_\ell)$}.$$
  \end{itemize}
\end{defi}

\begin{lem}\label{lem:blue_skeleton_bij}
  The functions $\Fbske$ and $\Fbskestar$ (resp. $\Tbske$ and $\Tbskestar$) are
  two converse bijections between packed words and blue-skeleton forests
  (resp. irreducible packed words and blue-skeleton trees). That is to
  say $\Fbskeinv = \Fbskestar$ and $\Tbskeinv = \Tbskestar$.
\end{lem}

\begin{proof}
  The proof structure is the same as the one of
  \cref{lem:red_skeleton_bij} with use of statements comming from this
  subsection. We can see in this table some of the main statements that are
  exchanged for this dual part:
  \begin{center}
    \begin{tabular}{|l|l|l|}
      \hline
      \cref{red_fact} & \cref{blue_fact} & red-factorization and blue-factorization.\\
      \hline
      \cref{red_irreducible} & \cref{blue_irreducible} & red-irreducible words and blue-irreducible words.\\
      \hline
      \cref{ins_irr} &                 & $u \ins v$ irreducible $\iff$ $v$ irreducible \\
                     & \cref{insl_irr} & $u \insl v$ irreducible $\iff$ $v$ irreducible.\\
      \hline
      \cref{red_fact_irr} &  & $u \ins v$ red-factorization $\iff$ $v$ red-irreducible \\
      & \cref{blue_fact_irr} & $u \insl v$ blue-factorization $\iff$ $v$ blue-irreducible.\\
      \hline
      \cref{def:construction_red_skeleton} & \cref{def:construction_blue_skeleton} & $\Frske$, $\Trske$ and $\Fbske$, $\Tbske$.\\
      \hline
      \cref{def:red_skeleton} & \cref{def:blue_skeleton} & red-skeleton forest and blue-skeleton forest.\\
      \hline
      \cref{def:construction_red_skeleton_star} & \cref{def:construction_blue_skeleton_star} & $\Frskestar$, $\Trskestar$ and $\Fbskestar$, $\Tbskestar$.\\
      \hline
    \end{tabular}\qedhere
  \end{center}
\end{proof}

Now that we have the blue-skeleton, we will add right forests to every nodes to
have biplane trees. For every nodes, if $v$ is the blue-irreducible word in the
label, then with \cref{psi_bij} $v = \psi_{i^\alpha}(v')$, $i^\alpha$ is the new label and
$\Fb(v')$ is the new right forest.

Here is the formal definition of $\Fb(w)$ and $\Tb(w)$ which is very similar to
\cref{def:construction_blue_skeleton}, only the third item is different. The
labels are now pairs of an integer and a sign $\alpha \in \{\circ, \bullet\}$.

\begin{defi} \label{def:construction_blue_packed}
  The forest~$\Fb(w)$~(resp. tree~$\Tb(w)$) associated to a packed word (resp.
  irreducible packed word) $w$ are defined in a mutual recursive way
  as follows:
  \begin{itemize}
  \item $\Fb(\epsilon) = []$ (empty forest),
  \item for any packed word $w$, let $w_1\gcdot w_2\gcdot\dots\gcdot w_k$ be the
    global descent decomposition of $w$, then
    $\Fb(w) \eqdef [\Tb(w_k), \Tb(w_{k-1}), \dots, \Tb(w_1)]$ (notice the
    inversion compared to~\cref{def:construction_red_packed}).
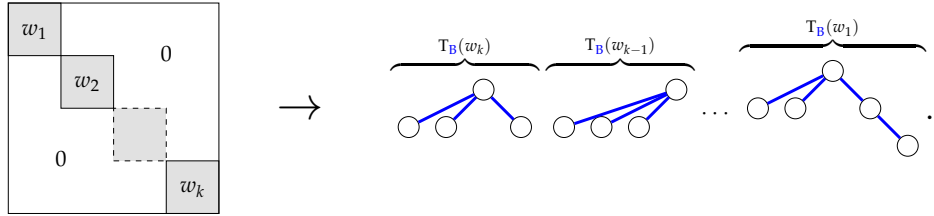
    $$\scalebox{0.7}{\input{figures/box_d/gd_fact} \qquad\scalebox{2}{$\to$}\qquad\quad $\overbrace{{ \newcommand{\nodea}{\node[draw,ellipse] (a) {}
;}\newcommand{\nodeb}{\node[draw,ellipse] (b) {}
;}\newcommand{\nodec}{\node[draw,ellipse] (c) {}
;}\newcommand{\noded}{\node[draw,ellipse] (d) {}
;}\begin{tikzpicture}[baseline={([yshift=-1ex]current bounding box.center)}]
\matrix[column sep=.3cm,row sep=.3cm,ampersand replacement=\&]{
         \&         \& \nodea \& \\ 
 \nodeb  \& \nodec  \&        \& \noded  \\
};
\path[ultra thick,blue] (a) edge (b) edge (c) edge (d);
\end{tikzpicture}}}^{\Tb(w_k)}$
      $\overbrace{\input{figures/arbres_pw/arbrew2_left}}^{\Tb(w_{k-1})}$
      $\cdots$ $\overbrace{{ \newcommand{\nodea}{\node[draw,ellipse] (a) {}
;}\newcommand{\nodeb}{\node[draw,ellipse] (b) {}
;}\newcommand{\nodec}{\node[draw,ellipse] (c) {}
;}\newcommand{\noded}{\node[draw,ellipse] (d) {}
;}\newcommand{\nodee}{\node[draw,ellipse] (e) {}
;}\begin{tikzpicture}[baseline={([yshift=-1ex]current bounding box.center)}]
\matrix[column sep=.3cm,row sep=.3cm,ampersand replacement=\&]{
         \&         \& \nodea \& \& \\ 
 \nodeb  \& \nodec  \&        \& \noded \& \\
 \&  \&        \& \& \nodee  \\
};
\path[ultra thick,blue] (a) edge (b) edge (c) edge (d)
(d) edge (e);
\end{tikzpicture}}}^{\Tb(w_1)}$}.$$
  \item for any irreducible packed word $w$, we define
    $\Tb(w) \eqdef \Nodeblue(i^\alpha, \Fb(u), \Fb(v'))$ where
    $w=u\insl\psi_{i^\alpha}(v')$ is the blue-factorization of $w$.
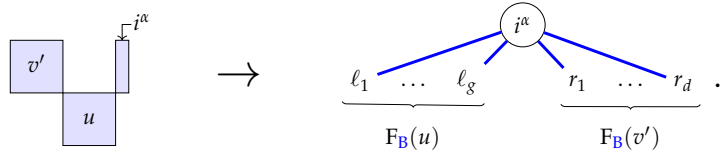
    $$\scalebox{0.7}{\input{figures/box_d/insl} \qquad\scalebox{2}{$\to$}\qquad\quad \input{figures/arbres_pw/arbre_left_F}}.$$
  \end{itemize}
\end{defi}

\begin{ex}
  Consider again $w = 8967647523314$. We start from the blue-skeleton forest from
  \cref{ex:blue_skeleton}.
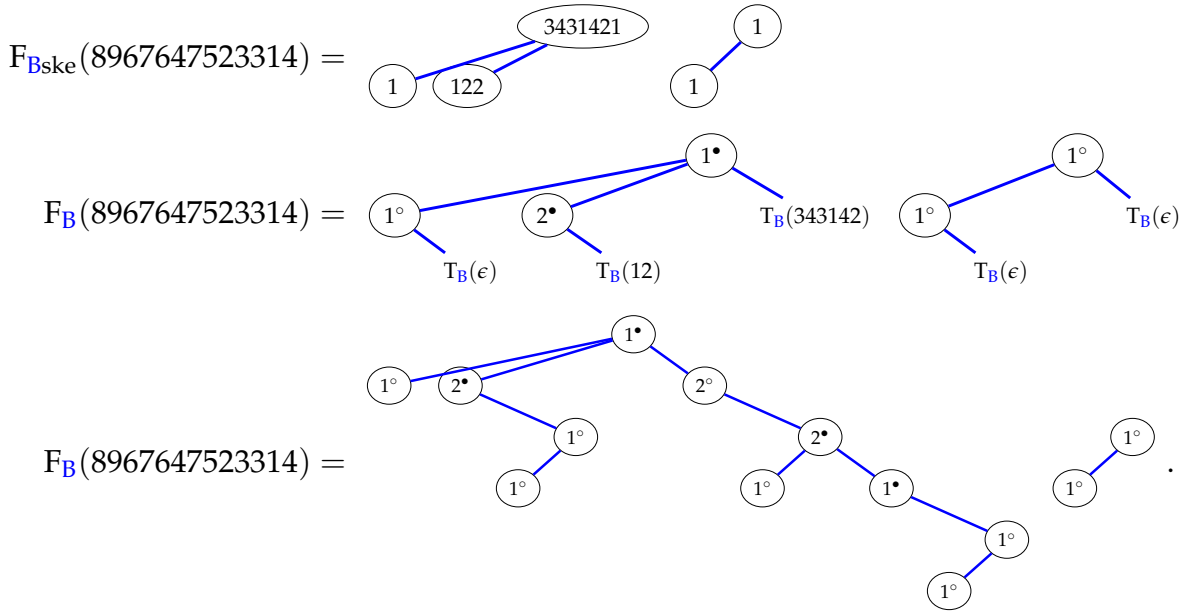
  \begin{align*}
    \Fbske(8967647523314) & = \scalebox{0.7}{\input{figures/arbres_pw/arbre8967647523314leftske}}\\
    \Fb(8967647523314) & = \scalebox{0.7}{\input{figures/arbres_pw/arbre8967647523314left_1}}\\
    \Fb(8967647523314) & = \scalebox{0.6}{\input{figures/arbres_pw/arbre8967647523314left}}.
  \end{align*}
\end{ex}

\begin{defi} \label{blue_weight}
  Let $t$ be a labeled biplane tree. We write
  $t = \Nodeblue(i^\alpha, f_\ell, f_r)$ where $i\in\NN_{>0}$,
  $\alpha \in \{\circ, \bullet\}$, $f_\ell = [\ell_1, \ldots, \ell_g]$, and
  $f_r = [r_1, \ldots, r_d]$, which is depicted as follows:
  \begin{equation*}
  	t = \scalebox{0.8}{{ \newcommand{\nodea}{\node[draw,ellipse] (a) {$i^\alpha$}
;}\newcommand{\nodeb}{\node (b) {$\ell_1$}
;}\newcommand{\nodebc}{\node (bc) {$\ldots$}
;}\newcommand{\nodec}{\node (c) {$\ell_g$}
;}\newcommand{\nodee}{\node (e) {$r_1$}
;}\newcommand{\nodeef}{\node (ef) {$\ldots$}
;}\newcommand{\nodef}{\node (f) {$r_d$}
;}
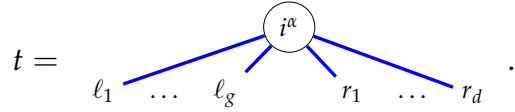
\begin{tikzpicture}[auto]
\matrix[column sep=.3cm, row sep=.3cm,ampersand replacement=\&]{
         \&         \&        \& \nodea  \&        \&         \&        \\ 
 \nodeb  \& \nodebc \& \nodec \&         \& \nodee \& \nodeef \& \nodef \\
};

\path[ultra thick, blue] (a) edge (b) edge (c) edge (e) edge (f);

\end{tikzpicture}}}.
  \end{equation*}
  
  The \textbf{weight} of $t$ ($\omega(t)$) is the number of nodes with $\circ$ in
  $t$. By extension, the \textbf{weight} of a forest is the sum of the weight of
  its trees.
\end{defi}

\begin{lem} \label{Fbweight}
  The weight of a forest (resp. a tree) obtain by the functions $\Fb$
  (resp. $\Tb$) is equal to the maximum value of the word. \textit{i.e.}
  $\forall w \in \PW, \omega(\Fb(w)) = \max(w)$, $\forall w \in \PW$ with $w$ irreducible,
  $\omega(\Tb(w)) = \max(w)$.
\end{lem}

\begin{proof}
  We prove by induction with the following hypothesis, for $n \in \NN$:
  \begin{equation}
    \begin{array}{rl}
      \forall w \in \PW_n, & \omega(\Fb(w)) = \max(w),\\
      \forall w \in \PW_n \text{ with $w$ irreducible}, & \omega(\Tb(w)) = \max(w).
    \end{array}
    \label{H_nsize_dual}
  \end{equation}
  The base case is given by the first item of
  \cref{def:construction_blue_packed} as $\Fb(\epsilon) = []$ and
  $\omega([]) = \max(\epsilon) = 0$ by convention.

  \noindent Let us fix $n \geq 1$ and suppose that the hypothesis (\ref{H_nsize_dual})
  holds. Let $w \in \PW_{n+1}$ and ${w = w_1\gcdot w_2\gcdot\dots\gcdot w_k}$ be
  the global descent decomposition of $w$.

  \noindent $\bullet$ If $k = 1$, we have $\Fb(w) = [\Tb(w)]$. Let
  $w = u \insl \psi_{i^\alpha}(v)$ with $i\in\NN_{>0}$ and
  $\alpha \in \{\circ, \bullet\}$ be the blue-factorization of $w$, then, depending of
  $\alpha$ the node is counted or not: \newline
  ${\omega(\Tb(w)) = \omega(\Nodeblue(i^\alpha, \Fb(u), \Fb(v))) = \overbrace{(1 +)}^{\text{if
      }\alpha=\circ} \omega(\Fb(u)) + \omega(\Fb(v))}$.  The sizes of $u$ and
  $v$ are at most $n$ so by induction
  ${\omega(\Tb(w)) = (1 +) \max(u) + \max(v) = \max(w)}$.

  \noindent $\bullet$ If $k \geq 2$, by induction on each factors, we have that
  \[\omega(\Fb(w)) = \omega(\Tb(w_1)) + \dots + \omega(\Tb(w_k)) = \max(w_1) + \dots +
    \max(w_k) = \max(w).\qedhere\]
\end{proof}



\begin{defi}\label{blue_packed_forest}
  Using the same notations as in previous \cref{blue_weight}, we say that $t$ is
  a \textbf{blue-packed tree} if it satisfies:
  \begin{equation}
    \left\{
        \begin{array}{l@{\,}l@{\,}l}
        d &=& 0, \\
        i^\alpha &=& 1^\circ,\\
        \multicolumn{3}{l}{\text{$\ell_1,\ldots, \ell_g$ are blue-packed trees.}}
      \end{array}
    \right.
    \text{or }
    \left\{
      \begin{array}{l@{\,}l@{\,}l}
        d &=& 1,\\
        i^\alpha &\neq& 1^\circ,\\
        \multicolumn{3}{l}{1 \leq i \leq \omega(r_1),}\\
        \multicolumn{3}{l}{\text{$\ell_1,\ldots, \ell_g$ and $r_1$ are blue-packed trees.}}
      \end{array}
    \right.
  \end{equation}
  An ordered list of blue-packed trees is a \textbf{blue-packed forest}.
\end{defi}

\begin{rem} \label{blue_ske_flat}
  The same remark as \cref{red_ske_flat} can be done with blue-skeleton trees
  that can be interpreted as flattened representations of blue-packed
  trees. Symmetrically, blue-packed trees can be interpreted as unfolded
  representations of blue-skeleton trees. We use the operation $\psi_{i^\alpha}$ to
  change between blue-packed and blue-skeleton trees.
\end{rem}

\begin{nota}\label{notation_blue}
  In the same way as~\cref{notation_red} we add the following notations:
  \begin{itemize}
  \item $\PForestbluen$ the set of blue-packed forests of size $n$,
    ($\PForestbluen \eqdef \{\Fb(w)\}_{w\in\PW_n}$),
  \item $\PTreebluen$ the set of blue-packed trees of size $n$,
    ($\PTreebluen \eqdef \{\Tb(w)\}_{w\in\PW_n}$ with $w$ irreducible),
  \item $\PParticularbluen$ the set of blue-packed trees of size $n$ such that
    the left forest of the root is empty,
    ($\PParticularbluen \eqdef \{\Tb(w)\}_{w\in\PW_n}$ with $w$
    blue-irreducible). In particular, the blue-skeleton of a tree of
    $\PParticularbluen$ consist of a single node labeled by a blue-irreducible
    word.
  \end{itemize}
\end{nota}

\begin{rem} \label{rem:PParticularbluen}
  The set $\PParticularbluen$ can be described as a disjointed union of sets
  depending on $i$ and $\alpha$. Let $\PForestbluen^{i^\alpha}$ denote the set of
  blue-packed forests of weight $n$ that can be right children of a node labeled
  by $i^\alpha$ (see~\cref{blue_packed_forest} for conditions), we have the following
  description:
  \begin{align} \label{eq:PParticularbluen}
    \PParticularbluen = \bigsqcup_{i,\alpha}~\{\Nodeblue(i^\alpha,[],f_r)~|~f_r \in \PForestbluenmp^{i^\alpha}\}
  \end{align}
\end{rem}

Analogously, we use $\PForestblueske$, $\PTreeblueske$ and $\PParticularblueske$
for blue-skeleton forests, trees and trees with only one node.

\noindent We can remark that for $n = 1$ we have
$\PParticularblueone = \PTreeblueone = \PForestblueone$ and
$\forall n > 1, \PParticularbluen \subsetneq \PTreebluen \subsetneq \PForestbluen$.

Once again we define two functions in order to prove that $\Fb$ and $\Tb$ are
bijections.

\begin{defi}\label{def:Fb-1}
  We define here the functions $\Fbstar$ (resp. $\Tbstar$) that transform blue-packed
  forest $f$ (resp. tree $t$) into a packed word. We reverse all instructions of
  \cref{def:construction_blue_packed} as follows:
  \begin{itemize}
  \item $\Fbstar([]) = \epsilon$,
  \item for any non empty blue-packed forest $f = [t_1, t_2\dots, t_k]$, then \newline
    $\Fbstar(f) = \Tbstar(t_k)\gcdot \Tbstar(t_{k-1})\gcdot\dots\gcdot
    \Tbstar(t_1)$ \newline (notice the inversion compared to~\cref{def:Fr-1}).
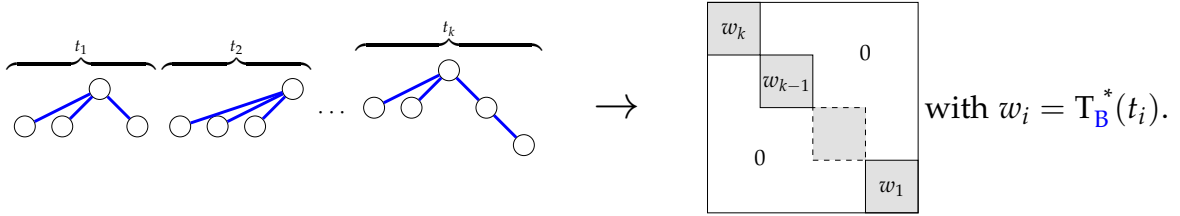
     $$\scalebox{0.7}{$\overbrace{{ \newcommand{\nodea}{\node[draw,ellipse] (a) {}
;}\newcommand{\nodeb}{\node[draw,ellipse] (b) {}
;}\newcommand{\nodec}{\node[draw,ellipse] (c) {}
;}\newcommand{\noded}{\node[draw,ellipse] (d) {}
;}\begin{tikzpicture}[baseline={([yshift=-1ex]current bounding box.center)}]
\matrix[column sep=.3cm,row sep=.3cm,ampersand replacement=\&]{
         \&         \& \nodea \& \\ 
 \nodeb  \& \nodec  \&        \& \noded  \\
};
\path[ultra thick,blue] (a) edge (b) edge (c) edge (d);
\end{tikzpicture}}}^{t_1}$
       $\overbrace{\input{figures/arbres_pw/arbrew2_left}}^{t_2}$
       $\cdots$ $\overbrace{{ \newcommand{\nodea}{\node[draw,ellipse] (a) {}
;}\newcommand{\nodeb}{\node[draw,ellipse] (b) {}
;}\newcommand{\nodec}{\node[draw,ellipse] (c) {}
;}\newcommand{\noded}{\node[draw,ellipse] (d) {}
;}\newcommand{\nodee}{\node[draw,ellipse] (e) {}
;}\begin{tikzpicture}[baseline={([yshift=-1ex]current bounding box.center)}]
\matrix[column sep=.3cm,row sep=.3cm,ampersand replacement=\&]{
         \&         \& \nodea \& \& \\ 
 \nodeb  \& \nodec  \&        \& \noded \& \\
 \&  \&        \& \& \nodee  \\
};
\path[ultra thick,blue] (a) edge (b) edge (c) edge (d)
(d) edge (e);
\end{tikzpicture}}}^{t_k}$
       \qquad\scalebox{2}{$\to$}\qquad\quad \input{figures/box_d/gd_fact_left}} \text{with
       $w_i = \Tbstar(t_i)$}.$$
  \item for any non empty blue-packed tree $t=\Nodeblue(i^\alpha,f_\ell, f_r)$, then \newline
    $\Tbstar(t) = \Fbstar(f_\ell)\insl\psi_{i^\alpha}(\Fbstar(f_r))$.
    $$\scalebox{0.7}{{ \newcommand{\nodea}{\node[draw,ellipse] (a) {$i^\alpha$}
;}\newcommand{\nodeb}{\node (b) {$\ell_1$}
;}\newcommand{\nodebc}{\node (bc) {$\ldots$}
;}\newcommand{\nodec}{\node (c) {$\ell_g$}
;}\newcommand{\noded}{\node (d) {$\bullet$}
;}\newcommand{\nodee}{\node (e) {$r_1$}
;}\newcommand{\nodeef}{\node (ef) {$\ldots$}
;}\newcommand{\nodef}{\node (f) {$r_d$}
;}\begin{tikzpicture}[auto]
\matrix[column sep=.3cm, row sep=.3cm,ampersand replacement=\&]{
      \&  \&  \& \nodea  \& \&    \&    \\ 
 \nodeb  \& \nodebc \& \nodec \&  \& \nodee \& \nodeef \& \nodef  \\
};

\path[ultra thick, blue] (a) edge (b) edge (c) edge (e) edge (f);
\draw[snake=brace,raise snake = 4mm,mirror snake] (b.west) -- (c.east);
\node[below of = bc] (F) {$f_\ell$};
\draw[snake=brace,raise snake = 4mm,mirror snake] (e.west) -- (f.east);
\node[below of = ef] (F) {$f_r$};

\end{tikzpicture}} \qquad\scalebox{2}{$\to$}\qquad\quad \input{figures/box_d/insl}} \text{with $v' = \Fbstar(f_r)$ and $u = \Fbstar(f_\ell)$}.$$
  \end{itemize}
  \vspace{-0.7cm}
  As $\psi_{i^\alpha}(\Fbstar(f_r))$ is only defined if
  $i \leq \max(\Fbstar(f_r)) \overbrace{(+ 1)}^{if \alpha = \circ}$ (see \cref{psi}), there
  might be a problem with this definition. We prove in the following
  \cref{Fb-1welldef} that this is the case if $t \in \PTreeblue$.
\end{defi}

\begin{lem}\label{Fb-1welldef}
  For any blue-packed forest $f$, $\Fbstar(f)$ is a well defined word and its
  maximum value is $\omega(f)$. For any blue-packed tree $t$, $\Tbstar(t)$ is a well
  defined word and its maximum is $\omega(t)$.
\end{lem}

\begin{proof}
  We prove by induction with the following hypothesis, for $n \in \NN$:
  \begin{equation}
    \begin{array}{rl}
      \forall f \in \PForestblueinfn,
      & \Fbstar(f) \text{ is well defined and } \max(\Fbstar(f)) = \omega(f),\\
      \forall t \in \PTreeblueinfn,
      & \Tbstar(t) \text{ is well defined and } \max(\Tbstar(t)) = \omega(t).
    \end{array}
    \label{H_nwelldefblue}
  \end{equation}
  The base case is given by the first item of
  \cref{def:Fb-1} as $\Fbstar([]) = \epsilon$ and $\omega([]) = \max(\epsilon) = 0$.

  \noindent Let us fix $n \geq 1$ and suppose that the hypothesis
  (\ref{H_nwelldefblue}) holds. Let
  $f = [t_1, \dots, t_k] \in \PForestbluenplusone$.

  \noindent $\bullet$ If $k = 1$, it is sufficient to prove the second item of
  (\ref{H_nwelldefblue}). Let $t = \Nodeblue(i^\alpha, f_\ell, f_r) \in
  \PTreebluenplusone$. According to \cref{blue_packed_forest} with notations of
  \cref{blue_weight} there are two cases:

  $\centerdot$ $d = 0$ and $i^\alpha = 1^\circ$. We have that
  $\max(\Fbstar(f_r)) = \max(\epsilon) = 0$ so $i \leq 0 + 1$ and
  $\psi_{i^\alpha}(\epsilon) = 1$. Now by induction on $f_\ell$, we have that
  $\Fbstar(f_\ell)$ is a well defined word and its maximum value is
  $\omega(f_\ell)$. Finally
  $\Tbstar(t) = \Fbstar(f_\ell) \insl \psi_{i^\alpha}(\epsilon) = \Fbstar(f_\ell) \dcdot \psi_{i^\alpha}(\epsilon)$
  is a well defined word and its maximum is the last value:
  $\max(\Fbstar(f_\ell)) + 1 = \omega(t)$.

  $\centerdot$ $d = 1$. In this case, we can directly apply the hypothesis
  (\ref{H_nwelldefblue}) on $f_r$ and $f_\ell$. According to \cref{blue_packed_forest}
  we have that
  \begin{eqnarray*}
    i &\leq& \omega(r_1),\\
    i &\leq& \max(\Fbstar(f_r)) .
  \end{eqnarray*}
  So $\Tbstar(t)$ is well defined. Moreover
  \begin{eqnarray*}
    \max(\Tbstar(t)) &=& \max(\Fbstar(f_\ell)\insl\psi_{i^\alpha}\Fbstar(f_r))\\
                     &=& \max(\Fbstar(f_\ell)) + \max(\Fbstar(f_r)) (+1)\\
                     &=& \omega(f_\ell) + \omega(f_r) (+1) = \omega(t).
  \end{eqnarray*}
  \noindent $\bullet$ If $k \geq 2$, the weight of trees are at least $1$ so we can apply
  (\ref{H_nwelldefblue}) on trees of $f$.
\end{proof}

\begin{theorem} \label{bij_blue}
  The functions $\Fb$ and $\Fbstar$ (resp. $\Tb$ and $\Tbstar$) are two converse
  bijections between packed words of size $n$ and blue-packed forests (resp.
  irreducible packed words and blue-packed trees) of size $n$. That is to
  say $\Fbinv = \Fbstar$ and $\Tbinv = \Tbstar$.
\end{theorem}

\begin{proof}
  The proof structure is the same as the one of \cref{bij_red} which is similar
  to the one of \cref{lem:red_skeleton_bij,lem:blue_skeleton_bij}. But with use
  of statements comming from this subsection. We can see in this table some of
  the main statements that are exchanged with their counterpart:
  \begin{center}
    \begin{tabular}{|l|l|l|}
      \hline
      \cref{red_fact} & \cref{blue_fact} & red-factorization and blue-factorization.\\
      \hline
      \cref{red_irreducible} & \cref{blue_irreducible} & red-irreducible words and blue-irreducible words.\\
      \hline
      \cref{ins_irr} &                 & $u \ins v$ irreducible $\iff$ $v$ irreducible \\
                     & \cref{insl_irr} & $u \insl v$ irreducible $\iff$ $v$ irreducible.\\
      \hline
      \cref{red_fact_irr} &  & $u \ins v$ red-factorization $\iff$ $v$ red-irreducible \\
      & \cref{blue_fact_irr} & $u \insl v$ blue-factorization $\iff$ $v$ blue-irreducible.\\
      \hline
      \cref{def:construction_red_packed} & \cref{def:construction_blue_packed} & $\Fr$, $\Tr$ and $\Fb$, $\Tb$.\\
      \hline
      \cref{red_weight} & & weight of red-forests ($\sum$ size of labels)\\
      & \cref{blue_weight} & weight of blue-forests ($\sum$ nodes with $\circ$).\\
      \hline
      \cref{Frweight} & \cref{Fbweight} & $\omega(\Fr(w)) = |w|$ and $\omega(\Fb(w)) = \max(w)$.\\
      \hline
      \cref{red_packed_forest} & \cref{blue_packed_forest} & red-packed forest and blue-packed forest.\\
      \hline
      \cref{def:Fr-1} & \cref{def:Fb-1} & $\Frstar$, $\Trstar$ and $\Fbstar$, $\Tbstar$.\\
      \hline
    \end{tabular}
  \end{center}
\end{proof}

\begin{rem}
  As we can see in~\cref{Primal_Red,Dual_Blue}, the role of size and weight are
  exchanged for red and blue-forests. For red forests, the size (number of
  nodes) is equal to the maximum letter of the word associated while the weight
  (\cref{red_weight}) is the number of letter of the associated word. For
  blue-forests, it is the opposite, the number of letters of the associated word
  is equal to the size of the forest while le maximum letter is equal to the
  weight (\cref{blue_weight}) of the forest. That is why we denote the set
  of red packed forests of \textbf{weight} $n$ by $\PForestredn$ and the set of
  blue-packed forests of \textbf{size} $n$ by $\PForestbluen$.
\end{rem}

\begin{ex}
  There is a unique forest in $\PForestblueone$, namely
  \scalebox{0.5}{{ \newcommand{\nodea}{\node[draw,ellipse] (a) {$1^\circ$}
;}\begin{tikzpicture}[baseline={([yshift=-1ex]current bounding box.center)}]
\matrix[column sep=.3cm,row sep=.3cm,ampersand replacement=\&]{
 \nodea  \\
};
\path[ultra thick,blue] ;
\end{tikzpicture}}}, here are the 3 forests
  of $\PForestbluetwo$with the associated packed
  word:$\Fb(12)=$\scalebox{0.5}{{ \newcommand{\nodea}{\node[draw,ellipse] (a) {$1^\circ$}
;}\newcommand{\nodeb}{\node[draw,ellipse] (b) {$1^\circ$}
;}\begin{tikzpicture}[baseline={([yshift=-1ex]current bounding box.center)}]
\matrix[column sep=.3cm,row sep=.3cm,ampersand replacement=\&]{
         \& \nodea  \\ 
 \nodeb  \&         \\
};
\path[ultra thick,blue] (a) edge (b);
\end{tikzpicture}}},
  $\Fb(21)=$\scalebox{0.5}{\input{figures/arbres_pw/arbre21left}},
  $\Fb(11)=$\scalebox{0.5}{{ \newcommand{\nodea}{\node[draw,ellipse] (a) {$1^\bullet$}
;}\newcommand{\nodeb}{\node[draw,ellipse] (b) {$1^\circ$}
;}\begin{tikzpicture}[baseline={([yshift=-1ex]current bounding box.center)}]
\matrix[column sep=.3cm,row sep=.3cm,ampersand replacement=\&]{
         \& \nodea  \&         \\ 
         \&         \& \nodeb  \\
};
\path[ultra thick,blue] (a) edge (b);
\end{tikzpicture}}}. We show below
  the forests of $\PForestbluethree$:

  \noindent
  $\Fb(123)=$ \scalebox{0.4}{\input{figures/arbres_pw/arbre123left}},
  $\Fb(132)=$ \scalebox{0.4}{\input{figures/arbres_pw/arbre132left}},
  $\Fb(213)=$ \scalebox{0.4}{\input{figures/arbres_pw/arbre213left}},
  $\Fb(231)=$ \scalebox{0.4}{\input{figures/arbres_pw/arbre231left}},
  $\Fb(312)=$ \scalebox{0.4}{\input{figures/arbres_pw/arbre312left}},
  $\Fb(321)=$ \scalebox{0.4}{\input{figures/arbres_pw/arbre321left}},
  $\Fb(122)=$ \scalebox{0.4}{\input{figures/arbres_pw/arbre122left}},
  $\Fb(212)=$ \scalebox{0.4}{\input{figures/arbres_pw/arbre212left}},
  $\Fb(221)=$ \scalebox{0.4}{\input{figures/arbres_pw/arbre221left}},
  $\Fb(112)=$ \scalebox{0.4}{\input{figures/arbres_pw/arbre112left}},
  $\Fb(121)=$ \scalebox{0.4}{\input{figures/arbres_pw/arbre121left}},
  $\Fb(211)=$ \scalebox{0.4}{\input{figures/arbres_pw/arbre211left}},
  $\Fb(111)=$ \scalebox{0.4}{\input{figures/arbres_pw/arbre111left}}.
  
  More examples can be found in the annexes section with
  \cref{F(w)123,F(1234),F(1233),F(1223),F(1123),F(1111)}.
\end{ex}

We conclude by the main theorem of this subsection. It is the dual of \cref{thm:bij_primal}.
\begin{theorem} \label{thm:bij_dual} For all $n\in\NN$ we have the three following
  equalities :
  \begin{equation*}
    \dim(\WQSym_n) = \#\PForestbluen
    \quad\text{and}\quad
    \dim(\Prim_n) = \#\PTreebluen
    \quad\text{and}\quad
    \dim(\TPrim_n) = \#\PParticularbluen
  \end{equation*}
\end{theorem}

\begin{proof}
  The proof is similar to the one of~\cref{thm:bij_primal} thanks to \cref{bij_blue}
  instead of \cref{bij_red}.
\end{proof}


\section{Bases for totally primitive elements}\label{sect3_bases}
In this section we construct two bases of primitive and totally primitive
elements of $\WQSym$ and $\WQSym^*$. Thanks
to~\cref{thm:bij_primal,thm:bij_dual} we now have the combinatorial objects to
index those bases and we know that their numbers agree with the dimensions. We
therefore only need to show that they are linearly independent. We will proceed
by showing that the decompositions through maximum and through last letter
preserve the total primitivity.

As in \cref{sect2_pw_trees}, we start by working on $\WQSym^*$ associated to the
color \red{red} and do the same work on the primal $\WQSym$ associated to the
color \blue{blue}.

\subsection{Dual (\red{Red})}

\subsubsection{Decomposition through maximums and totally primitive elements}

\begin{defi}\label{Phi}
  Let $I = [i_1, \ldots, i_p]$ with $0 <i_1< \ldots< i_p$. We define a linear map
  $\Phi_I:\WQSym^* \to \WQSym^*$ as follows: for all $n\in\NN$ and
  $w = w_1\cdot w_2 \cdots w_n \in\PW_n$,
  \begin{equation}
    \Phi_I(\RR_w) \eqdef \left\{
      \begin{array}{rl}
        \RR_{\phi_I(w)} & \text{if } i_p \leq n + p,\\
        0 & \text{if } i_p > n + p\,.
      \end{array}\right.
  \end{equation}
\end{defi}

\begin{defi}\label{tau_I}
  Let $I = (i_1, \ldots, i_p)$ with $0 <i_1< \ldots< i_p$. We define a projector
  $\tau_I: \WQSym^* \to \WQSym^*$ as follows: for all $n\in\NN$ and
  $w = w_1\cdot w_2 \cdots w_n \in\PW_n$,
  \begin{equation}
    \tau_I(\RR_w) \eqdef \left\{
      \begin{array}{rl}
        \RR_w & \text{if } w_i = \max(w) \text{ if and only if } i \in I,\\
        0 & \text{else.}
      \end{array}\right.
  \end{equation}
  These are orthogonal projectors in the sense that $\tau_I^2 = \tau _I$ and $\tau_I \circ \tau_J = 0$ ($I \neq J$).
\end{defi}

\begin{lem}\label{im_eq_phi}
  For any $I$, we have $\Imm(\Phi_I) = \Imm(\tau_I)$ where
  $\Imm(f)$ denotes the image of $f$.
\end{lem}
\begin{proof}
  For any $I$, the inclusion $\Imm(\Phi_I) \subset \Imm(\tau_I)$ is automatic by definition
  of $\Phi_I$ and $\tau_I$. Indeed, for any $w\in\PW_n$ if $i_p \leq n + p$ then
  $\Phi_I(\RR_w) = \RR_{\phi_I(w)}$ and
  $\tau_I(\RR_{\phi_I(w)}) = \RR_{\phi_I(w)}$ and $\Phi_I(\RR_w) = 0$ otherwise. By
  linearity $\Imm(\Phi_I) \subset \Imm(\tau_I)$.

  For any $I$, the inclusion $\Imm(\Phi_I) \supset \Imm(\tau_I)$ is a consequence
  of~\cref{phi_bij} and linearity. Indeed, for any $w \in \PW$,
  $\tau_I(\RR_w) = \RR_w \Leftrightarrow (w_{i} = \max(w) \Leftrightarrow i \in I)$. If
  $\tau_I(\RR_w) = \RR_w$ let $w'$ be such that $\phi_I(w') = w$ using \cref{phi_bij},
  then $\Phi_I(\RR_{w'}) = \RR_w = \tau_I(\RR_w)$. By linearity $\Imm(\Phi_I) \supset \Imm(\tau_I)$.
\end{proof}

\begin{lem}\label{tot_stable_phi}
  For any~$I$, the projection by $\tau_I$ of a totally primitive element is still a
  totally primitive element, so that
  $\tau_I(\TPrim^*) = \Imm(\tau_I) \bigcap \TPrim^*$. Moreover,
  \begin{equation}\label{eq:direct_sum}
    \TPrim^* = \bigoplus_I \Imm(\tau_I) \cap \TPrim^*\,.
  \end{equation}
\end{lem}

\begin{proof}
  Let $w$ a packed word. We have
  $\Delta_\prec(\tau_I(\RR_w)) = (\tau_I \otimes \Id) \circ \Delta_\prec(\RR_w)$ by definition of
  $\tau_I$ and $\Delta_\prec$. Indeed, in $ \Delta_\prec(\RR_w)$, the deconcatenations cannot be done
  before the last maximum letter of $w$. By linearity, for all $p \in \TPrim^*$, we have
  $\Delta_\prec(\tau_I(p)) = (\tau_I \otimes \Id) \circ \Delta_\prec(p) = 0$. The same argument works on the right
  so that~$\tau_I(p)\in\TPrim^*$. Morevover $\tau_I$ are orthogonal projectors so
  $\TPrim^* = \bigoplus_I \tau_I(\TPrim^*) = \bigoplus_I \Imm(\tau_I)\cap\TPrim^*$.
\end{proof}

  





\subsubsection{The new basis ~\texorpdfstring{$\PP$}{PP}}
\begin{defi}\label{P}
  Let $t_1, \ldots, t_k$ be $k$ red-packed trees, $I= [i_1, \dots, i_p]$,
  $f_l = [\ell_1, \ldots, \ell_g]$ be a red-packed forest and
  $f_r \in \PForestred^I$ be a red-packed forest that can be right children of a
  node labeled by $I$,
  \begin{align}
    \label{PE}  \PP_{[]}  & \eqdef \RR_\epsilon,\\
    \label{PF} \PP_{t_1, \ldots, t_k} & \eqdef \PP_{t_k} \prec (\PP_{t_{k-1}} \prec (\ldots  \prec \PP_{t_1})\ldots),\\
    \label{PT} \PP_{\Nodered(I, f_l = [\ell_1, \ldots, \ell_g], f_r)} & \eqdef \langle\PP_{\ell_1}, \PP_{\ell_2}, \ldots, \PP_{\ell_g}; \PP_{\Nodered(I, [], f_r)}\rangle,\\
    \label{PP} \PP_{\Nodered(I, [], f_r)} & \eqdef \Phi_{I}(\PP_{f_r}).
  \end{align}
\end{defi}

\begin{ex}\label{ex:P}
  \begin{align*}
    \PP_{\scalebox{0.5}{\input{figures/arbres_pw/arbre43412}}} &= \PP_{\scalebox{0.5}{{ \newcommand{\nodea}{\node[draw,ellipse] (a) {$\,1\,$}
;}\newcommand{\nodeb}{\node[draw,ellipse] (b) {$\,1\,$}
;}\begin{tikzpicture}[baseline={([yshift=-1ex]current bounding box.center)}]
\matrix[column sep=.3cm,row sep=.3cm,ampersand replacement=\&]{
         \& \nodea  \\ 
 \nodeb  \&         \\
};
\path[ultra thick,red] (a) edge (b);
\end{tikzpicture}}}} \prec \PP_{\scalebox{0.5}{{ \newcommand{\nodea}{\node[draw,ellipse] (a) {$1, 3$}
;}\newcommand{\nodeb}{\node[draw,ellipse] (b) {$\,1\,$}
;}\begin{tikzpicture}[baseline={([yshift=-1ex]current bounding box.center)}]
\matrix[column sep=.3cm,row sep=.3cm,ampersand replacement=\&]{
         \& \nodea  \&         \\ 
         \&         \& \nodeb  \\
};
\path[ultra thick,red] (a) edge (b);
\end{tikzpicture}}}} = (\PP_{\scalebox{0.5}{{ \newcommand{\nodea}{\node[draw,ellipse] (a) {$\,1\,$}
;}\begin{tikzpicture}[baseline={([yshift=-1ex]current bounding box.center)}]
\matrix[column sep=.3cm,row sep=.3cm,ampersand replacement=\&]{
 \nodea  \\
};
\path[ultra thick,red] ;
\end{tikzpicture}}}} \succ \PP_{\scalebox{0.5}{{ \newcommand{\nodea}{\node[draw,ellipse] (a) {$\,1\,$}
;}\begin{tikzpicture}[baseline={([yshift=-1ex]current bounding box.center)}]
\matrix[column sep=.3cm,row sep=.3cm,ampersand replacement=\&]{
 \nodea  \\
};
\path[ultra thick,red] ;
\end{tikzpicture}}}} - \PP_{\scalebox{0.5}{{ \newcommand{\nodea}{\node[draw,ellipse] (a) {$\,1\,$}
;}\begin{tikzpicture}[baseline={([yshift=-1ex]current bounding box.center)}]
\matrix[column sep=.3cm,row sep=.3cm,ampersand replacement=\&]{
 \nodea  \\
};
\path[ultra thick,red] ;
\end{tikzpicture}}}} \prec \PP_{\scalebox{0.5}{{ \newcommand{\nodea}{\node[draw,ellipse] (a) {$\,1\,$}
;}\begin{tikzpicture}[baseline={([yshift=-1ex]current bounding box.center)}]
\matrix[column sep=.3cm,row sep=.3cm,ampersand replacement=\&]{
 \nodea  \\
};
\path[ultra thick,red] ;
\end{tikzpicture}}}}) \prec \Phi_{1,3}(\PP_{\scalebox{0.5}{{ \newcommand{\nodea}{\node[draw,ellipse] (a) {$\,1\,$}
;}\begin{tikzpicture}[baseline={([yshift=-1ex]current bounding box.center)}]
\matrix[column sep=.3cm,row sep=.3cm,ampersand replacement=\&]{
 \nodea  \\
};
\path[ultra thick,red] ;
\end{tikzpicture}}}})\\
                                                               &= \RR_{14342} + \RR_{41342} + \RR_{43142} + \RR_{43412} - \RR_{24341} - \RR_{42341} - \RR_{43241} - \RR_{43421}
  \end{align*}
  More examples can be found in the annexes section with \cref{RPOQ123} and \cref{matrPR3,matrPR4}.
\end{ex}

\begin{theorem}\label{thm:P}
  For all $n \in \NN_{>0}$
  \begin{multicols}{2}
  \begin{enumerate}
  \item \label{PForest} $(\PP_f)_{f\in\PForestredn}$ is a basis of $\WQSym^*_n$,
  \item \label{PTree} $(\PP_t)_{t\in\PTreeredn}$ is a basis of $\Prim^*_n$,
  \item \label{PParticular} $(\PP_t)_{t\in\PParticularredn}$ is a basis of $\TPrim^*_n$.
  \end{enumerate}
  \end{multicols}
\end{theorem}

\begin{proof}
  We do a mutually recursive induction on $n$ to prove these three
  items. \linebreak As~${\dim(\WQSym^*_1) = \dim(\Prim^*_1) = \dim(\TPrim^*_1) = 1}$ the base
  case is trivial. By \cref{left_prod}, \cref{PTree} up to degree~$n$ implies
  \cref{PForest} up to degree~$n$. Similarly, \cref{brace} shows that
  \cref{PParticular} up to degree~$n$ implies \cref{PTree} up to degree~$n$. By
  induction it is sufficient to show that \cref{PForest,PTree} up to degree
  $n-1$ implies \cref{PParticular}
  for~$n$. 
  
  For all $k \in \NN$, let $\pi_k$ be the canonical projector on the homogeneous
  component of degree $k$ of $\WQSym^*$. We define
  $\pi_{<k} \eqdef \sum_{i=0}^{k-1}\pi_i$. Fix $I = [i_1, \ldots, i_p]$ with
  $p \leq n$ and $u$ a packed word of size~$n-p$. Notice that if $p = n$ we
  immediately have $u=\epsilon$ and $\Delta_\prec(\Phi_I(\RR_\epsilon)) = 0$. We suppose now that
  $p < n$. By \cref{RR_left_coprod}, in the half coproduct
  $\Delta_\prec(\Phi_I(\RR_u))$ all the maximums must be in the left tensor factor, which
  therefore must be at least of degree $i_p$. By linearity, for all
  $x \in \WQSym^*_{n-p}$,
  \begin{align} \label{left_coprod_phi}
    \Delta_\prec(\Phi_I(x)) = \left( \sum_{j = i_p}^{n-1} \Phi_I \circ \pi_j \otimes \pi_{n - 1 - j} \right) \circ \tilde{\Delta}(x).
  \end{align}
  Thanks to \cref{coro:prim_tot} for $f =[r_1, \ldots, r_d]\in \PForestrednmp$, the
  coproduct $\tilde{\Delta}(\PP_{f})$ is computed by deconcatenation of forests. So
  if $f \in \PForestrednmp^I$, in particular $r_d$ is of weight at least
  $n - i_p + 1$ then for $j \geq i_p$ we have $n-1-j < \omega(r_d)$ so that all the
  terms in the previous sum vanishe. A similar reasoning applies to
  $\Delta_\succ(\Phi_I(x))$ using the fact that $f \in \PForestrednmp^I$ imply
  $1 \leq i_1 \leq \omega(r_1)$.

  So for all $t = \Nodered(I, [], f_r) \in \PParticularredn$, we have that
  $\Delta_\prec(\PP_t) = \Delta_\prec(\Phi_I(\PP_{f_r})) = 0$ and
  $\Delta_\succ(\PP_t) = \Delta_\succ(\Phi_I(\PP_{f_r})) = 0$.

  Moreover, by induction we have that $\{\PP_f\mid f\in\PForestrednmp^I\}$ are
  linearity independent as $\{\PP_f\mid f\in\PForestrednmp\}$ is a basis of
  $\WQSym^*_{n-p}$. Recall the description of $\PParticularredn$ in
  \cref{rem:PParticularredn} as a disjointed union of sets depending on $I$:
  \begin{align}
    \PParticularredn = \bigsqcup_I~\{\Nodered(I,[],f_r) ~\mid~ f_r \in \PForestrednmp^I\}
  \end{align}
  Since $\Phi_I$ is injective on $\WQSym^*_{n-p}$ then
  $\{\Phi_I(\PP_f)\mid f\in\PForestrednmp^I\}$ are linearly independent. According to
  \cref{im_eq_phi}, for all $f\in\PForestrednmp^I$ we have
  $\Phi_I(\PP_f) \in \Imm(\tau_I) \cap \TPrim^*_n$. Moreover, thanks to the direct sum of
  \cref{eq:direct_sum}:
  $$\TPrim^*_n = \bigoplus_I \Imm(\tau_I) \cap \TPrim^*_n\,$$ and by
  definition of $\PP$, in particular
  $\PP_{\Nodered(I, [], f_r)} \eqdef \Phi_{I}(\PP_{f_r}),$ the family
  $\{\PP_t\mid t\in\PParticularredn\}$ are linearly independent. Finally, by
  cardinalities of \cref{thm:bij_primal} it is a basis of $\TPrim^*_n$.
\end{proof}
  


\begin{rem}\label{rem:Pske}
  The basis $\PP$ is indexed by red-packed forests ($\PForestred$). We will also
  use red-skeletons ($\PForestredske$) or packed words ($\PW$) as index thanks
  to the bijections of \cref{red_ske_flat} and $\Fr$ of
  \cref{def:construction_red_packed}.
\end{rem}


\subsection{Primal (\blue{Blue})}\label{Obasis}

\subsubsection{Decomposition through last letter and totally primitive elements}

\begin{defi}\label{Psi}
  Let $i \in \NN_{>0}$ and $\alpha \in \{\circ, \bullet\}$. We define a linear map
  $\Psi_{i^\alpha}:\WQSym \to \WQSym$ as follows: for all $n\in\NN$ and $w \in\PW_n$,
  \begin{equation}
    \Psi_{i^\alpha}(\QQ_w) \eqdef \left\{
      \begin{array}{rl}
        \QQ_{\psi_{i^\circ}(w)} & \text{if } \alpha = \circ \text{ and } 1 \leq i \leq \max(w) + 1,\\
        \QQ_{\psi_{i^\bullet}(w)} & \text{if } \alpha = \bullet \text{ and } 1 \leq i \leq \max(w),\\
        0 & \text{else.}
      \end{array}\right.
  \end{equation}
\end{defi}


\begin{defi}\label{tau_ialpha}
  Let $i \in \NN_{>0}$ and $\alpha \in \{\circ, \bullet\}$. We define a projector
  $\tau_{i^\alpha}: \WQSym \to \WQSym$ as follows: for all $n\in\NN$ and
  $w = w_1\cdot w_2 \cdots w_n \in\PW_n$,
  \begin{equation}
    \tau_{i^\alpha}(\QQ_w) \eqdef \left\{
      \begin{array}{rl}
        \QQ_w & \text{if } w_n = i \text{ and } \alpha = \bullet \text{ and } i \in [w_1, \ldots, w_{n-1}],\\
        \QQ_w & \text{if } w_n = i \text{ and } \alpha = \circ \text{ and } i \notin [w_1, \ldots, w_{n-1}],\\
        0 & \text{else.}
      \end{array}\right.
  \end{equation}
  These are orthogonal projectors in the sense that $\tau_{i^\alpha}^2 = \tau_{i^\alpha}$ and $\tau_{i^\alpha} \circ \tau_{j^\beta} = 0$ ($i \neq j$ or $\alpha \neq \beta$).
\end{defi}

\begin{lem}\label{im_eq_psi} For any $i$ and $\alpha$, we have $Im(\Psi_{i^\alpha}) = Im(\tau_{i^\alpha})$ where
  $Im(f)$ denotes the image of $f$.
\end{lem}
\begin{proof}
  For any $i$ and $\alpha$, the inclusion
  $Im(\Psi_{i^\alpha}) \subset Im(\tau_{i^\alpha})$ is automatic by definition of
  $\Psi_{i^\alpha}$ and $\tau_{i^\alpha}$ and linearity. Indeed, for any $w\in\PW_n$
  $\tau_{i^\alpha}(\Psi_{i^\alpha}(\QQ_w)) = \Psi_{i^\alpha}(\QQ_w)$.

  For any $i$ and any $\alpha$, the inclusion
  $Im(\Psi_{i^\alpha}) \supset Im(\tau_{i^\alpha})$ is a consequence of~\cref{psi_bij} and
  linearity. Indeed, for any $w \in \PW$, if $\tau_{i^\alpha}(\QQ_w) = \QQ_w$ then
  $w_n = i$. With $w' = \pack(w_1\dots w_{n-1})$ we have
  $\Psi_{i^\alpha}(\QQ_{w'}) = \QQ_w = \tau_{i^\alpha}(\QQ_w)$.
\end{proof}

\begin{lem} \label{tot_stable_psi} For any~$i$ and $\alpha$, the projection by
  $\tau_{i^\alpha}$ of a totally primitive element is still a totally primitive element, so
  that $\tau_{i^\alpha}(\TPrim) = Im(\tau_{i^\alpha}) \bigcap \TPrim$. Moreover,
  \begin{equation}\label{eq:direct_sum_dual}
    \TPrim = \bigoplus_{\alpha, i} Im(\tau_{i^\alpha}) \cap \TPrim\,.
  \end{equation}
\end{lem}
\begin{proof}
  Let $w$ a packed word. We have
  $\Delta_{\valprec}(\tau_{i^\alpha}(\QQ_w)) = (\tau_{i^\alpha} \otimes Id) \circ \Delta_{\valprec}(\QQ_w)$ by definition of
  $\tau_{i^\alpha}$ and $\Delta_{\valprec}$. Indeed, in $ \Delta_{\valprec}(\QQ_w)$, the decomposition can't be done
  under the last letter of $w$. By linearity, for all $p \in \TPrim$, we have
  $\Delta_{\valprec}(\tau_{i^\alpha}(p)) = (\tau_{i^\alpha} \otimes Id) \circ \Delta_{\valprec}(p) = 0$. The same argument works on the right so
  that~$\tau_{i^\alpha}(p)\in\TPrim$. Morevover $\tau_{i^\alpha}$ are orthogonal projectors so
  $\TPrim = \bigoplus_{\alpha, i} \tau_{i^\alpha}(\TPrim) = \bigoplus_{\alpha, i} Im(\tau_{i^\alpha})\cap\TPrim$.
\end{proof}





\subsubsection{The new basis  ~\texorpdfstring{$\OO$}{OO}}

\begin{defi} \label{O} Let $t_1, \ldots, t_k, r \in \PTreeblue$,
  $f_r \in \{[], [r]\}$ and $f_l = [\ell_1, \ldots, \ell_g] \in \PForestblue$,
  \begin{align}
    \label{OE}  \OO_{[]}  & \eqdef \QQ_\epsilon,\\
    \label{OF} \OO_{t_1, \ldots, t_k} & \eqdef \OO_{t_k} \valprec (\OO_{t_{k-1}} \valprec (\ldots  \valprec \OO_{t_1})\ldots),\\
    \label{OT} \OO_{\Nodeblue(i^\alpha, f_l = [\ell_1, \ldots, \ell_g], f_r)} & \eqdef \langle\OO_{\ell_1}, \OO_{\ell_2}, \ldots, \OO_{\ell_g}; \OO_{\Nodeblue(i^\alpha, [], f_r)}\rangle,\\
    \label{OP} \OO_{\Nodeblue(i^\alpha, [], f_r)} & \eqdef \Psi_{{i^\alpha}}(\OO_{f_r}).
  \end{align}
\end{defi}

\begin{ex}\label{ex:O}
  \begin{align*}
    \OO_{\scalebox{0.5}{\input{figures/arbres_pw/arbre34122left}}} =& \OO_{\scalebox{0.5}{{ \newcommand{\nodea}{\node[draw,ellipse] (a) {$1^\circ$}
;}\newcommand{\nodeb}{\node[draw,ellipse] (b) {$1^\circ$}
;}\begin{tikzpicture}[baseline={([yshift=-1ex]current bounding box.center)}]
\matrix[column sep=.3cm,row sep=.3cm,ampersand replacement=\&]{
         \& \nodea  \\ 
 \nodeb  \&         \\
};
\path[ultra thick,blue] (a) edge (b);
\end{tikzpicture}}}} \valprec \OO_{\scalebox{0.5}{\input{figures/arbres_pw/arbre122left}}} \\
    =& (\OO_{\scalebox{0.5}{{ \newcommand{\nodea}{\node[draw,ellipse] (a) {$1^\circ$}
;}\begin{tikzpicture}[baseline={([yshift=-1ex]current bounding box.center)}]
\matrix[column sep=.3cm,row sep=.3cm,ampersand replacement=\&]{
 \nodea  \\
};
\path[ultra thick,blue] ;
\end{tikzpicture}}}} \valsucc \OO_{\scalebox{0.5}{{ \newcommand{\nodea}{\node[draw,ellipse] (a) {$1^\circ$}
;}\begin{tikzpicture}[baseline={([yshift=-1ex]current bounding box.center)}]
\matrix[column sep=.3cm,row sep=.3cm,ampersand replacement=\&]{
 \nodea  \\
};
\path[ultra thick,blue] ;
\end{tikzpicture}}}} - \OO_{\scalebox{0.5}{{ \newcommand{\nodea}{\node[draw,ellipse] (a) {$1^\circ$}
;}\begin{tikzpicture}[baseline={([yshift=-1ex]current bounding box.center)}]
\matrix[column sep=.3cm,row sep=.3cm,ampersand replacement=\&]{
 \nodea  \\
};
\path[ultra thick,blue] ;
\end{tikzpicture}}}} \valprec \OO_{\scalebox{0.5}{{ \newcommand{\nodea}{\node[draw,ellipse] (a) {$1^\circ$}
;}\begin{tikzpicture}[baseline={([yshift=-1ex]current bounding box.center)}]
\matrix[column sep=.3cm,row sep=.3cm,ampersand replacement=\&]{
 \nodea  \\
};
\path[ultra thick,blue] ;
\end{tikzpicture}}}}) \valprec \Psi_{2^\bullet}(\OO_{\scalebox{0.5}{{ \newcommand{\nodea}{\node[draw,ellipse] (a) {$1^\circ$}
;}\newcommand{\nodeb}{\node[draw,ellipse] (b) {$1^\circ$}
;}\begin{tikzpicture}[baseline={([yshift=-1ex]current bounding box.center)}]
\matrix[column sep=.3cm,row sep=.3cm,ampersand replacement=\&]{
         \& \nodea  \\ 
 \nodeb  \&         \\
};
\path[ultra thick,blue] (a) edge (b);
\end{tikzpicture}}}})\\
    =& \quad \QQ_{34122} + \QQ_{24133} + \QQ_{14233} + \QQ_{43212} + \QQ_{42313} + \QQ_{41323}\\
                                                                    & - \QQ_{34212} - \QQ_{24313} - \QQ_{14323} - \QQ_{43122} - \QQ_{42133} - \QQ_{41233}
  \end{align*}
  More examples can be found in the annexes section with
  \cref{RPOQ123,matrOQ4}.
\end{ex}

\begin{theorem} \label{thm:O} For all $n \in \NN_{>0}$
  \begin{multicols}{2}
  \begin{enumerate}
  \item \label{OForest} $(\OO_f)_{f\in\PForestbluen}$ is a basis of $\WQSym_n$,
  \item \label{OTree} $(\OO_t)_{t\in\PTreebluen}$ is a basis of $\Prim_n$,
  \item \label{OParticular} $(\OO_t)_{t\in\PParticularbluen}$ is a basis of $\TPrim_n$.
  \end{enumerate}
  \end{multicols}
\end{theorem}

\begin{proof}
  The proof structure is the same as the one of \cref{thm:P} except for some
  statements that are exchanged as we can see in this table:
  \begin{center}
    \begin{tabular}{|l|l|l|}
      \hline
      \cref{RR_left_coprod} & \cref{QQ_left_coprod} & left corproduct in basis $\RR$ and $\QQ$.\\
      \hline
      \cref{rem:PParticularredn} & & $\PParticularredn = \bigsqcup_{I}\{\Nodered(I,[],f_r) | f_r \in \PForestrednmp^{I}\}$\\
      & \cref{rem:PParticularbluen} & $\PParticularbluen = \bigsqcup_{i,\alpha}\{\Nodeblue(i^\alpha,[],f_r) | f_r \in \PForestbluenmp^{i^\alpha}\}$.\\
      \hline
      \cref{im_eq_phi} & \cref{im_eq_psi} & $\Imm(\Phi_I) = \Imm(\tau_I)$ and $\Imm(\Psi_{i^\alpha}) = \Imm(\tau_{i^\alpha})$.\\
      \hline
      \cref{tot_stable_phi} & & $\TPrim^* = \bigoplus_{I} Im(\tau_{I}) \cap \TPrim^*\,.$\\
      & \cref{tot_stable_psi} & $\TPrim = \bigoplus_{\alpha, i} Im(\tau_{i^\alpha}) \cap \TPrim\,.$\\
      \hline
      \cref{thm:bij_primal} & \cref{thm:bij_dual} & $\dim(\TPrim^*_n) = \#\PParticularredn$ and $\dim(\TPrim_n) = \#\PParticularbluen$.\\
      \hline
    \end{tabular}
  \end{center}

\end{proof}

\begin{rem}\label{rem:Oske}
  The same remark as \cref{rem:Pske} can be done on the basis $\OO$. Indeed, it
  is defined with blue-packed forests ($\PForestblue$), it is nevertheless
  possible to use the blue-skeletons ($\PForestblueske$) or packed words($\PW$)
  thanks to the bijections of \cref{blue_ske_flat} and $\Fb$ of
  \cref{def:construction_blue_packed}.
\end{rem}


\section{Isomorphism between $\WQSym$ and $\WQSym^*$}
\label{sect4_bij}

According to \cref{prim_tot} $\WQSym$ (resp. $\WQSym^*$) is freely generated as
a dendriform algebra by $\TPrim$ (resp.  $\TPrim^*$). Therefore, any linear
isomorphism between $\TPrim$ and $\TPrim^*$ would lead to a bidendriform
isomorphism between $\WQSym$ and its dual. Thanks to the two bases $\PP$ and
$\OO$ any graded bijection between red-irreducible and blue-irreducible packed
words leads to such an isomorphism. We first make explicit how this is
done. Then the bijection is actually obtained as the restriction to
red-irreducibles of an involution on all packed words. The definition of the
bijection requires a new kind of forest mixing red and blue factorizations, namely
bicolored-packed forests. 

\subsection{A combinatorial solution to an algebraic problem} \label{sect4.1}

In this \cref{sect4.1}, we use the skeleton representation for bases $\PP$ and
$\OO$ as said in \cref{rem:Pske,rem:Oske}. Moreover we fix a graded bijection
$\bijirr$ between red-irreducible and blue-irreducible packed words.

\begin{defi}
Recall that $(\PP_t)_{t\in\PParticularredn}$ is a basis of $\TPrim^*_n$
(\cref{thm:P}) and $(\OO_t)_{t\in\PParticularbluen}$ is a basis of $\TPrim_n$
(\cref{thm:O}). By linearity, setting
  \begin{equation}
    t' \eqdef \scalebox{0.7}{{ \newcommand{\nodea}{\node[draw,ellipse] (a) {$\mu(v)$}
;}\begin{tikzpicture}[baseline={([yshift=-1ex]current bounding box.center)}]
\matrix[column sep=.3cm,row sep=.3cm,ampersand replacement=\&]{
 \nodea  \\
};
\path[ultra thick,red] ;
\end{tikzpicture}}} \in \PParticularblue,
    \qquad\text{and}\qquad
    \Bijirr(\PP_{t}) \eqdef \OO_{t'}.
  \end{equation}
  for all
  $t = \scalebox{0.7}{{ \newcommand{\nodea}{\node[draw,ellipse] (a) {$\,v\,$}
;}\begin{tikzpicture}[baseline={([yshift=-1ex]current bounding box.center)}]
\matrix[column sep=.3cm,row sep=.3cm,ampersand replacement=\&]{
 \nodea  \\
};
\path[ultra thick,red] ;
\end{tikzpicture}}} \in \PParticularred$,
  defines a linear isomorphism between the vector spaces $\TPrim^*_n$ and $\TPrim_n$
\end{defi}

\begin{defi}
  We define $\bijt$ as the extension of $\bijirr$ from red-skeleton to
  blue-skeleton forests by:
  \begin{align}
    \forall f = [t_1, \dots, t_k] \in \PForestred, \quad &\bijt(f) \eqdef [\bijt(t_1), \dots, \bijt(t_k)]\\
    \forall t = \Nodered(v, f_\ell, []) \in \PTreered, \quad &\bijt(t) \eqdef \Nodeblue(\bijirr(v),\bijt(f_\ell), []).
  \end{align}
\end{defi}

\begin{defi}
  We denote $\Bij$ the unique bidendriform isomorphism from $\WQSym^*$ to
  $\WQSym$ which verify for all $f \in \PForestred$:
  \begin{equation}
    \Bij(\PP_f) \eqdef \OO_{\bijt(f)}
  \end{equation}
\end{defi}

The existence and unicity are guaranteed by \cref{prim_tot}.

\begin{ex}
  \begin{align*}
    \bijt\left(\overbrace{\scalebox{0.7}{\input{figures/arbres_pw/arbre_bij_ex}}}^{f_{\red{R}}}\right)
    &= \overbrace{\scalebox{0.7}{\input{figures/arbres_pw/arbre_bij_ex_}}}^{{f_{\blue{B}}}}\\
    \Bij\left(\PP_{f_{\red{R}}}\right)
    &= \OO_{f_{\blue{B}}}\\
  \end{align*}
\end{ex}

In the following example, we take $\bijirr(212) \eqdef 122$ and
$\bijirr(w) \eqdef w$ for other words with a size less than $3$ as they are
simultaneously red-irreducible and blue-irreducible packed words.

Here are all red and blue-irreducible packed words of size 1, 2 and 3:
\[ \text{\red{Red} :\qquad}
  1,\quad 11,\quad
  132\ 121\ 212\ 111,
  \qquad\qquad
  \text{\blue{Blue} :\qquad}
  1,\quad 11,\quad
  132\ 121\ 122\ 111.
\]

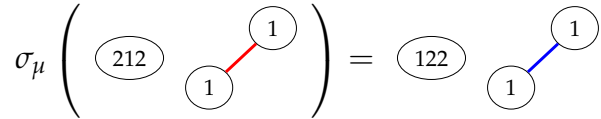
\begin{ex}
  \begin{align*}
    \bijt\left(\scalebox{0.7}{\input{figures/arbres_pw/arbre43412ske}}\right)
    &= \scalebox{0.7}{\input{figures/arbres_pw/arbre34122leftske}}\\
  \end{align*}
  These two forests are the same as those used in \cref{ex:P,ex:O}. So we have
  here the first example of the isomorphism from the basis $\RR$ to the basis
  $\QQ$:

  $\Bij\left(
    \begin{matrix}
      \RR_{14342} + \RR_{41342} + \RR_{43142} + \RR_{43412}\\
      - \RR_{24341} - \RR_{42341} - \RR_{43241} - \RR_{43421}\\
    \end{matrix}
  \right) =
  \left(
  \begin{matrix}
    \QQ_{34122} + \QQ_{24133} + \QQ_{14233} + \QQ_{43212}\\
    + \QQ_{42313} + \QQ_{41323} - \QQ_{34212} - \QQ_{24313}\\
    - \QQ_{14323} - \QQ_{43122} - \QQ_{42133} - \QQ_{41233}\\
  \end{matrix}
  \right)$
\end{ex}

We now have a construction of a bidendriform isomorphism for any graded
bijection~$\bijirr$ between red-irreducible and blue-irreducible packed words.

\subsection{Full decomposition of packed words into bicolored forests}

In order to define a bijection between red-irreducible and blue-irreducible
packed words, we need a new kind of forests that mixes up the red and blue
factorizations. More precisely, we will recursively alternate these
factorizations. We start with an unexpected lemma which implies that starting by
red or blue does not matter.

\begin{lem}\label{ins_insl}
  For all $a, b, c \in \PW$, with $c \neq \epsilon$, the following relations hold:
  \begin{equation}
    a \ins ( b \insl c) = b \insl (a \ins c) \text{\quad and\quad}
    a \ins 1 = a \insl 1,
  \end{equation}
  where $1$ is the packed word of size $1$.
  $$\scalebox{1}{\input{figures/box_d/ins_of_insl} \qquad\scalebox{2}{=}\qquad\quad \input{figures/box_d/insl_of_ins}}$$
\end{lem}

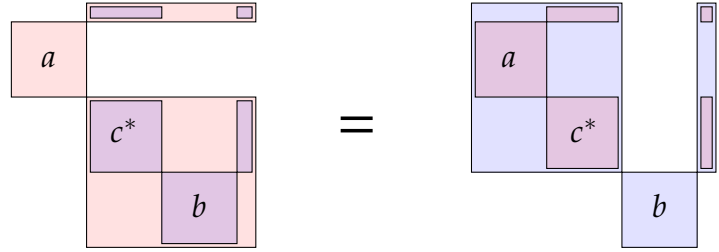
\begin{proof}
  Let~$a \in \PW$ and $1$ the packed word of size 1, by \cref{ins,insl}
  $a \ins 1 = a \dcdot 1$ and $a \insl 1 = a \dcdot 1$.

  Let~$a, b, c \in \PW$, with $c = c_1\cdots c_n$ of size $n>0$. We start by assuming
  that $c_n \neq \max(c)$ which implies that
  $c = \phi_I(\psi_{i^\alpha}(c^*)) = \psi_{i^\alpha}(\phi_I(c^*))$ with
  $I, i, \alpha, c^*$ unique by \cref{phi_bij,psi_bij}. With this
  relation we can deduce:
  \begin{align*}
    a \ins ( b \insl c) &= a \ins ( b \insl \psi_{i^\alpha}(\phi_{I}(c^*)))\\
                        &= a \ins ( \psi_{i+\max(b)^\alpha}(\phi_{I}(c^*) \gcdot b))\\
                        &= a \ins ( \psi_{i+\max(b)^\alpha}(\phi_{I}(c^* \gcdot b)))\\
                        &= a \ins ( \phi_{I}(\psi_{i+\max(b)^\alpha}(c^* \gcdot b)))\\
                        &= \phi_{I+|a|}(a \gcdot \psi_{i+\max(b)^\alpha}(c^* \gcdot b))\\
                        &= \phi_{I+|a|}(\psi_{i+\max(b)^\alpha}(a \gcdot c^* \gcdot b))
  \end{align*}
  \begin{align*}
    b \insl ( a \ins c) &= b \insl ( a \ins \phi_I(\psi_{i^\alpha}(c^*)))\\
                        &= b \insl ( \phi_{I+|a|}(a \gcdot \psi_{i^\alpha}(c^*)))\\
                        &= b \insl ( \phi_{I+|a|}(\psi_{i^\alpha}(a \gcdot c^*)))\\
                        &= b \insl ( \psi_{i^\alpha}(\phi_{I+|a|}(a \gcdot c^*)))\\
                        &= \psi_{i+\max(b)^\alpha}(\phi_{I+|a|}(a \gcdot c^*) \gcdot b)\\
                        &= \psi_{i+\max(b)^\alpha}(\phi_{I+|a|}(a \gcdot c^* \gcdot b))\\
                        &= \phi_{I+|a|}(\psi_{i+\max(b)^\alpha}(a \gcdot c^* \gcdot b)).
  \end{align*}
  The case where $c_n = \max(c)$ can be decomposed into different particular
  cases. In each of these cases, it is possible to find a relation with two
  different writings of $c$ that begin with $\phi$ or $\psi$ just like
  $c = \phi_I(\psi_{i^\alpha}(c^*)) = \psi_{i^\alpha}(\phi_I(c^*))$. These cases with the associated
  relation are:
  \begin{itemize}
  \item the case where $c$ is the packed word $1$ then $c = \phi_{1}(\epsilon) = \psi_{1^\circ}(\epsilon)$,
  \item the case where $c$ is of the form $c' \dcdot 1$ then
    $c = \phi_{1}(c') = \psi_{1^\circ}(c')$,
  \item the more general case where there is more than 1
    maximum including the one at the end then $c = \phi_{I\cdot c_n}(c^*) = \psi_{c_n^\bullet}(\phi_{I}(c^*))$.
  \end{itemize}
  In each of these cases it is possible to prove with a similar method that
  $a \ins ( b \insl c) = b \insl (a \ins c)$.
\end{proof}

\begin{ex}\label{ex:rel_L-alg}
  Here are some examples of this relation:
  \begin{align*}
    1 \ins (1 \insl 1) &= 213 \;\,\qquad= 1 \insl (1 \ins 1), \\
    11 \ins (12 \insl 2111) &= 44533123 = 12 \insl (11 \ins 2111), \\
    11 \ins (21 \insl 123) &= 5534216 \;\,= 21 \insl (11 \ins 123), \\
    1 \ins (112 \insl 3132) &= 56361124 = 112 \insl (1 \ins 3132).
  \end{align*}
\end{ex}

\begin{defi}\label{rb_br_fact}
  Let~$w$ be an irreducible packed word. Let $w = x \ins u$ be the
  red-factorization of $w$ and let $u = y \insl z$ be the blue-factorization of
  $u$. Then $w = x \ins (y \insl z)$ is called the
  \textbf{red-blue-factorization} of $w$. Symmetrically we define
  $w = y' \insl (x' \ins z')$ the \textbf{blue-red-factorization} of $w$.
\end{defi}

\begin{lem}\label{bi_fact}
  Let $w = x \ins (y \insl z)$ be the red-blue-factorization and let
  $w = y' \insl (x' \ins z')$ be the blue-red-factorization of an irreducible
  packed word $w$.

  With these two factorizations, we have that $z = z'$ and it is both
  red-irreducible and blue-irreducible packed word. Moreover,
  \begin{itemize}
  \item either $z = z' = 1$, $y = x' = \epsilon$ and $x = y'$
  \item or $x = x'$, $y = y'$.
  \end{itemize}
\end{lem}

\begin{ex}
  Here are some examples of red-blue-factorization and blue-red-factorization:
  \begin{align*}
    12 \ins (\epsilon \insl 1) &= 213 \;\,\qquad= 12 \insl (\epsilon \ins 1), \\
    11 \ins (12 \insl 1211) &= 44353123 = 12 \insl (11 \ins 1211), \\
    553421 \ins (\epsilon \insl 1) &= 5534216 \;\,= 553421 \insl (\epsilon \ins 1), \\
    1 \ins (112 \insl 3132) &= 56361124 = 112 \insl (1 \ins 3132).
  \end{align*}
\end{ex}

\begin{proof}
  We start by prooving the case where $z = z' = 1$, $y = x' = \epsilon$ and
  $x = y'$.  Let~$w'$ be an irreducible packed word and $w = w' \dcdot 1$. We
  have that $w = w' \ins 1$ is the red-factorization of $w$ and $w = w' \insl 1$
  is the blue-factorization of $w$. In this case we immediately have that
  $w = w' \ins (\epsilon \insl 1)$ is the red-blue-factorization of $w$ and that
  $w = w' \insl (\epsilon \ins 1)$ is the blue-red-factorization of $w$.

  Now let $w$ be an irreducible packed word of size $n$ that cannot be written
  as $w' \dcdot 1$. In other words, there is a maximum strictly before the last
  letter of $w$ ($\exists i<n, w_i = \max(w)$).

  We define the two sets of triplet of packed words that verify equations of the
  factorizations for $w$: 
  \begin{align*}
    S_{\red{R}\blue{B}}(w) &\eqdef \{(a,b,c) \in \PW, c \neq \epsilon, w = a \ins (b \insl c) \},\\
    S_{\blue{B}\red{R}}(w) &\eqdef \{(a,b,c) \in \PW, c \neq \epsilon, w = b \insl (a \ins c) \}.
  \end{align*}
  Thanks to \cref{ins_insl} these two sets are equal, we define
  $S(w) \eqdef S_{\red{R}\blue{B}}(w) = S_{\blue{B}\red{R}}(w)$.

  In the red-blue-factorization $w = x \ins (y \insl z)$, we maximize the size
  of $x$, then we maximize the size of $y$ in the remaining word. In the
  blue-red-factorization we commute the order of maximizations. We will
  caracterize $S(w)$ and see the limit of the two maximizations to prove that
  they can commute.

  Let $w^*$ be the packed word comming from $w$ where the last letter and all
  occurences of the maximum are removed and let
  $I = [i_1, \dots, i_p], i, \alpha$ such that
  $w = \psi_{i^\alpha}(\phi_I(w^*))$. By hypothesis, we have that
  $I \neq \emptyset$. Let $w^* = w_1\gcdot \cdots \gcdot w_k$ be the global descent
  decomposition of $w^*$. Let $\ell$ be the maximum such that
  $|w_1\gcdot\cdots \gcdot w_{\ell}| \leq i_1$, $i_1$ being the position of the first
  maximum of $w$. Let $r$ be the minimum such that
  $|w_{r}\gcdot\cdots \gcdot w_{k}| \leq n - i_p - 1$, $i_p$ being the position of the
  last maximum of $w$ before the last letter. As $i_1 \leq i_p$ by definition, we
  have that $\ell < r$. We can caracterize the set $S(w)$:
  \begin{align*}
    S(w) = \{(&a = w_1\gcdot \cdots \gcdot w_{r^0}, \\
              &b = w_{\ell^0}\gcdot \cdots \gcdot w_{k}, \\
              &c = \psi_{(i-\max(b))^\alpha}(\phi_{I-|a|}(w_{r^0+1}\gcdot \cdots \gcdot w_{\ell^0-1}))),\\
              & \text{with }r^0\leq r \text{ and } \ell \leq \ell^0\}.\qedhere
  \end{align*}




\end{proof}

Here are all packed words that are both red-irreducible and blue-irreducible of
size less than $4$:
  \[
    1,\quad 11,\quad
    111\ 121\ 132,
  \]
  \[
    1111\ 1121\ 1132\ 1211\ 1212\ 1221\ 1231\ 1232\ 1243\ 1312\ 1321\
  \]
  \[
    1322\ 1323\ 1332\ 1342\ 1423\ 1432\ 2121\ 2122\ 2132\ 2143\ 3132\
  \]

\begin{table}[!h]
  \[
    \begin{tabular}{|c|c|c|c|c|c|c|c|c|c||l|}
      \hline
      $n$ & 1 & 2 & 3 & 4 & 5 & 6 & 7 & 8 & 9 \\
      \hline
      $\Si_n\in\PW_n$ & 1 & 1 & 3 & 22 & 196 & 2 008 & 23 184 & 297 456 & 4 199 216 \\
      \hline
      $\Si_n\in\Sn$ & 1 & 0 & 1 & 5 & 32 & 236 & 1 951 & 17 827 & 178 418 \\
      \hline
    \end{tabular} 
  \]
  \caption{Number of both red-and-blue-irreducible packed words and permutations.}
  \label{oeis_rbirr}
\end{table}

Recall our notations for Hilbert series of an algebra $A$,
$\SA(z) \eqdef \sum_{n=1}^{+\infty}\dim(A_n)z^n$,
$\SP(z) \eqdef \sum_{n=1}^{+\infty}\dim(\Primof(A_n))z^n$ and
$\ST(z) \eqdef \sum_{n=1}^{+\infty}\dim(\TPrimof(A_n))z^n.$

Recall the relations between these series:
$$\SP = \SA / (1 + \SA) \text{ or equivalently } \SA = \SP / (1 - \SP) \text{ (see~\cref{left_prod}),}$$
$$\ST = \SA / (1 + \SA)^2 \text{ or equivalently } \SP = \ST (1 + \SA) \text{ (see~\cref{brace}).}$$

If we define the serie $\SI = \sum_{n=1}^{+\infty}\Si_nz^n$ where $\Si_n$ is the number
of both red-and-blue-irreducible words of size $n$, then we have the following
relation:
$$\SI = \SA / (1 + \SA)^3 + z \SP \text{ or equivalently } \ST = (\SI - z)(1 + \SA) + z.$$

\bigskip

So far we have seen red-biplane trees and blue-biplane trees. In this section we
define red-blue-biplane trees and blue-red-biplane trees, the edges of these
trees are of two different colors and the labels are red-and-blue-irreducible
packed words. We denote by $\Noderb(x, f_\ell, f_r)$ (resp.
$\Nodebr(x, f_\ell, f_r)$) the biplane tree whose edges between the root and the
left forest $f_\ell$ are \red{red} (resp. \blue{blue}) and edges between the root
and the right forest $f_r$ are \blue{blue} (resp. \red{red}).



\begin{defi}\label{def:construction_bicolor}
  The bicolored forests~$\Frb(w)$ and $\Fbr(w)$~(resp. trees~$\Trb(w)$
  and~$\Tbr(w)$) associated to a packed word (resp. irreducible packed word) $w$
  are defined in a mutual recursive way as follows:
  \begin{itemize}
  \item $\Frb(\epsilon) = \Fbr(\epsilon) = []$ (empty forest),
  \item for any packed word $w$, let $ w = w_1\gcdot w_2\gcdot\dots\gcdot w_k$
    be the global descent decomposition, then
    $\Frb(w) \eqdef [\Trb(w_1), \Trb(w_2), \dots, \Trb(w_k)],$
    $$\scalebox{0.7}{\input{figures/box_d/gd_fact} \qquad\scalebox{2}{$\to$}\qquad\quad
      $\overbrace{\input{figures/arbres_pw/arbrew1_rb}}^{\Trb(w_1)}$
      $\overbrace{{ \newcommand{\nodea}{\node[draw,ellipse] (a) {}
;}\newcommand{\nodeb}{\node[draw,ellipse] (b) {}
;}\newcommand{\nodec}{\node[draw,ellipse] (c) {}
;}\newcommand{\noded}{\node[draw,ellipse] (d) {}
;}\newcommand{\nodee}{\node[draw,ellipse] (e) {}
;}\begin{tikzpicture}[baseline={([yshift=-1ex]current bounding box.center)}]
\matrix[column sep=.3cm,row sep=.3cm,ampersand replacement=\&]{
         \&         \& \nodea \& \& \\ 
 \nodeb  \& \nodec  \&        \& \noded \& \nodee \\
};
\path[ultra thick,red] (a) edge (b) edge (c);
\path[ultra thick,blue] (a) edge (d) edge (e);
\end{tikzpicture}}}^{\Trb(w_{2})}$
      $\cdots$ $\overbrace{{ \newcommand{\nodea}{\node[draw,ellipse] (a) {}
;}\newcommand{\nodeb}{\node[draw,ellipse] (b) {}
;}\newcommand{\nodec}{\node[draw,ellipse] (c) {}
;}\newcommand{\noded}{\node[draw,ellipse] (d) {}
;}\newcommand{\nodee}{\node[draw,ellipse] (e) {}
;}\newcommand{\nodef}{\node[draw,ellipse] (f) {}
;}\begin{tikzpicture}[baseline={([yshift=-1ex]current bounding box.center)}]
\matrix[column sep=.3cm,row sep=.3cm,ampersand replacement=\&]{
         \& \nodea \&        \&        \&        \&\\ 
 \nodeb  \&        \&        \& \nodec \& \noded \& \\
         \&        \& \nodee \&        \&        \& \nodef\\
};
\path[ultra thick,red] (a) edge (b)
(c) edge (e);
\path[ultra thick,blue] (a) edge (c) edge (d)
(d) edge (f);
\end{tikzpicture}}}^{\Trb(w_k)}$}.$$
  \item for any packed word $w$, let $ w = w_1\gcdot w_2\gcdot\dots\gcdot w_k$
    be the global descent decomposition, then
    $\Fbr(w) \eqdef [\Tbr(w_k), \Tbr(w_{k-1}), \dots, \Tbr(w_1)].$\\
    $$\scalebox{0.7}{\input{figures/box_d/gd_fact} \qquad\scalebox{2}{$\to$}\qquad\quad
      $\overbrace{{ \newcommand{\nodea}{\node[draw,ellipse] (a) {}
;}\newcommand{\nodeb}{\node[draw,ellipse] (b) {}
;}\newcommand{\nodec}{\node[draw,ellipse] (c) {}
;}\newcommand{\noded}{\node[draw,ellipse] (d) {}
;}\newcommand{\nodee}{\node[draw,ellipse] (e) {}
;}\newcommand{\nodef}{\node[draw,ellipse] (f) {}
;}\begin{tikzpicture}[baseline={([yshift=-1ex]current bounding box.center)}]
\matrix[column sep=.3cm,row sep=.3cm,ampersand replacement=\&]{
         \&        \&        \&        \& \nodea \& \\
         \& \noded \& \nodec \&        \&        \& \nodeb \\
 \nodef  \&        \&        \& \nodee \&        \& \\
};
\path[ultra thick,red] (a) edge (b)
(c) edge (e);
\path[ultra thick,blue] (a) edge (c) edge (d)
(d) edge (f);
\end{tikzpicture}}}^{\Tbr(w_k)}$
      $\cdots$ $\overbrace{{ \newcommand{\nodea}{\node[draw,ellipse] (a) {}
;}\newcommand{\nodeb}{\node[draw,ellipse] (b) {}
;}\newcommand{\nodec}{\node[draw,ellipse] (c) {}
;}\newcommand{\noded}{\node[draw,ellipse] (d) {}
;}\newcommand{\nodee}{\node[draw,ellipse] (e) {}
;}\begin{tikzpicture}[baseline={([yshift=-1ex]current bounding box.center)}]
\matrix[column sep=.3cm,row sep=.3cm,ampersand replacement=\&]{
         \&         \& \nodea \& \& \\ 
 \nodeb  \& \nodec  \&        \& \noded \& \nodee \\
};
\path[ultra thick,blue] (a) edge (b) edge (c);
\path[ultra thick,red] (a) edge (d) edge (e);
\end{tikzpicture}}}^{\Tbr(w_{2})}$
      $\overbrace{\input{figures/arbres_pw/arbrew1_br}}^{\Tbr(w_1)}$}.$$
    (notice the same inversion as
    in~\cref{def:construction_blue_skeleton,def:construction_blue_packed} for
    $\Fbr$.)
  \item for any irreducible packed word $w$, let
    $w = a \ins (b \insl c)$ be the red-blue-factorization, then
    $\Trb(w) \eqdef \Noderb(c, \Frb(a), \Frb(b)).$
    $$\scalebox{0.7}{\input{figures/box_d/ins_of_insl} \qquad\scalebox{2}{$\to$}\qquad\quad
      \input{figures/arbres_pw/arbre_rb_F}}.$$
  \item for any irreducible packed word $w$, let
    $w = b \insl (a \ins c)$ be the blue-red-factorization, then
    $\Tbr(w) \eqdef \Nodebr(c, \Fbr(b), \Fbr(a)).$
    $$\scalebox{0.7}{\input{figures/box_d/insl_of_ins} \qquad\scalebox{2}{$\to$}\qquad\quad
      \input{figures/arbres_pw/arbre_br_F}}.$$
  \end{itemize}
\end{defi}
\begin{ex}
  For this example we write the word $w$ in hexadecimal in order to have a big
  and clear example, let $w = DDDCCCEBBE9FA587653213449$. The word $w$ is
  irreducible so there is only one tree in the forest.
  To have $\Trb(w)$, we start by the blue-red-factorization,
  \begin{equation*}
    w = \redul{DDDCCCEBBE}\,{9FA}\,\blueul{58765321344}\,{9} = \redul{3332224114} \ins (\blueul{58765321344} \insl {1321}).
  \end{equation*}
  Then we decompose each sub-word according to their global descents and do
  blue-red-factorizations recursively until we have only both red-irreducible
  and blue-irreducible packed words:
  \begin{align*}
    w &= \redul{3332224114} \ins (\blueul{\underline{58765}\,\underline{321344}} \insl {1321}),\\
      &= \redul{\redul{333222}\,{4}\,\blueul{11}\,{4}} \ins (\blueul{(\underline{{14321}} \gcdot \underline{\redul{3213}\,{44}})} \insl {1321}),\\
      &[\dots]\\
      &= ((111 \gcdot 111) \ins (11 \insl 11)) \ins ((14321\gcdot(((1 \gcdot 1)\insl 11)\ins 11))\insl 1321)
  \end{align*}
  $$\scalebox{0.3}{\input{figures/box_d/pw_rb_ex}}
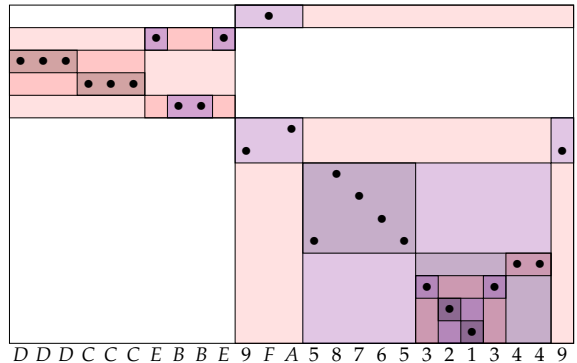
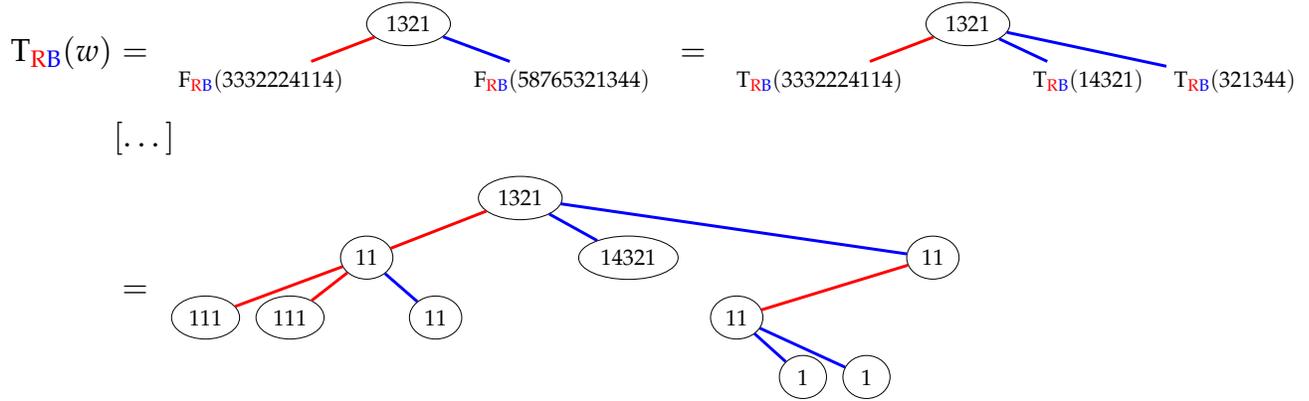
$$
  \begin{align*}
    \Trb(w) & = \scalebox{0.7}{\input{figures/arbres_pw/arbre_construct_DDD1}} = \scalebox{0.7}{\input{figures/arbres_pw/arbre_construct_DDD2}}\\
    &[\dots]\\
            & = \scalebox{0.7}{\input{figures/arbres_pw/arbre13,13,13,12,12,12,14,11,11,14,9,15,10,5,8,7,6,5,3,2,1,3,4,4,9ske_bi}}\\
  \end{align*}
  More examples can be found in the annexes section with
  \cref{F(w)123,F(1234),F(1233),F(1223),F(1123),F(1111)}.
\end{ex}





\begin{defi}\label{def:Frb}
  There are two types of bicolored trees, the only difference is that colors
  \red{red} and \blue{blue} are inverted. Let $t$ be a labeled biplane tree. We
  write $t = \Noderb(w, f_\ell, f_r)$ where ${w \in \PW}$,
  $f_\ell = [\ell_1, \dots, \ell_g]$ is the left forest of $t$ and
  $f_r = [r_1, \dots, r_d]$ is the right forest of $t$. We depict $t$ as
  follows:
  \begin{equation*}
    t = \scalebox{0.8}{{ \newcommand{\nodea}{\node[draw,ellipse] (a) {$w$}
;}\newcommand{\nodeb}{\node (b) {$\ell_1$}
;}\newcommand{\nodebc}{\node (bc) {$\ldots$}
;}\newcommand{\nodec}{\node (c) {$\ell_g$}
;}\newcommand{\noded}{\node (d) {$r_1$}
;}\newcommand{\nodede}{\node (de) {$\ldots$}
;}\newcommand{\nodee}{\node (e) {$r_d$}
;}
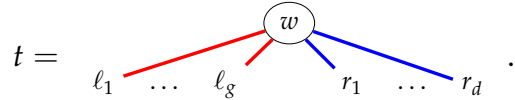
\begin{tikzpicture}[auto]
\matrix[column sep=.3cm, row sep=.3cm,ampersand replacement=\&]{
      \&  \&  \& \nodea  \& \\ 
 \nodeb  \& \nodebc \& \nodec \&  \& \noded \& \nodede \& \nodee \\
};

\path[ultra thick, red] (a) edge (b) edge (c);
\path[ultra thick, blue] (a) edge (d) edge (e);

\end{tikzpicture}}}.
  \end{equation*}

  We say that $t$ is a \textbf{red-blue-packed tree} if it
  satisfies:
  \begin{equation*}
    \left\{
      \begin{array}{l}
        w = 1,\\ 
        f_r = [],\\
        f_\ell \text{ is a red-blue-packed forest}.
      \end{array}\right.
    \text{ or }
    \left\{
      \begin{array}{l}
        w \neq 1 \text{ is red-irreducible and blue-irreducible},\\
        f_r \text{ and } f_\ell \text{ are red-blue-packed forest}.
      \end{array}\right.
  \end{equation*}
\end{defi}

\begin{defi}
  The \textbf{weight} of a bicolored-packed tree is the sum of the size of
  packed words in the nodes.
\end{defi}

We have already done it four times
(\cref{def:construction_red_skeleton_star,def:Fr-1,def:construction_blue_skeleton_star,def:Fb-1})
and we will do it one last time, to prove that $\Frb, \Trb, \Fbr, \Tbr$ are
bijections, we define $\Frbstar, \Trbstar, \Fbrstar, \Tbrstar$ and prove that
they are the inverse maps.

\begin{defi}\label{def:Frb-1}
  We define here the maps $\Frbstar, \Trbstar, \Fbrstar, \Tbrstar$ that
  transforms bicolored-packed forests and trees into packed words. We reverse all
  instructions of \cref{def:Frb} as follows:
  \begin{itemize}
  \item $\Frbstar([]) = \Fbrstar([]) = \epsilon$ (empty packed word),
  \item for any red-blue-packed forest $f = [t_1, t_2\dots, t_k]$, we have\\
    $\Frbstar(f) \eqdef [\Trbstar(t_1), \Trbstar(t_2), \dots, \Trbstar(t_k)].$
    $$\scalebox{0.7}{$\overbrace{\input{figures/arbres_pw/arbrew1_rb}}^{t_1}$
      $\overbrace{{ \newcommand{\nodea}{\node[draw,ellipse] (a) {}
;}\newcommand{\nodeb}{\node[draw,ellipse] (b) {}
;}\newcommand{\nodec}{\node[draw,ellipse] (c) {}
;}\newcommand{\noded}{\node[draw,ellipse] (d) {}
;}\newcommand{\nodee}{\node[draw,ellipse] (e) {}
;}\begin{tikzpicture}[baseline={([yshift=-1ex]current bounding box.center)}]
\matrix[column sep=.3cm,row sep=.3cm,ampersand replacement=\&]{
         \&         \& \nodea \& \& \\ 
 \nodeb  \& \nodec  \&        \& \noded \& \nodee \\
};
\path[ultra thick,red] (a) edge (b) edge (c);
\path[ultra thick,blue] (a) edge (d) edge (e);
\end{tikzpicture}}}^{t_{2}}$
      $\cdots$ $\overbrace{{ \newcommand{\nodea}{\node[draw,ellipse] (a) {}
;}\newcommand{\nodeb}{\node[draw,ellipse] (b) {}
;}\newcommand{\nodec}{\node[draw,ellipse] (c) {}
;}\newcommand{\noded}{\node[draw,ellipse] (d) {}
;}\newcommand{\nodee}{\node[draw,ellipse] (e) {}
;}\newcommand{\nodef}{\node[draw,ellipse] (f) {}
;}\begin{tikzpicture}[baseline={([yshift=-1ex]current bounding box.center)}]
\matrix[column sep=.3cm,row sep=.3cm,ampersand replacement=\&]{
         \& \nodea \&        \&        \&        \&\\ 
 \nodeb  \&        \&        \& \nodec \& \noded \& \\
         \&        \& \nodee \&        \&        \& \nodef\\
};
\path[ultra thick,red] (a) edge (b)
(c) edge (e);
\path[ultra thick,blue] (a) edge (c) edge (d)
(d) edge (f);
\end{tikzpicture}}}^{t_k}$ \qquad\scalebox{2}{$\to$}\qquad\quad
      \input{figures/box_d/gd_fact}} \text{ with $w_i = \Trbstar(t_i)$}.$$
  \item for any blue-red-packed forest $f = [t_1, t_2\dots, t_k]$, we have\\
    $\Fbrstar(f) \eqdef [\Tbrstar(t_k), \Tbrstar(t_{k-1}), \dots, \Tbrstar(t_1)].$
    $$\scalebox{0.7}{
      $\overbrace{{ \newcommand{\nodea}{\node[draw,ellipse] (a) {}
;}\newcommand{\nodeb}{\node[draw,ellipse] (b) {}
;}\newcommand{\nodec}{\node[draw,ellipse] (c) {}
;}\newcommand{\noded}{\node[draw,ellipse] (d) {}
;}\newcommand{\nodee}{\node[draw,ellipse] (e) {}
;}\newcommand{\nodef}{\node[draw,ellipse] (f) {}
;}\begin{tikzpicture}[baseline={([yshift=-1ex]current bounding box.center)}]
\matrix[column sep=.3cm,row sep=.3cm,ampersand replacement=\&]{
         \&        \&        \&        \& \nodea \& \\
         \& \noded \& \nodec \&        \&        \& \nodeb \\
 \nodef  \&        \&        \& \nodee \&        \& \\
};
\path[ultra thick,red] (a) edge (b)
(c) edge (e);
\path[ultra thick,blue] (a) edge (c) edge (d)
(d) edge (f);
\end{tikzpicture}}}^{t_1}$
      $\cdots$ $\overbrace{{ \newcommand{\nodea}{\node[draw,ellipse] (a) {}
;}\newcommand{\nodeb}{\node[draw,ellipse] (b) {}
;}\newcommand{\nodec}{\node[draw,ellipse] (c) {}
;}\newcommand{\noded}{\node[draw,ellipse] (d) {}
;}\newcommand{\nodee}{\node[draw,ellipse] (e) {}
;}\begin{tikzpicture}[baseline={([yshift=-1ex]current bounding box.center)}]
\matrix[column sep=.3cm,row sep=.3cm,ampersand replacement=\&]{
         \&         \& \nodea \& \& \\ 
 \nodeb  \& \nodec  \&        \& \noded \& \nodee \\
};
\path[ultra thick,blue] (a) edge (b) edge (c);
\path[ultra thick,red] (a) edge (d) edge (e);
\end{tikzpicture}}}^{t_{k-1}}$
      $\overbrace{\input{figures/arbres_pw/arbrew1_br}}^{t_k}$ \qquad\scalebox{2}{$\to$}\qquad\quad
      \input{figures/box_d/gd_fact_left}} \text{ with $w_i = \Tbrstar(t_i)$}.$$
  \item for any red-blue-packed tree $t = \Noderb(c, f_\ell, f_r)$, we have\\
    $\Trbstar(t) \eqdef \Frbstar(f_\ell) \ins (\Frbstar(f_r) \insl c).$
    $$\scalebox{0.7}{ \input{figures/arbres_pw/arbre_rb}\qquad\scalebox{2}{$\to$}\qquad\quad
      \input{figures/box_d/ins_of_insl}}
    \begin{matrix}
      \text{with}&a = \Frbstar(f_\ell)\\
      \text{and}&b = \Frbstar(f_r)
    \end{matrix}.
    $$
  \item for any blue-red-packed tree $t = \Nodebr(c, f_\ell, f_r)$, we have\\
    $\Tbrstar(t) \eqdef \Fbrstar(f_r) \insl (\Fbrstar(f_\ell) \ins c).$
    $$\scalebox{0.7}{\input{figures/arbres_pw/arbre_br} \qquad\scalebox{2}{$\to$}\qquad\quad
      \input{figures/box_d/insl_of_ins}}
    \begin{matrix}
      \text{with}&a = \Frbstar(f_r)\\
      \text{and}&b = \Frbstar(f_\ell)
    \end{matrix}.$$
  \end{itemize}
\end{defi}

\begin{theorem}\label{bij_bicolored} 
  The maps $\Frb$ and $\Frbstar$ (resp. $\Trb$ and $\Trbstar$) are two
  converse bijections between packed words of size $n$ and red-blue-packed
  forests (resp. irreducible packed words and red-blue-packed trees)
  of weight $n$. That is to say $\Frbinv = \Frbstar$ and $\Trbinv =
  \Trbstar$. We have the same result with inversions of red and blue.
\end{theorem}

\begin{proof}
  It is simple to prove by induction on the size of the trees that domain and
  codomain are as announced and that the functions are inverse to each
  other. The proof is similar to the proofs of
  \cref{lem:red_skeleton_bij,bij_red,lem:blue_skeleton_bij,bij_blue} using
  \cref{def:construction_bicolor} of $\Frb, \Trb, \Fbr, \Tbr$, \cref{def:Frb} of
  bicolored-packed forests and trees and \cref{def:Frb-1} of
  $\Frbstar, \Trbstar, \Fbrstar, \Tbrstar$.
\end{proof}

\begin{rem}\label{rem:OP_bicolored}
  We now have two new families of forests $\PForestrb$ and $\PForestbr$ that are
  in bijection with packed words and therefore in bijection with red-packed and
  blue-packed forests. As in \cref{rem:Pske,rem:Oske}, this gives us two other
  way to index bases $\OO$ and $\PP$ of $\WQSym$ and $\WQSym^*$.
\end{rem}

\subsection{An involution on packed words}

We are now in position to define a bijection between red-irreducible and
blue-irreducible packed words. This bijection is actually the restriction of an
involution defined on all packed words. Precisely, we will define two
transformations on bicolored forests

We need to define the notion of mirror transformation of bicolored-packed
forests and trees. This transformation is defined from a red-blue to blue-red or
from blue-red to red-blue, so in the notations we will use $XY$ instead of
$\red{R}\blue{B}$ or $\blue{B}\red{R}$ to point out where the swap is made.
\begin{defi}\label{mirror}
  The mirror transformation of a bicolored-packed forest $f = [t_1, \dots, t_k]$
  is given by $\widetilde{f} \eqdef [\widetilde{t_k}, \dots, \widetilde{t_1}]$
  where $\widetilde{t_i}$ is the mirror transformation of $t_i$ recursively defined as
  follows. For any $t = \Nodexy(z, f_\ell, f_r)$ then
  $$
  \widetilde{t} \eqdef \left\{
    \begin{array}{ll}
        \Nodeyx(z, \widetilde{f_r},  \widetilde{f_\ell})& \mbox{if } z \neq 1,\\
        \Nodeyx(1, \widetilde{f_\ell}, []) & \mbox{if } z = 1.
    \end{array}\right.
  $$
\end{defi}
Note that when $z\neq1$, the left and right forests are swapped whereas they are
not when $z=1$. But in the latter case, we have necessarily
$f_r = \widetilde{f_r} = []$. These two cases correspond to the two cases of
\cref{def:Frb} so the mirror transformation of a red-blue-packed forest is
indeed a blue-red-packed forest.

\begin{ex}
  Here are two examples of mirror transformations.
  \begin{align*}
    \scalebox{0.6}{\input{figures/arbres_pw/arbre55632124ske_bi}}
    &\text{\quad and \quad}
      \scalebox{0.6}{\input{figures/arbres_pw/arbre55632124leftske_bi}}\\
    \scalebox{0.47}{\input{figures/arbres_pw/arbre13,13,13,12,12,12,14,11,11,14,9,15,10,5,8,7,6,5,3,2,1,3,4,4,9ske_bi}}
    &\text{\quad and \quad}
      \scalebox{0.47}{{ \newcommand{\nodea}{\node[draw,ellipse] (a) {$1321$}
;}\newcommand{\nodeb}{\node[draw,ellipse] (b) {$11$}
;}\newcommand{\nodec}{\node[draw,ellipse] (c) {$11$}
;}\newcommand{\noded}{\node[draw,ellipse] (d) {$\,1\,$}
;}\newcommand{\nodee}{\node[draw,ellipse] (e) {$\,1\,$}
;}\newcommand{\nodef}{\node[draw,ellipse] (f) {$14321$}
;}\newcommand{\nodeg}{\node[draw,ellipse] (g) {$11$}
;}\newcommand{\nodeh}{\node[draw,ellipse] (h) {$11$}
;}\newcommand{\nodei}{\node[draw,ellipse] (i) {$111$}
;}\newcommand{\nodej}{\node[draw,ellipse] (j) {$111$}
;}
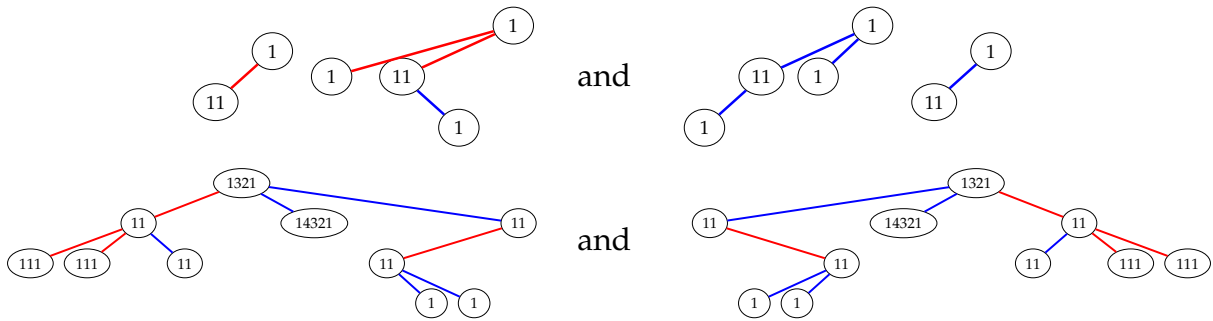
\begin{tikzpicture}[baseline={([yshift=-1ex]current bounding box.center)}]
\matrix[column sep=.3cm,row sep=.3cm,ampersand replacement=\&]{
         \&         \&         \&         \&         \&         \& \nodea  \&         \&         \&         \&         \\ 
         \& \nodeb  \&         \&         \&         \& \nodef  \&         \&         \& \nodeg  \&         \&         \\ 
         \&         \&         \&         \& \nodec  \&         \&         \& \nodeh  \&         \& \nodei  \& \nodej  \\ 
         \&         \& \noded  \& \nodee  \&         \&         \&         \&         \&         \&         \&         \\
};
\path[ultra thick,blue] (c) edge (d) edge (e)
	(b)
	(g) edge (h)
	(a) edge (b) edge (f);

\path[ultra thick,red] (b) edge (c)
	(g) edge (i)
	(g) edge (j)
	(a) edge (g);
\end{tikzpicture}}}
  \end{align*}
\end{ex}

\begin{prop}\label{lem:mirror}
  For all packed words $w$, both associated bicolored-packed forests $\Frb(w)$
  and $\Fbr(w)$ are mirror image of each other.
\end{prop}

\begin{proof} 
  The proof is a computation of mirror transformation (\cref{mirror}) on each
  items of \cref{def:construction_bicolor} of $\Frb$. For the first three items
  the computation of
  $\widetilde{f} \eqdef [\widetilde{t_k}, \dots, \widetilde{t_1}]$ is
  sufficient. Thanks to the relation of \cref{bi_fact}
  ($x \ins (y \insl z) = y \insl (x \ins z)$ in the case $z \neq 1$) the two
  remaining items are also simple computation of $\widetilde{t}$.
\end{proof}

\begin{defi}
  The color swap of a bicolored-packed forest $f = [t_1, \dots, t_k]$ is given
  by $\widehat{f} \eqdef [\widehat{t_k}, \dots, \widehat{t_1}]$ where
  $\widehat{t_i}$ is the color swap of $t_i$ recursively defined as follows. For
  any $t = \Nodexy(z, f_\ell, f_r)$ then
  $\widehat{t} \eqdef \Nodeyx(z, \widehat{f_\ell}, \widehat{f_r})$.

  In other words, it is a recoloration of each edges using the other
  color. Every blue edges become red and vice versa.
\end{defi}

\begin{ex}
  Here are two examples of color swaps.
  \begin{align*}
    \scalebox{0.6}{\input{figures/arbres_pw/arbre55632124ske_bi}}
    &\text{\quad and \quad}
      \scalebox{0.6}{\input{figures/arbres_pw/arbre45536112leftske_bi}}\\
    \scalebox{0.47}{\input{figures/arbres_pw/arbre13,13,13,12,12,12,14,11,11,14,9,15,10,5,8,7,6,5,3,2,1,3,4,4,9ske_bi}}
    &\text{\quad and \quad}
      \scalebox{0.47}{\input{figures/arbres_pw/arbre14,12,11,13,13,14,7,10,9,8,7,5,15,6,3,3,4,2,2,2,1,1,1,4,5leftske_bi}}
  \end{align*}
  More examples can be found in the annexes section with \cref{RPOQ123}.
\end{ex}

When we focus on the packed words associated to these forest, the color swap
correspond to the swap of the two operations $\ins$ and $\insl$ in a
bicolored-factorization. More precisely, if $w$ is an irreducible packed word
and $w = x \ins (y \insl z)$ is the red-blue-factorization of $w$, then the
color swap on the associated forest correspond to $w' = x \insl (y \ins z)$.

\begin{lem}\label{lem:commut}
  Mirror transformation and color swap commute. It means that for all
  bicolored-packed forest $f$, we have
  $\widehat{\left(\widetilde{f}\right)} = \widetilde{\left(\widehat{f}\right)}$.
\end{lem}

\begin{proof}
  The proof is immediate. Indeed, the definition of mirror transformation is
  independant of color swap and symmetrically, the color swap is independant of
  the tree shape.
\end{proof}

\begin{coro}\label{coro:commut_diag}
  The diagram on \cref{commut_diag} is commutative.
  \begin{figure}[!h]
    $$\input{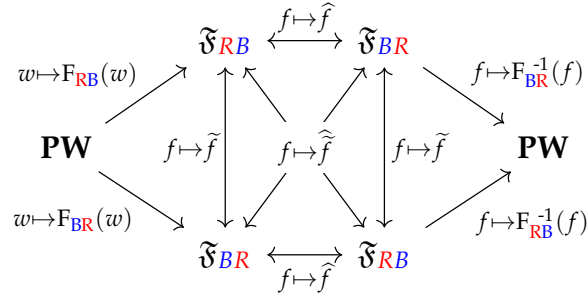}$$
    \caption{Commutative diagram of maps on bicolored-packed forests}
    \label{commut_diag}
  \end{figure}
  So $\widehat{w}:= \Frbinv(\widehat{\Fbr(w)})$ is an involution on packed words.
\end{coro}

\begin{proof}
  Thanks to \cref{lem:mirror,lem:commut} the diagram is immediately
  commutative. The mirror transformation and the color swap are independant
  involutions so the conjunction is an involution.
\end{proof}

\begin{coro}\label{coro:red_to_blue}
  The application $w \mapsto \widehat{w}$ send blue (resp. red) irreducibles packed
  words to red (resp. blue) irreducibles packed words.
\end{coro}

Some examples can be found in th annexes section with \cref{inv4,inv5}.

\begin{proof}
  If $w$ is a red-irreducible packed word, then the red-blue-factorization of
  $w$ is of the form $w = \epsilon \ins (y \insl z)$. Then the color swap correspond to
  the words $w' = \epsilon \insl (y \ins z)$ which is blue-irreducible.
\end{proof}

\subsection{Main theorem}\label{main_thm}

In \cref{sect4.1} we fixed a graded bijection $\bijirr$ between red-irreducible
and blue-irreducible packed words. After that, we extend it to all red-skeleton
forests as $\bijt$. We finished by defining $\Bij$ as a bidendriform isomorphism
from $\WQSym^*$ to $\WQSym$. Now we can set
$\bijirr: w \mapsto \widehat{w}$ as a graded bijection. The extension
$\bijt$ correspond to the color swap on red-packed forests (\textit{i.e.}
$\bijt: f \mapsto \widehat{f}$). Finally we have the following theorem:

\begin{theorem}
  The linear map $\Sigma: \WQSym^* \to \WQSym$ defined as for all packed forest
  $f$, $$\Sigma(\PP_f) \eqdef \OO_{\widehat{f}}$$ is a bidendriform isomorphism
  between $\WQSym^*$ and $\WQSym$.
\end{theorem}

\begin{proof}
  This theorem is a direct consequence of \cref{prim_tot,coro:red_to_blue}.
\end{proof}


\section*{Conclusion}

The main contribution of this paper is the combinatorial construction of biplane
trees. They are the combinatorial ingredient which completes the algebraic
theory of Foissy \cite{Foissy_2011} and allows us to describe the explicit
isomorphism. Besides, they are also an innovative combinatorial family and open
promising research perspective.

\subsection*{Generalization of the inversion of permutation to packed words}

The inherent difficulty of finding an explicit isomorphism between $\WQSym$ and
its dual lies in the fact that there is no ``inversion'' operation on packed
words. Indeed, in the case of $\FQSym$, the Hopf algebra indexed by
permutations, the isomorphism is given by the inversion of permutations. The
solution we offer, using biplane trees, is actually not a generalization of
$\FQSym$ in this sense. Indeed, even though permutations are a subset of packed
words, the restriction of our involution on packed words to permutation is not
the inversion. In particular, if a permutation $\sigma$ is both red-irreducible and
blue-irreducible, its image is itself and not its inverse. This is the case for
all packed words which are both red-irreducible and
blue-irreducible. Nevertheless, our involution is somehow ``compatible'' with
the inversion in the sense that if we arbitrary decide that the image of
$\sigma$ is $\sigma^{-1}$ for all $\sigma$ such that $\sigma$ is a red-blue-irreducible permutation,
then the rest of construction ensures that the image of $\sigma$ is
$\sigma^{-1}$ for all permutations (not necessarily irreducible anymore). But we
don't know how to define the inversion on red-blue-irreducible elements which
are not permutations, which is why to stick with the identity in all case,
including permutations.

Stays the open question: is there a generalization of the inversion of
permutations on packed words? In other words, one would want an involution on
packed words which restricts to the inversion on permutations and gives a
bidendriform isomorphism between $\WQSym$ and its dual. A consequence of our work
is that it is sufficient to find such an involution on red-blue-irreducible
packed words.

\subsection*{Generalization of the biplane trees to parking functions}

A long term goal would be to somehow generalize the structure of biplane trees
to all bidendriform Hopf algebra. The first step would be to look at the Hopf
algebra indexed by parking functions $\PQSym$. Indeed $\PQSym$ is also a
bidendriform bialgebra and parking functions are a superset of the packed
words. The question of generalizing the structure to parking functions involves
both combinatorics and algebra. The first thing is to compute bases of $\PQSym$
in which the shuffle product is not shifted. It can be done with a
generalization of \cref{BerZab_order}\cite{BerZab}.

The lines of research induced by this work are the following:
\begin{itemize}[parsep=0cm, itemsep=0cm, topsep=0cm]
\item How to generalize biplane tree structure to $\PQSym$?
\item We will then look for what are the necessary and sufficient ingredients to
  develop biplane tree structures and obtain bidendriform automorphisms on all
  bidendriform bialgebras.
\end{itemize}

\subsection*{Link between bidendriform bialgebras and skew-duplicial operad}

As said in \cref{rem:skew_dup_red,rem:skew_dup_blue} the operations $\ins$ and
$\gcdot$ (resp.$\insl$ and $\gcdot$) unexpectedly verify relations of the
skew-duplicial operad \cite{BurDel_dup}. These relations reveal a new
application of the skew-duplical operad applied on packed words.

\begin{itemize}[parsep=0cm, itemsep=0cm, topsep=0cm]
\item Can we find a skew-duplical structure on $\WQSym$ which is linked to the
  bidendriform structure?
\item More generally, is there a link between bidendriform bialgebra and
  skew-duplical?
\end{itemize}

\subsection*{Link between bidendriform bialgebras and $L$-algebras}

As said in \cref{rem:L}, the sequence that count unlabeled biplane trees is the
dimensions of the free L-algebra on one generator (see~\cite{Leroux}). It would
be interesting to investigate the link between $L$-algebra and bidendriform
bialgebras through the use of biplane trees.

\bigskip

The study of the operad on the three operations $\{\ins,\insl,\gcdot\}$ is a
start in order to study the link between bidendriform bialgebras and the
skew-duplical operade or $L$-algebras.

\paragraph*{Acknowledgments}
The computation and tests needed for this research were done using the
open-source mathematical software \textsc{SageMath} and its combinatorics
features developed by the \textsc{Sage-combinat} community. I
particularly thank F. Hivert and V. Pons for all comments on the
writing.


\bibliography{biblio}
\bibliographystyle{ieeetr}

\newpage
\section*{Annexes}\label{annexes}

In \cref{F(w)123,F(1234),F(1233),F(1223),F(1123),F(1111)} we have red-packed
forests, blue-packed forest and bicolored-packed forests associated to all
packed words of size smaller than 4.
\bigskip

In \cref{RPOQ123} we have the isomorphism between $\WQSym$ (bases $\OO$ and
$\QQ$) and its dual (bases $\PP$ and $\RR$) for size smaller than 3. Basis $\OO$
and $\PP$ are indexed by bicolored-packed forests. This illustrates the main
theorem of \cref{main_thm}.  \bigskip

In \cref{inv4} we have the involution of \cref{coro:red_to_blue} for all packed
words of size 4. They are organized by evaluations. Red-irreducible
(resp. bleu-irreducible) packed words are underlined in \redul{red} (resp
\blueul{blue}) in the first (resp. second) column.
\bigskip

In \cref{inv5} we have the involution of \cref{coro:red_to_blue} for all
red-irreducible packed words that are not blue-irreducible. It correspond to
words underlined in red in front of a word underlined in blue in \cref{inv4}.
\bigskip

The matrix of \cref{matrPR3} is redundant with the column $\RR$ and $\PP$ of
\cref{RPOQ123}. Note that even though the matrix of \cref{matrQR3} is symmetric,
it is not the case anymore on \cref{matrQR4}. Even if we restrict to
permtuations, the matrix is not symmetric for size 5.

\begin{table}[!h]
  \begin{center}
    \begin{tabular}{|l|c|c|c|c|}
      \hline
      $w$ & $\Fr(w)$ & $\Fb(w)$ & $\Frb(w)$ & $\Fbr(w)$\\
      \hline
      \hline
      $1$ & \scalebox{0.5}{\input{figures/arbres_pw/arbre1}} & \scalebox{0.5}{\input{figures/arbres_pw/arbre1left}} & \scalebox{0.5}{\input{figures/arbres_pw/arbre1ske_bi}} & \scalebox{0.5}{\input{figures/arbres_pw/arbre1leftske_bi}}\\
      \hline
    \end{tabular}
  \end{center}
  
  \begin{center}
    \begin{tabular}{|l|c|c|c|c|}
      \hline
      $w$ & $\Fr(w)$ & $\Fb(w)$ & $\Frb(w)$ & $\Fbr(w)$\\
      \hline
      \hline
      $12$ & \scalebox{0.5}{\input{figures/arbres_pw/arbre12}} & \scalebox{0.5}{\input{figures/arbres_pw/arbre12left}} & \scalebox{0.5}{\input{figures/arbres_pw/arbre12ske_bi}} & \scalebox{0.5}{\input{figures/arbres_pw/arbre12leftske_bi}}\\
      \hline
      $21$ & \scalebox{0.5}{\input{figures/arbres_pw/arbre21}} & \scalebox{0.5}{\input{figures/arbres_pw/arbre21left}} & \scalebox{0.5}{\input{figures/arbres_pw/arbre21ske_bi}} & \scalebox{0.5}{\input{figures/arbres_pw/arbre21leftske_bi}}\\
      \hline
      $11$ & \scalebox{0.5}{\input{figures/arbres_pw/arbre11}} & \scalebox{0.5}{\input{figures/arbres_pw/arbre11left}} & \scalebox{0.5}{\input{figures/arbres_pw/arbre11ske_bi}} & \scalebox{0.5}{\input{figures/arbres_pw/arbre11leftske_bi}}\\
      \hline
    \end{tabular}
  \end{center}
  
  \begin{center}
    \begin{tabular}{|l|c|c|c|c|}
      \hline
      $w$ & $\Fr(w)$ & $\Fb(w)$ & $\Frb(w)$ & $\Fbr(w)$\\
      \hline
      \hline
      $123$ & \scalebox{0.4}{\input{figures/arbres_pw/arbre123}} & \scalebox{0.4}{\input{figures/arbres_pw/arbre123left}} & \scalebox{0.4}{\input{figures/arbres_pw/arbre123ske_bi}} & \scalebox{0.4}{\input{figures/arbres_pw/arbre123leftske_bi}}\\
      \hline
      $132$ & \scalebox{0.4}{\input{figures/arbres_pw/arbre132}} & \scalebox{0.4}{\input{figures/arbres_pw/arbre132left}} & \scalebox{0.4}{\input{figures/arbres_pw/arbre132ske_bi}} & \scalebox{0.4}{\input{figures/arbres_pw/arbre132leftske_bi}}\\
      \hline
      $213$ & \scalebox{0.4}{\input{figures/arbres_pw/arbre213}} & \scalebox{0.4}{\input{figures/arbres_pw/arbre213left}} & \scalebox{0.4}{\input{figures/arbres_pw/arbre213ske_bi}} & \scalebox{0.4}{\input{figures/arbres_pw/arbre213leftske_bi}}\\
      \hline
      $231$ & \scalebox{0.4}{\input{figures/arbres_pw/arbre231}} & \scalebox{0.4}{\input{figures/arbres_pw/arbre231left}} & \scalebox{0.4}{\input{figures/arbres_pw/arbre231ske_bi}} & \scalebox{0.4}{\input{figures/arbres_pw/arbre231leftske_bi}}\\
      \hline
      $312$ & \scalebox{0.4}{\input{figures/arbres_pw/arbre312}} & \scalebox{0.4}{\input{figures/arbres_pw/arbre312left}} & \scalebox{0.4}{\input{figures/arbres_pw/arbre312ske_bi}} & \scalebox{0.4}{\input{figures/arbres_pw/arbre312leftske_bi}}\\
      \hline
      $321$ & \scalebox{0.4}{\input{figures/arbres_pw/arbre321}} & \scalebox{0.4}{\input{figures/arbres_pw/arbre321left}} & \scalebox{0.4}{\input{figures/arbres_pw/arbre321ske_bi}} & \scalebox{0.4}{\input{figures/arbres_pw/arbre321leftske_bi}}\\
      \hline
      $122$ & \scalebox{0.4}{\input{figures/arbres_pw/arbre122}} & \scalebox{0.4}{\input{figures/arbres_pw/arbre122left}} & \scalebox{0.4}{\input{figures/arbres_pw/arbre122ske_bi}} & \scalebox{0.4}{\input{figures/arbres_pw/arbre122leftske_bi}}\\
      \hline
      $212$ & \scalebox{0.4}{\input{figures/arbres_pw/arbre212}} & \scalebox{0.4}{\input{figures/arbres_pw/arbre212left}} & \scalebox{0.4}{\input{figures/arbres_pw/arbre212ske_bi}} & \scalebox{0.4}{\input{figures/arbres_pw/arbre212leftske_bi}}\\
      \hline
      $221$ & \scalebox{0.4}{\input{figures/arbres_pw/arbre221}} & \scalebox{0.4}{\input{figures/arbres_pw/arbre221left}} & \scalebox{0.4}{\input{figures/arbres_pw/arbre221ske_bi}} & \scalebox{0.4}{\input{figures/arbres_pw/arbre221leftske_bi}}\\
      \hline
      $112$ & \scalebox{0.4}{\input{figures/arbres_pw/arbre112}} & \scalebox{0.4}{\input{figures/arbres_pw/arbre112left}} & \scalebox{0.4}{\input{figures/arbres_pw/arbre112ske_bi}} & \scalebox{0.4}{\input{figures/arbres_pw/arbre112leftske_bi}}\\
      \hline
      $121$ & \scalebox{0.4}{\input{figures/arbres_pw/arbre121}} & \scalebox{0.4}{\input{figures/arbres_pw/arbre121left}} & \scalebox{0.4}{\input{figures/arbres_pw/arbre121ske_bi}} & \scalebox{0.4}{\input{figures/arbres_pw/arbre121leftske_bi}}\\
      \hline
      $211$ & \scalebox{0.4}{\input{figures/arbres_pw/arbre211}} & \scalebox{0.4}{\input{figures/arbres_pw/arbre211left}} & \scalebox{0.4}{\input{figures/arbres_pw/arbre211ske_bi}} & \scalebox{0.4}{\input{figures/arbres_pw/arbre211leftske_bi}}\\
      \hline
      $111$ & \scalebox{0.4}{\input{figures/arbres_pw/arbre111}} & \scalebox{0.4}{\input{figures/arbres_pw/arbre111left}} & \scalebox{0.4}{\input{figures/arbres_pw/arbre111ske_bi}} & \scalebox{0.4}{\input{figures/arbres_pw/arbre111leftske_bi}}\\
      \hline
    \end{tabular}
  \end{center}
  \caption{All packed words of size smaller than 3 and forests associated to it.}
  \label{F(w)123}
\end{table}

\begin{table}[!h]
  \begin{center}
    \begin{tabular}{|l|c|c|c|c|}
      \hline
      $w$ & $\Tr(w)$ & $\Tb(w)$ & $\Trb(w)$ & $\Tbr(w)$\\
      \hline
      \hline
      $1234$ & \scalebox{0.4}{\input{figures/arbres_pw/arbre1234}} & \scalebox{0.4}{\input{figures/arbres_pw/arbre1234left}} & \scalebox{0.4}{\input{figures/arbres_pw/arbre1234ske_bi}} & \scalebox{0.4}{\input{figures/arbres_pw/arbre1234leftske_bi}}\\
      \hline
      $1243$ & \scalebox{0.4}{\input{figures/arbres_pw/arbre1243}} & \scalebox{0.4}{\input{figures/arbres_pw/arbre1243left}} & \scalebox{0.4}{\input{figures/arbres_pw/arbre1243ske_bi}} & \scalebox{0.4}{\input{figures/arbres_pw/arbre1243leftske_bi}}\\
      \hline
      $1324$ & \scalebox{0.4}{\input{figures/arbres_pw/arbre1324}} & \scalebox{0.4}{\input{figures/arbres_pw/arbre1324left}} & \scalebox{0.4}{\input{figures/arbres_pw/arbre1324ske_bi}} & \scalebox{0.4}{\input{figures/arbres_pw/arbre1324leftske_bi}}\\
      \hline
      $1342$ & \scalebox{0.4}{\input{figures/arbres_pw/arbre1342}} & \scalebox{0.4}{\input{figures/arbres_pw/arbre1342left}} & \scalebox{0.4}{\input{figures/arbres_pw/arbre1342ske_bi}} & \scalebox{0.4}{\input{figures/arbres_pw/arbre1342leftske_bi}}\\
      \hline
      $1423$ & \scalebox{0.4}{\input{figures/arbres_pw/arbre1423}} & \scalebox{0.4}{\input{figures/arbres_pw/arbre1423left}} & \scalebox{0.4}{\input{figures/arbres_pw/arbre1423ske_bi}} & \scalebox{0.4}{\input{figures/arbres_pw/arbre1423leftske_bi}}\\
      \hline
      $1432$ & \scalebox{0.4}{\input{figures/arbres_pw/arbre1432}} & \scalebox{0.4}{\input{figures/arbres_pw/arbre1432left}} & \scalebox{0.4}{\input{figures/arbres_pw/arbre1432ske_bi}} & \scalebox{0.4}{\input{figures/arbres_pw/arbre1432leftske_bi}}\\
      \hline
      $2134$ & \scalebox{0.4}{\input{figures/arbres_pw/arbre2134}} & \scalebox{0.4}{\input{figures/arbres_pw/arbre2134left}} & \scalebox{0.4}{\input{figures/arbres_pw/arbre2134ske_bi}} & \scalebox{0.4}{\input{figures/arbres_pw/arbre2134leftske_bi}}\\
      \hline
      $2143$ & \scalebox{0.4}{\input{figures/arbres_pw/arbre2143}} & \scalebox{0.4}{\input{figures/arbres_pw/arbre2143left}} & \scalebox{0.4}{\input{figures/arbres_pw/arbre2143ske_bi}} & \scalebox{0.4}{\input{figures/arbres_pw/arbre2143leftske_bi}}\\
      \hline
      $2314$ & \scalebox{0.4}{\input{figures/arbres_pw/arbre2314}} & \scalebox{0.4}{\input{figures/arbres_pw/arbre2314left}} & \scalebox{0.4}{\input{figures/arbres_pw/arbre2314ske_bi}} & \scalebox{0.4}{\input{figures/arbres_pw/arbre2314leftske_bi}}\\
      \hline
      $2413$ & \scalebox{0.4}{\input{figures/arbres_pw/arbre2413}} & \scalebox{0.4}{\input{figures/arbres_pw/arbre2413left}} & \scalebox{0.4}{\input{figures/arbres_pw/arbre2413ske_bi}} & \scalebox{0.4}{\input{figures/arbres_pw/arbre2413leftske_bi}}\\
      \hline
      $3124$ & \scalebox{0.4}{\input{figures/arbres_pw/arbre3124}} & \scalebox{0.4}{\input{figures/arbres_pw/arbre3124left}} & \scalebox{0.4}{\input{figures/arbres_pw/arbre3124ske_bi}} & \scalebox{0.4}{\input{figures/arbres_pw/arbre3124leftske_bi}}\\
      \hline
      $3142$ & \scalebox{0.4}{\input{figures/arbres_pw/arbre3142}} & \scalebox{0.4}{\input{figures/arbres_pw/arbre3142left}} & \scalebox{0.4}{\input{figures/arbres_pw/arbre3142ske_bi}} & \scalebox{0.4}{\input{figures/arbres_pw/arbre3142leftske_bi}}\\
      \hline
      $3214$ & \scalebox{0.4}{\input{figures/arbres_pw/arbre3214}} & \scalebox{0.4}{\input{figures/arbres_pw/arbre3214left}} & \scalebox{0.4}{\input{figures/arbres_pw/arbre3214ske_bi}} & \scalebox{0.4}{\input{figures/arbres_pw/arbre3214leftske_bi}}\\
      \hline
    \end{tabular}
  \end{center}
  \caption{Packed words of size 4 and associated forests (part 1).}
  \label{F(1234)}
\end{table}

\begin{table}[!h]
  \begin{center}
    \begin{tabular}{|l|c|c|c|c|}
      \hline
      $w$ & $\Tr(w)$ & $\Tb(w)$ & $\Trb(w)$ & $\Tbr(w)$\\
      \hline
      \hline
      $1233$ & \scalebox{0.4}{\input{figures/arbres_pw/arbre1233}} & \scalebox{0.4}{\input{figures/arbres_pw/arbre1233left}} & \scalebox{0.4}{\input{figures/arbres_pw/arbre1233ske_bi}} & \scalebox{0.4}{\input{figures/arbres_pw/arbre1233leftske_bi}}\\
      \hline
      $1323$ & \scalebox{0.4}{\input{figures/arbres_pw/arbre1323}} & \scalebox{0.4}{\input{figures/arbres_pw/arbre1323left}} & \scalebox{0.4}{\input{figures/arbres_pw/arbre1323ske_bi}} & \scalebox{0.4}{\input{figures/arbres_pw/arbre1323leftske_bi}}\\
      \hline
      $1332$ & \scalebox{0.4}{\input{figures/arbres_pw/arbre1332}} & \scalebox{0.4}{\input{figures/arbres_pw/arbre1332left}} & \scalebox{0.4}{\input{figures/arbres_pw/arbre1332ske_bi}} & \scalebox{0.4}{\input{figures/arbres_pw/arbre1332leftske_bi}}\\
      \hline
      $2133$ & \scalebox{0.4}{\input{figures/arbres_pw/arbre2133}} & \scalebox{0.4}{\input{figures/arbres_pw/arbre2133left}} & \scalebox{0.4}{\input{figures/arbres_pw/arbre2133ske_bi}} & \scalebox{0.4}{\input{figures/arbres_pw/arbre2133leftske_bi}}\\
      \hline
      $2313$ & \scalebox{0.4}{\input{figures/arbres_pw/arbre2313}} & \scalebox{0.4}{\input{figures/arbres_pw/arbre2313left}} & \scalebox{0.4}{\input{figures/arbres_pw/arbre2313ske_bi}} & \scalebox{0.4}{\input{figures/arbres_pw/arbre2313leftske_bi}}\\
      \hline
      $3123$ & \scalebox{0.4}{\input{figures/arbres_pw/arbre3123}} & \scalebox{0.4}{\input{figures/arbres_pw/arbre3123left}} & \scalebox{0.4}{\input{figures/arbres_pw/arbre3123ske_bi}} & \scalebox{0.4}{\input{figures/arbres_pw/arbre3123leftske_bi}}\\
      \hline
      $3132$ & \scalebox{0.4}{\input{figures/arbres_pw/arbre3132}} & \scalebox{0.4}{\input{figures/arbres_pw/arbre3132left}} & \scalebox{0.4}{\input{figures/arbres_pw/arbre3132ske_bi}} & \scalebox{0.4}{\input{figures/arbres_pw/arbre3132leftske_bi}}\\
      \hline
      $3213$ & \scalebox{0.4}{\input{figures/arbres_pw/arbre3213}} & \scalebox{0.4}{\input{figures/arbres_pw/arbre3213left}} & \scalebox{0.4}{\input{figures/arbres_pw/arbre3213ske_bi}} & \scalebox{0.4}{\input{figures/arbres_pw/arbre3213leftske_bi}}\\
      \hline
    \end{tabular}
  \end{center}
  \caption{Packed words of size 4 and associated forests (part 2).}
  \label{F(1233)}
\end{table}

\begin{table}[!h]
  \begin{center}
    \begin{tabular}{|l|c|c|c|c|}
      \hline
      $w$ & $\Tr(w)$ & $\Tb(w)$ & $\Trb(w)$ & $\Tbr(w)$\\
      \hline
      \hline
      $1223$ & \scalebox{0.4}{\input{figures/arbres_pw/arbre1223}} & \scalebox{0.4}{\input{figures/arbres_pw/arbre1223left}} & \scalebox{0.4}{\input{figures/arbres_pw/arbre1223ske_bi}} & \scalebox{0.4}{\input{figures/arbres_pw/arbre1223leftske_bi}}\\
      \hline
      $1232$ & \scalebox{0.4}{\input{figures/arbres_pw/arbre1232}} & \scalebox{0.4}{\input{figures/arbres_pw/arbre1232left}} & \scalebox{0.4}{\input{figures/arbres_pw/arbre1232ske_bi}} & \scalebox{0.4}{\input{figures/arbres_pw/arbre1232leftske_bi}}\\
      \hline
      $1322$ & \scalebox{0.4}{\input{figures/arbres_pw/arbre1322}} & \scalebox{0.4}{\input{figures/arbres_pw/arbre1322left}} & \scalebox{0.4}{\input{figures/arbres_pw/arbre1322ske_bi}} & \scalebox{0.4}{\input{figures/arbres_pw/arbre1322leftske_bi}}\\
      \hline
      $2123$ & \scalebox{0.4}{\input{figures/arbres_pw/arbre2123}} & \scalebox{0.4}{\input{figures/arbres_pw/arbre2123left}} & \scalebox{0.4}{\input{figures/arbres_pw/arbre2123ske_bi}} & \scalebox{0.4}{\input{figures/arbres_pw/arbre2123leftske_bi}}\\
      \hline
      $2132$ & \scalebox{0.4}{\input{figures/arbres_pw/arbre2132}} & \scalebox{0.4}{\input{figures/arbres_pw/arbre2132left}} & \scalebox{0.4}{\input{figures/arbres_pw/arbre2132ske_bi}} & \scalebox{0.4}{\input{figures/arbres_pw/arbre2132leftske_bi}}\\
      \hline
      $2213$ & \scalebox{0.4}{\input{figures/arbres_pw/arbre2213}} & \scalebox{0.4}{\input{figures/arbres_pw/arbre2213left}} & \scalebox{0.4}{\input{figures/arbres_pw/arbre2213ske_bi}} & \scalebox{0.4}{\input{figures/arbres_pw/arbre2213leftske_bi}}\\
      \hline
      $2312$ & \scalebox{0.4}{\input{figures/arbres_pw/arbre2312}} & \scalebox{0.4}{\input{figures/arbres_pw/arbre2312left}} & \scalebox{0.4}{\input{figures/arbres_pw/arbre2312ske_bi}} & \scalebox{0.4}{\input{figures/arbres_pw/arbre2312leftske_bi}}\\
      \hline
    \end{tabular}
  \end{center}
  \caption{Packed words of size 4 and associated forests (part 3).}
  \label{F(1223)}
\end{table}
\begin{table}[!h]
  \begin{center}
    \begin{tabular}{|l|c|c|c|c|}
      \hline
      $w$ & $\Tr(w)$ & $\Tb(w)$ & $\Trb(w)$ & $\Tbr(w)$\\
      \hline
      \hline
      $1123$ & \scalebox{0.4}{\input{figures/arbres_pw/arbre1123}} & \scalebox{0.4}{\input{figures/arbres_pw/arbre1123left}} & \scalebox{0.4}{\input{figures/arbres_pw/arbre1123ske_bi}} & \scalebox{0.4}{\input{figures/arbres_pw/arbre1123leftske_bi}}\\
      \hline
      $1132$ & \scalebox{0.4}{\input{figures/arbres_pw/arbre1132}} & \scalebox{0.4}{\input{figures/arbres_pw/arbre1132left}} & \scalebox{0.4}{\input{figures/arbres_pw/arbre1132ske_bi}} & \scalebox{0.4}{\input{figures/arbres_pw/arbre1132leftske_bi}}\\
      \hline
      $1213$ & \scalebox{0.4}{\input{figures/arbres_pw/arbre1213}} & \scalebox{0.4}{\input{figures/arbres_pw/arbre1213left}} & \scalebox{0.4}{\input{figures/arbres_pw/arbre1213ske_bi}} & \scalebox{0.4}{\input{figures/arbres_pw/arbre1213leftske_bi}}\\
      \hline
      $1231$ & \scalebox{0.4}{\input{figures/arbres_pw/arbre1231}} & \scalebox{0.4}{\input{figures/arbres_pw/arbre1231left}} & \scalebox{0.4}{\input{figures/arbres_pw/arbre1231ske_bi}} & \scalebox{0.4}{\input{figures/arbres_pw/arbre1231leftske_bi}}\\
      \hline
      $1312$ & \scalebox{0.4}{\input{figures/arbres_pw/arbre1312}} & \scalebox{0.4}{\input{figures/arbres_pw/arbre1312left}} & \scalebox{0.4}{\input{figures/arbres_pw/arbre1312ske_bi}} & \scalebox{0.4}{\input{figures/arbres_pw/arbre1312leftske_bi}}\\
      \hline
      $1321$ & \scalebox{0.4}{\input{figures/arbres_pw/arbre1321}} & \scalebox{0.4}{\input{figures/arbres_pw/arbre1321left}} & \scalebox{0.4}{\input{figures/arbres_pw/arbre1321ske_bi}} & \scalebox{0.4}{\input{figures/arbres_pw/arbre1321leftske_bi}}\\
      \hline
      $2113$ & \scalebox{0.4}{\input{figures/arbres_pw/arbre2113}} & \scalebox{0.4}{\input{figures/arbres_pw/arbre2113left}} & \scalebox{0.4}{\input{figures/arbres_pw/arbre2113ske_bi}} & \scalebox{0.4}{\input{figures/arbres_pw/arbre2113leftske_bi}}\\
      \hline
      $2131$ & \scalebox{0.4}{\input{figures/arbres_pw/arbre2131}} & \scalebox{0.4}{\input{figures/arbres_pw/arbre2131left}} & \scalebox{0.4}{\input{figures/arbres_pw/arbre2131ske_bi}} & \scalebox{0.4}{\input{figures/arbres_pw/arbre2131leftske_bi}}\\
      \hline
    \end{tabular}
  \end{center}
  \caption{Packed words of size 4 and associated forests (part 4).}
  \label{F(1123)}
\end{table}

\begin{table}[!h]
  \begin{center}
    \begin{tabular}{|l|c|c|c|c|}
      \hline
      $w$ & $\Tr(w)$ & $\Tb(w)$ & $\Trb(w)$ & $\Tbr(w)$\\
      \hline
      \hline
      $1222$ & \scalebox{0.4}{\input{figures/arbres_pw/arbre1222}} & \scalebox{0.4}{\input{figures/arbres_pw/arbre1222left}} & \scalebox{0.4}{\input{figures/arbres_pw/arbre1222ske_bi}} & \scalebox{0.4}{\input{figures/arbres_pw/arbre1222leftske_bi}}\\
      \hline
      $2122$ & \scalebox{0.4}{\input{figures/arbres_pw/arbre2122}} & \scalebox{0.4}{\input{figures/arbres_pw/arbre2122left}} & \scalebox{0.4}{\input{figures/arbres_pw/arbre2122ske_bi}} & \scalebox{0.4}{\input{figures/arbres_pw/arbre2122leftske_bi}}\\
      \hline
      $2212$ & \scalebox{0.4}{\input{figures/arbres_pw/arbre2212}} & \scalebox{0.4}{\input{figures/arbres_pw/arbre2212left}} & \scalebox{0.4}{\input{figures/arbres_pw/arbre2212ske_bi}} & \scalebox{0.4}{\input{figures/arbres_pw/arbre2212leftske_bi}}\\
      \hline

      $1122$ & \scalebox{0.4}{\input{figures/arbres_pw/arbre1122}} & \scalebox{0.4}{\input{figures/arbres_pw/arbre1122left}} & \scalebox{0.4}{\input{figures/arbres_pw/arbre1122ske_bi}} & \scalebox{0.4}{\input{figures/arbres_pw/arbre1122leftske_bi}}\\
      \hline
      $1212$ & \scalebox{0.4}{\input{figures/arbres_pw/arbre1212}} & \scalebox{0.4}{\input{figures/arbres_pw/arbre1212left}} & \scalebox{0.4}{\input{figures/arbres_pw/arbre1212ske_bi}} & \scalebox{0.4}{\input{figures/arbres_pw/arbre1212leftske_bi}}\\
      \hline
      $1221$ & \scalebox{0.4}{\input{figures/arbres_pw/arbre1221}} & \scalebox{0.4}{\input{figures/arbres_pw/arbre1221left}} & \scalebox{0.4}{\input{figures/arbres_pw/arbre1221ske_bi}} & \scalebox{0.4}{\input{figures/arbres_pw/arbre1221leftske_bi}}\\
      \hline
      $2112$ & \scalebox{0.4}{\input{figures/arbres_pw/arbre2112}} & \scalebox{0.4}{\input{figures/arbres_pw/arbre2112left}} & \scalebox{0.4}{\input{figures/arbres_pw/arbre2112ske_bi}} & \scalebox{0.4}{\input{figures/arbres_pw/arbre2112leftske_bi}}\\
      \hline
      $2121$ & \scalebox{0.4}{\input{figures/arbres_pw/arbre2121}} & \scalebox{0.4}{\input{figures/arbres_pw/arbre2121left}} & \scalebox{0.4}{\input{figures/arbres_pw/arbre2121ske_bi}} & \scalebox{0.4}{\input{figures/arbres_pw/arbre2121leftske_bi}}\\
      \hline
      $1112$ & \scalebox{0.4}{\input{figures/arbres_pw/arbre1112}} & \scalebox{0.4}{\input{figures/arbres_pw/arbre1112left}} & \scalebox{0.4}{\input{figures/arbres_pw/arbre1112ske_bi}} & \scalebox{0.4}{\input{figures/arbres_pw/arbre1112leftske_bi}}\\
      \hline
      $1121$ & \scalebox{0.4}{\input{figures/arbres_pw/arbre1121}} & \scalebox{0.4}{\input{figures/arbres_pw/arbre1121left}} & \scalebox{0.4}{\input{figures/arbres_pw/arbre1121ske_bi}} & \scalebox{0.4}{\input{figures/arbres_pw/arbre1121leftske_bi}}\\
      \hline
      $1211$ & \scalebox{0.4}{\input{figures/arbres_pw/arbre1211}} & \scalebox{0.4}{\input{figures/arbres_pw/arbre1211left}} & \scalebox{0.4}{\input{figures/arbres_pw/arbre1211ske_bi}} & \scalebox{0.4}{\input{figures/arbres_pw/arbre1211leftske_bi}}\\
      \hline
      $1111$ & \scalebox{0.4}{\input{figures/arbres_pw/arbre1111}} & \scalebox{0.4}{\input{figures/arbres_pw/arbre1111left}} & \scalebox{0.4}{\input{figures/arbres_pw/arbre1111ske_bi}} & \scalebox{0.4}{\input{figures/arbres_pw/arbre1111leftske_bi}}\\
      \hline
    \end{tabular}
  \end{center}
  \caption{Packed words of size 4 and associated forests (part 5).}
  \label{F(1111)}
\end{table}

\begin{table}[!h]
  \begin{center}
    \begin{tabular}{|l|c||c|l|}
      \hline
      $\RR_1$ & $\PP_{\scalebox{0.5}{\input{figures/arbres_pw/arbre1ske_bi}}}$ & $\OO_{\scalebox{0.5}{\input{figures/arbres_pw/arbre1leftske_bi}}}$ & $\QQ_1$\\
      \hline
    \end{tabular}
  \end{center}

  \begin{center}
    \begin{tabular}{|l|c||c|l|}
      \hline
      $\RR_{12} - \RR_{21}$ & $\PP_{\scalebox{0.5}{\input{figures/arbres_pw/arbre12ske_bi}}}$ & $\OO_{\scalebox{0.5}{\input{figures/arbres_pw/arbre12leftske_bi}}}$ & $\QQ_{12} - \QQ_{21}$\\
      \hline
      $\RR_{21}$ & $\PP_{\scalebox{0.5}{\input{figures/arbres_pw/arbre21ske_bi}}}$ & $\OO_{\scalebox{0.5}{\input{figures/arbres_pw/arbre21leftske_bi}}}$ & $\QQ_{21}$\\
      \hline
      $\RR_{11}$ & $\PP_{\scalebox{0.5}{\input{figures/arbres_pw/arbre11ske_bi}}}$ & $\OO_{\scalebox{0.5}{\input{figures/arbres_pw/arbre11leftske_bi}}}$ & $\QQ_{11}$\\
      \hline
    \end{tabular}
  \end{center}

  \begin{center}
    \begin{tabular}{|l|c||c|l|}
      \hline
      $\RR$ & $\PP$ & $\OO$ & $\QQ$\\
      \hline
      \hline
      $ 123 - 213 - 231 + 321 $ & $ \scalebox{0.5}{\input{figures/arbres_pw/arbre123ske_bi}} $ & $ \scalebox{0.5}{\input{figures/arbres_pw/arbre123leftske_bi}} $ & $ 123 - 213 - 312 + 321 $\\
      \hline
      $ 132 - 231 $ & $ \scalebox{0.5}{\input{figures/arbres_pw/arbre132ske_bi}} $ & $ \scalebox{0.5}{\input{figures/arbres_pw/arbre132leftske_bi}} $ & $ 132 - 312 $\\
      \hline
      $ 213 - 312 + 231 - 132 $ & $ \scalebox{0.5}{\input{figures/arbres_pw/arbre213ske_bi}} $ & $ \scalebox{0.5}{\input{figures/arbres_pw/arbre213leftske_bi}} $ & $ 213 + 312 - 231 - 132 $\\
      \hline
      $ 231 - 321 $ & $ \scalebox{0.5}{\input{figures/arbres_pw/arbre231ske_bi}} $ & $ \scalebox{0.5}{\input{figures/arbres_pw/arbre312leftske_bi}} $ & $ 312 - 321 $\\
      \hline
      $ 132 + 312 - 231 - 321 $ & $ \scalebox{0.5}{\input{figures/arbres_pw/arbre312ske_bi}} $ & $ \scalebox{0.5}{\input{figures/arbres_pw/arbre231leftske_bi}} $ & $ 231 + 132 - 321 - 312 $\\
      \hline
      $ 321 $ & $ \scalebox{0.5}{\input{figures/arbres_pw/arbre321ske_bi}} $ & $ \scalebox{0.5}{\input{figures/arbres_pw/arbre321leftske_bi}} $ & $ 321 $\\
      \hline
      $ 122 - 121 + 212 - 211 $ & $ \scalebox{0.5}{\input{figures/arbres_pw/arbre122ske_bi}} $ & $ \scalebox{0.5}{\input{figures/arbres_pw/arbre212leftske_bi}} $ & $ 122 - 221 $\\
      \hline
      $ 212 $ & $ \scalebox{0.5}{\input{figures/arbres_pw/arbre212ske_bi}} $ & $ \scalebox{0.5}{\input{figures/arbres_pw/arbre122leftske_bi}} $ & $ 122 - 212 $\\
      \hline
      $ 221 $ & $ \scalebox{0.5}{\input{figures/arbres_pw/arbre221ske_bi}} $ & $ \scalebox{0.5}{\input{figures/arbres_pw/arbre211leftske_bi}} $ & $ 211 $\\
      \hline
      $ 112 - 221 $ & $ \scalebox{0.5}{\input{figures/arbres_pw/arbre112ske_bi}} $ & $ \scalebox{0.5}{\input{figures/arbres_pw/arbre112leftske_bi}} $ & $ 112 - 211 $\\
      \hline
      $ 121 $ & $ \scalebox{0.5}{\input{figures/arbres_pw/arbre121ske_bi}} $ & $ \scalebox{0.5}{\input{figures/arbres_pw/arbre121leftske_bi}} $ & $ 121 - 211 $\\
      \hline
      $ 121 + 211 $ & $ \scalebox{0.5}{\input{figures/arbres_pw/arbre211ske_bi}} $ & $ \scalebox{0.5}{\input{figures/arbres_pw/arbre221leftske_bi}} $ & $ 221 $\\
      \hline
      $ 111 $ & $ \scalebox{0.5}{\input{figures/arbres_pw/arbre111ske_bi}} $ & $ \scalebox{0.5}{\input{figures/arbres_pw/arbre111leftske_bi}} $ & $ 111 $\\
      \hline
    \end{tabular}
  \end{center}
  \caption{The automorphism of $\WQSym_{\leq3}$.}
  \label{RPOQ123}
\end{table}

\begin{table}[h]
  \begin{center}
    \begin{tabular}{ccccc}
      \multirow{24}{*}{\begin{tabular}{|c|c|}
  \hline
  $1234$ & $1234$ \\
  \hline
  $\redul{\blueul{1243}}$ & $\redul{\blueul{1243}}$ \\
  \hline
  $\redul{\blueul{1324}}$ & $\redul{\blueul{1324}}$ \\
  \hline
  $1342$ & $1342$ \\
  \hline
  $\redul{\blueul{1423}}$ & $\redul{\blueul{1423}}$ \\
  \hline
  $\redul{\blueul{1432}}$ & $\redul{\blueul{1432}}$ \\
  \hline
  $2134$ & $2134$ \\
  \hline
  $\redul{\blueul{2143}}$ & $\redul{\blueul{2143}}$ \\
  \hline
  $2314$ & $3124$ \\
  \hline
  $2341$ & $4123$ \\
  \hline
  $\redul{\blueul{2413}}$ & $\redul{\blueul{3142}}$ \\
  \hline
  $2431$ & $4132$ \\
  \hline
  $3124$ & $2314$ \\
  \hline
  $3142$ & $2413$ \\
  \hline
  $3214$ & $3214$ \\
  \hline
  $3241$ & $4213$ \\
  \hline
  $3412$ & $3412$ \\
  \hline
  $3421$ & $4312$ \\
  \hline
  $4123$ & $2341$ \\
  \hline
  $4132$ & $2431$ \\
  \hline
  $4213$ & $3241$ \\
  \hline
  $4231$ & $4231$ \\
  \hline
  $4312$ & $3421$ \\
  \hline
  $4321$ & $4321$ \\
  \hline
\end{tabular}} & \multirow{12}{*}{\begin{tabular}{|c|c|}
  \hline
  $\blueul{1233}$ & $\redul{3123}$ \\
  \hline
  $\redul{\blueul{1323}}$ & $\redul{\blueul{1323}}$ \\
  \hline
  $\redul{\blueul{1332}}$ & $\redul{\blueul{1332}}$ \\
  \hline
  $2133$ & $3213$ \\
  \hline
  $2313$ & $2313$ \\
  \hline
  $2331$ & $3212$ \\
  \hline
  $\redul{3123}$ & $\blueul{1233}$ \\
  \hline
  $\redul{\blueul{3132}}$ & $\redul{\blueul{3132}}$ \\
  \hline
  $3213$ & $2133$ \\
  \hline
  $3231$ & $3122$ \\
  \hline
  $3312$ & $2311$ \\
  \hline
  $3321$ & $3211$ \\
  \hline
\end{tabular}} & \multirow{12}{*}{\begin{tabular}{|c|c|}
  \hline
  $1223$ & $2123$ \\
  \hline
  $\redul{\blueul{1232}}$ & $\redul{\blueul{1232}}$ \\
  \hline
  $\redul{\blueul{1322}}$ & $\redul{\blueul{1322}}$ \\
  \hline
  $2123$ & $1223$ \\
  \hline
  $\redul{\blueul{2132}}$ & $\redul{\blueul{2132}}$ \\
  \hline
  $2213$ & $2113$ \\
  \hline
  $2231$ & $3112$ \\
  \hline
  $\redul{2312}$ & $\blueul{2131}$ \\
  \hline
  $2321$ & $3121$ \\
  \hline
  $3122$ & $3231$ \\
  \hline
  $3212$ & $2331$ \\
  \hline
  $3221$ & $3221$ \\
  \hline
\end{tabular}} & \multirow{12}{*}{\begin{tabular}{|c|c|}
  \hline
  $1123$ & $1123$ \\
  \hline
  $\redul{\blueul{1132}}$ & $\redul{\blueul{1132}}$ \\
  \hline
  $1213$ & $1213$ \\
  \hline
  $\redul{\blueul{1231}}$ & $\redul{\blueul{1231}}$ \\
  \hline
  $\redul{\blueul{1312}}$ & $\redul{\blueul{1312}}$ \\
  \hline
  $\redul{\blueul{1321}}$ & $\redul{\blueul{1321}}$ \\
  \hline
  $2113$ & $2213$ \\
  \hline
  $\blueul{2131}$ & $\redul{2312}$ \\
  \hline
  $2311$ & $3312$ \\
  \hline
  $3112$ & $2231$ \\
  \hline
  $3121$ & $2321$ \\
  \hline
  $3211$ & $3321$ \\
  \hline
\end{tabular}} \\
                                              & & & \\
                                              & & & \\
                                              & & & \\
                                              & & & \\
                                              & & & \\
                                              & & & \\
                                              & & & \\
                                              & & & \\
                                              & & & \\
                                              & & & \\
                                              & & & \\
                                              & & & \\
                                              & & & \\
                                              & & & \\
                                              & & & \\
      \vspace{0.25cm}
                                              & & & \\
                                              & & & \\
                                              & \multirow{6}{*}{\begin{tabular}{|c|c|}
  \hline
  $\blueul{1122}$ & $\redul{2112}$ \\
  \hline
  $\redul{\blueul{1212}}$ & $\redul{\blueul{1212}}$ \\
  \hline
  $\redul{\blueul{1221}}$ & $\redul{\blueul{1221}}$ \\
  \hline
  $\redul{2112}$ & $\blueul{1122}$ \\
  \hline
  $\redul{\blueul{2121}}$ & $\redul{\blueul{2121}}$ \\
  \hline
  $2211$ & $2211$ \\
  \hline
\end{tabular}
} & \multirow{4}{*}{\begin{tabular}{|c|c|}
  \hline
  $\blueul{1222}$ & $\redul{2212}$ \\
  \hline
  $\redul{\blueul{2122}}$ & $\redul{\blueul{2122}}$ \\
  \hline
  $\redul{2212}$ & $\blueul{1222}$ \\
  \hline
  $2221$ & $2111$ \\
  \hline
\end{tabular}} & \multirow{4}{*}{
  \begin{tabular}{|c|c|}
    \hline
    $1112$ & $1112$ \\
    \hline
    $\redul{\blueul{1121}}$ & $\redul{\blueul{1121}}$ \\
    \hline
    $\redul{\blueul{1211}}$ & $\redul{\blueul{1211}}$ \\
    \hline
    $2111$ & $2221$ \\
    \hline
  \end{tabular}} \\
                                              & & & \\
                                              & & & \\
      \vspace{0.07cm}
                                              & & & \\
                                              & & & \\
                                              & & \begin{tabular}{|c|c|}
  \hline
  $\redul{\blueul{1111}}$ & $\redul{\blueul{1111}}$ \\
  \hline
\end{tabular} & \\
    \end{tabular}
    \caption{The involution $w \mapsto \widehat{w}$ on packed words of size 4.}
    \label{inv4}
  \end{center}
\end{table}

\begin{table}[h]
  \begin{center}
    \begin{tabular}{cccc}
      \multirow{7}{*}{\begin{tabular}{|c|c|}
  \hline
  $23514$ & $41253$ \\
  \hline
  $24513$ & $41352$ \\
  \hline
  $25314$ & $41523$ \\
  \hline
  $25413$ & $41532$ \\
  \hline
  $32514$ & $42153$ \\
  \hline
  $35124$ & $34152$ \\
  \hline
  $35214$ & $43152$ \\
  \hline
\end{tabular}} & \multirow{9}{*}{\begin{tabular}{|c|c|}
  \hline
  $24314$ & $31424$ \\
  \hline
  $24413$ & $31442$ \\
  \hline
  $41234$ & $12344$ \\
  \hline
  $41324$ & $13244$ \\
  \hline
  $42134$ & $21344$ \\
  \hline
  $42314$ & $31244$ \\
  \hline
  $42413$ & $34142$ \\
  \hline
  $43124$ & $23144$ \\
  \hline
  $43214$ & $32144$ \\
  \hline
\end{tabular}} & \multirow{8}{*}{\begin{tabular}{|c|c|}
  \hline
  $23413$ & $31242$ \\
  \hline
  $24313$ & $31422$ \\
  \hline
  $32413$ & $32142$ \\
  \hline
  $34123$ & $23141$ \\
  \hline
  $34213$ & $32141$ \\
  \hline
  $32313$ & $23133$ \\
  \hline
  $33123$ & $12333$ \\
  \hline
  $33213$ & $21333$ \\
  \hline
\end{tabular}} & \multirow{4}{*}{\begin{tabular}{|c|c|}
  \hline
  $22413$ & $31142$ \\
  \hline
  $23412$ & $31241$ \\
  \hline
  $24213$ & $31412$ \\
  \hline
  $24312$ & $31421$ \\
  \hline
\end{tabular}} \\
                                              & & & \\
                                              & & & \\
                                              & & & \\
                                              & & & \\
                                              & & & \\
                                              & & & \\
                                              & & & \\
                                              & & & \\
                                              & & & \\
      \multirow{6}{*}{\begin{tabular}{|c|c|}
  \hline
  $23213$ & $21313$ \\
  \hline
  $23312$ & $21331$ \\
  \hline
  $31223$ & $21233$ \\
  \hline
  $32123$ & $12233$ \\
  \hline
  $32213$ & $21133$ \\
  \hline
  $32312$ & $23131$ \\
  \hline
\end{tabular}} & \multirow{2}{*}{\begin{tabular}{|c|c|}
  \hline
  $22312$ & $21131$ \\
  \hline
  $23212$ & $21311$ \\
  \hline
\end{tabular}} & \multirow{1}{*}{\begin{tabular}{|c|c|}
  \hline
  $22212$ & $12222$ \\
  \hline
\end{tabular}} & \multirow{1}{*}{\begin{tabular}{|c|c|}
  \hline
  $24113$ & $33142$ \\
  \hline
\end{tabular}}\\
                                              & & & \\
                                              & & & \\
                                              & \multirow{3}{*}{\begin{tabular}{|c|c|}
  \hline
  $31123$ & $11233$ \\
  \hline
  $31213$ & $12133$ \\
  \hline
  $32113$ & $22133$ \\
  \hline
\end{tabular}} & \multirow{1}{*}{\begin{tabular}{|c|c|}
  \hline
  $23112$ & $22131$ \\
  \hline
\end{tabular}} & \multirow{1}{*}{\begin{tabular}{|c|c|}
  \hline
  $22112$ & $11222$ \\
  \hline
\end{tabular}}\\
                                              & & & \\
                                              & & \multirow{1}{*}{\begin{tabular}{|c|c|}
  \hline
  $21112$ & $11122$ \\
  \hline
\end{tabular}} & \\
    \end{tabular}
    \caption{The involution $w \mapsto \widehat{w}$ on red-irreducible packed words
      that are not blue-irreducible of size 5.}
    \label{inv5}
  \end{center}
\end{table}

\begin{figure}[!h]
  \centering
  {
    \begin{displaymath}
      \begin{array}{c|p{0.5cm}@{\,}p{0.5cm}@{\,}p{0.5cm}@{\,}p{0.5cm}@{\,}p{0.5cm}@{\,}p{0.5cm}@{\,}p{0.5cm}@{\,}p{0.5cm}@{\,}p{0.5cm}@{\,}p{0.5cm}@{\,}p{0.5cm}@{\,}p{0.5cm}@{\,}p{0.5cm}@{\,}p{0.5cm}@{\,}p{0.5cm}@{\,}p{0.5cm}@{\,}p{0.5cm}@{}|c}
        &\pc{123}&\pc{132}&\pc{213}&\pc{231}&\pc{312}&\pc{321}&\pc{122}&\pc{212}&\pc{221}&\pc{112}&\pc{121}&\pc{211}&\pc{111}\\
        \hline
        123& 1& .& .& .& .& .& .& .& .& .& .& .& .\\
        132& .& 1&-1& .& 1& .& .& .& .& .& .& .& .\\
        213&-1& .& 1& .& .& .& .& .& .& .& .& .& .\\
        231&-1&-1& 1& 1&-1& .& .& .& .& .& .& .& .\\
        312& .& .&-1& .& 1& .& .& .& .& .& .& .& .\\
        321& 1& .& .&-1&-1& 1& .& .& .& .& .& .& .\\
        122& .& .& .& .& .& .& 1& .& .& .& .& .& .\\
        212& .& .& .& .& .& .& 1& 1& .& .& .& .& .\\
        221& .& .& .& .& .& .& .& .& 1&-1& .& .& .\\
        112& .& .& .& .& .& .& .& .& .& 1& .& .& .\\
        121& .& .& .& .& .& .&-1& .& .& .& 1& 1& .\\
        211& .& .& .& .& .& .&-1& .& .& .& .& 1& .\\
        111& .& .& .& .& .& .& .& .& .& .& .& .& 1\\
        \hline
      \end{array}
    \end{displaymath}}
  \caption{Change-of-basis matrix from $\PP_3$ to $\RR_3$.}
  \label{matrPR3}
\end{figure}

\begin{figure}[!h]
  \centering
  {
    \begin{displaymath}
      \begin{array}{c|p{0.5cm}@{\,}p{0.5cm}@{\,}p{0.5cm}@{\,}p{0.5cm}@{\,}p{0.5cm}@{\,}p{0.5cm}@{\,}p{0.5cm}@{\,}p{0.5cm}@{\,}p{0.5cm}@{\,}p{0.5cm}@{\,}p{0.5cm}@{\,}p{0.5cm}@{\,}p{0.5cm}@{\,}p{0.5cm}@{\,}p{0.5cm}@{\,}p{0.5cm}@{\,}p{0.5cm}@{}|c}
        &\pc{123}&\pc{132}&\pc{213}&\pc{231}&\pc{312}&\pc{321}&\pc{122}&\pc{212}&\pc{221}&\pc{112}&\pc{121}&\pc{211}&\pc{111}\\
        \hline
        123&1&.&.&.&.&.&.&.&.&.&.&.&.\\
        132&.&1&.&.&.&.&.&.&.&.&.&.&.\\
        213&.&.&1&.&.&.&.&.&.&.&.&.&.\\
        231&.&.&.&.&1&.&.&.&.&.&.&.&.\\
        312&.&.&.&1&.&.&.&.&.&.&.&.&.\\
        321&.&.&.&.&.&1&.&.&.&.&.&.&.\\
        122&.&.&.&.&.&.&1&1&.&.&.&.&.\\
        212&.&.&.&.&.&.&1&.&.&.&.&.&.\\
        221&.&.&.&.&.&.&.&.&.&.&1&1&.\\
        112&.&.&.&.&.&.&.&.&.&1&.&.&.\\
        121&.&.&.&.&.&.&.&.&1&.&1&.&.\\
        211&.&.&.&.&.&.&.&.&1&.&.&.&.\\
        111&.&.&.&.&.&.&.&.&.&.&.&.&1\\
        \hline
      \end{array}
    \end{displaymath}}
  \caption{Change-of-basis matrix from $\QQ_3$ to $\RR_3$.}
  \label{matrQR3}
\end{figure}

\begin{figure}[!h]
  \hspace{-1.9cm}
  \scalebox{0.4}{
    \parbox{50cm}{
      \begin{displaymath}
        \input{annexes/matrPR4}
      \end{displaymath}
    }}
  \caption{Change-of-basis matrix from $\PP_4$ to $\RR_4$.}
  \label{matrPR4}
\end{figure}

\begin{figure}[!h]
  \hspace{-1.9cm}
  \scalebox{0.4}{
    \parbox{50cm}{
      \begin{displaymath}
        \input{annexes/matrOQ4}
      \end{displaymath}
    }}
  \caption{Change-of-basis matrix from $\OO_4$ to $\QQ_4$.}
  \label{matrOQ4}
\end{figure}

\begin{figure}[!h]
  \hspace{-1.9cm}
  \scalebox{0.4}{
    \parbox{50cm}{
      \begin{displaymath}
        \input{annexes/matrQR4}
      \end{displaymath}
    }}
  \caption{Change-of-basis matrix from $\QQ_4$ to $\RR_4$.}
  \label{matrQR4}
\end{figure}


\end{document}